\documentclass[preprint,authoryear,11pt]{elsarticle}
\usepackage{amssymb,amsmath,latexsym}
\usepackage[T1]{fontenc}
\usepackage{amsmath,amsfonts,amssymb,dsfont}
\usepackage[english]{babel}
\usepackage{latexsym}
\usepackage{epsfig}
\usepackage{graphics}

\usepackage{enumerate}
\usepackage{multirow}

\usepackage{tabularx}
\usepackage[left=2.5cm,right=2cm,top=2cm,bottom=2cm]{geometry}
\usepackage{graphicx}
 \usepackage{natbib}
\usepackage{here}%*****************************************************************************

\numberwithin{propo}{section}
\newtheorem{thm}{Theorem}
\numberwithin{thm}{section}
\newtheorem{rem}{Remark}
\numberwithin{rem}{section}

\numberwithin{exemple}{section}
\newtheorem{lem}{Lemma}
\numberwithin{lem}{section}

\newtheorem{cor}{Corollary}
\numberwithin{cor}{section}

\def\zak{\null\hfill{$\Box$}\par\vspace*{0.2cm}}

\makeatletter

\@addtoreset{equation}{section}
\makeatother

\begin{document}

%--------------------------------------------------------------------------

\begin{frontmatter}

\title{Portmanteau test for a class of multivariate asymmetric power GARCH model}

\author[Yacouba]{Y. Boubacar Ma{\"\i}nassara}
\address[Yacouba]{Universit\'e Bourgogne Franche-Comt\'e, \\
Laboratoire de math\'{e}matiques de Besan\c{c}on, \\ UMR CNRS 6623, \\
16 route de Gray, \\ 25030 Besan\c{c}on, France.}
\ead{mailto:yacouba.boubacar\_mainassara@univ-fcomte.fr}

\author[Yacouba]{O. Kadmiri}
\ead{mailto:othman.kadmiri@univ-fcomte.fr}

\author[Yacouba]{B. Saussereau}
\ead{mailto:bruno.saussereau@univ-fcomte.fr}

\begin{abstract}
We establish the asymptotic behaviour of the sum of squared residuals  autocovariances and autocorrelations   for the class of multivariate power transformed asymmetric models. We then derive  a portmanteau test. We establish the asymptotic distribution of the proposed statistics. These asymptotic results are illustrated by Monte Carlo experiments. An application to a bivariate real financial data is also proposed.
\end{abstract}

\end{frontmatter}

%--------------------------------------------------------------------------

\section{Introduction}\label{Introduction}
In econometric application, the univariate generalized autoregressive conditional heteroscedasticity (GARCH) framework is very restrictive. Consequently the class of multivariate models is commonly used
in time series analysis and econometrics. It describes the possible cross-relationships between the time series and not only the properties of the individual time series (see for instance \cite{FZ2019}, \cite{L2005}).  There are many extensions of multivariate GARCH models (MGARCH) with many approaches because the specification of the GARCH model does not suggest a natural extension to the multivariate framework. See  for instance  \cite{BLR-Survey} for a survey on MGARCH models.  See also \cite{Silvennoinen2009} and \cite{BHL-VEC-2012handbook} for recent surveys on MGARCH processes. The MGARCH model with conditional constant correlation (CCC-GARCH) introduced by \cite{B-MGARCH} and extended by \cite{jeantheau}, seems to be one of the most popular models used to model multivariate financial series. 
\cite{FZ-MAPGARCH} proposed an asymmetric CCC-GARCH (CCC-AGARCH) model that includes the CCC-GARCH introduced by \cite{B-MGARCH}  and its generalization by \cite{jeantheau}.
In all this work, we use the following notation $\underline{u}^{\underline{v}} := (u_1^{v_1},\ldots, u_d^{v_d})'$ for $\underline{u},\underline{v}\in\mathbb R^d$, and $x^+ = \max(0, x)$ and $x^- = \max(0, -x)$. We consider the asymmetric power GARCH model with constant conditional correlation (CCC-APGARCH$(p,q)$ for short) proposed by  \cite{Nous-MAPGARCH} and defined by 
%.  The $d$-dimensional process $\underline{\varepsilon}_t = (\varepsilon_{1,t}, \ldots, \varepsilon_{d,t})'$ is called a CCC-APGARCH$(p,q)$ if it verifies
%%
\begin{equation}\label{MAPGARCH}
\left\{
\begin{aligned}
&\underline{\varepsilon}_t = H_{0t}^{1/2}\eta_t, \\
&H_{0t} = D_{0t} R_0 D_{0t},\qquad D_{0t} = \mbox{diag}(\sqrt{h_{1,0t}},\dots, \sqrt{h_{d,0t}}),\\
&\underline{h}_{0t}^{\underline{\delta}_0/2} = \underline{\omega}_0 + \sum\limits_{i=1}^q\left\{ A_{0i}^+(\underline{\varepsilon}_{t-i}^+)^{\underline{\delta}_0/2} + A_{0i}^-(\underline{\varepsilon}_{t-i}^-)^{\underline{\delta}_0/2}\right\} + \sum\limits_{j=1}^p B_{0j} \underline{h}_{0t-j}^{\underline{\delta}_0/2},
\end{aligned}
\right.	
\end{equation}
%The innovation process $(\eta_t)_t$ is a vector of size $m\times 1$ and satisfies the assumption:\\[2mm]
%\indent $\textbf{A0} : (\eta_t)$ is an  independent and identically distributed (iid for short)  sequence of variables on $\mathbb{R}^m$ with identity covariance matrix and $\mathbb{E}[\eta_t] = 0$.\\[2mm]
%With this assumption, the  matrix $H_{0t}$ is interpreted as the conditional variance (volatility) of $\underline{\varepsilon}_t$.
%
%%
where $\underline{\varepsilon}_t = (\varepsilon_{1,t}, \ldots, \varepsilon_{d,t})'$ is a $d$-dimensional process, $\underline{h}_{0t} = (h_{1,0t}, \ldots, h_{d,0t})'$, $h_{i,0t}$ is the  conditional variance of  $\varepsilon_{i,t}= h_{i,0t}^{1/2}\eta_{i,t}$ for $i=1,\dots,d$,
%with $x^+ = \max(0, x)$ and $x^- = \max(0, -x)$
\[ \underline{\varepsilon}_t^+ = \left(\{\varepsilon_{1,t}^+\}^2, \ldots, \{\varepsilon_{d,t}^+\}^2\right)'\qquad \underline{\varepsilon}_t^- = \left(\{\varepsilon_{1,t}^-\}^2, \ldots, \{\varepsilon_{d,t}^-\}^2\right)',\]
$\underline{\omega}_0$ and $\underline{\delta}_0$ are vectors of size $d \times 1$ with strictly positive coefficients, $A_{0i}^+, A_{0i}^-$ and $B_{0j}$ are matrices of size $d\times d$ with positive coefficients and $R_0$ is a correlation matrix and where the  innovation process $(\eta_t)_t$  is an  independent and identically distributed (iid for short)  sequence of variables on $\mathbb{R}^d$ with identity covariance matrix and $\mathbb{E}[\eta_t] = 0$. The parameters of the model are the coefficients of the vectors $\underline{\omega}_0$, $\underline{\delta}_0$, the coefficients of the matrices $A_{0i}^+, A_{0i}^-, B_{0j}$ and the coefficients in the lower triangular part excluding the diagonal of the matrix $R_0$.

Model \eqref{MAPGARCH} includes various GARCH class models: for $\underline{\delta}_0=(2,\dots,2)'$, we obtain the CCC-AGARCH of \cite{FZ-MAPGARCH};  when $d=1$ and $\underline{\delta}_0=2$, we retrieve the  threshold GARCH (TGARCH) of \cite{RZ-TGARCH}. 
The asymptotic properties of the quasi-maximum likelihood (QML) estimation of the model \eqref{MAPGARCH} are established by \cite{Nous-MAPGARCH} when the power $\underline{\delta}_0$ is assumed to be known
or unknown.

As mentioned by  \cite{FZ-MAPGARCH}, the attractiveness of the CCC-AGARCH models follows from their tractability. They mention three main reasons :  the number of unknown coefficients is less than in other specifications and remains tractable in small dimension; the coefficients are easy to interpret; the conditions ensuring positive definiteness of the conditional variance are simple and explicit. There is also an advantage concerning the strict stationarity conditions which are explicit too. 
In our work, we pass from a constant power CCC-AGARCH to a component-varying power CCC-AGARCH model. In addition to the theoretical contribution, the numerical illustrations proposed  in Section \ref{donnees reelles} highlight the value of this work. 
To be more specific, our study on real dataset proves that a component-varying power is relevant for the daily exchange rates of the (Dollar,Yen) against the Euro (see Table \ref{tabexchange}).

In CCC-APGARCH$(p,q)$ models, the choice of $p$ and $q$ is particularly important because the number of parameters quickly increases with $p$ and $q$, which entails statistical difficulties.
After identification and estimation of the MGARCH processes, the next important step in the modeling consists in checking if the estimated model fits the data satisfactorily. This adequacy checking step allows to validate or invalidate the choice of the orders $p$ and $q$. Thus it is important to check the validity of a MGARCH$(p,q)$ model, for given orders $p$ and $q$. 
Based on the residuals empirical autocorrelations,  \cite{BP} derived a goodness-of-fit test, the portmanteau test, for univariate
strong autoregressive moving-average (ARMA) models (i.e. under the assumption that the error term is iid).  \cite{LB} proposed a modified portmanteau test which is nowadays one
of the most popular diagnostic checking tools in ARMA modeling of time series.
Since the articles by \cite{LB} and \cite{McL}, portmanteau tests have been important tools in time series analysis, in particular for testing the adequacy of an estimated ARMA$(p,q)$ model. See also \cite{Li-Book}, for a reference book on the portmanteau tests. The standard portmanteau tests consist in rejecting the adequacy of the model for large values of some quadratic form of the residuals autocorrelations. These tests  cannot be applied directly to conditional heteroscedasticity  or other processes displaying a second order dependence. Indeed such non-linearities may arise for instance when the observed process follows a GARCH representation. 
Consequently \cite{LM} and \cite{LL} proposed a portmanteau test based on the autocorrelations of the squared residuals. 
The intuition behind this portmanteau test is that when the model is correctly specified, the autocorrelations for squared residuals will be close to zero. Other situations where the standard  tests  do not give satisfactory results can also be found for instance in  \cite{RelvasP2016}, \cite{CaoLZ2010},  \cite{frz}, \cite{BMIA2020}, \cite{BMS2018}, \cite{yac}.

%The intuition behind these portmanteau tests is that if a given
%time series model with iid innovation $\eta_t$ is appropriate for the data at hand,
%the autocorrelations of the residuals $\hat{\eta}_t$ should be close to zero, which is the theoretical
%value of the autocorrelations of  $\eta_t$. The standard portmanteau tests thus consists in rejecting
%the adequacy of the model for large values of some quadratic form of the residual
%autocorrelations.

The asymptotic theory on MGARCH model diagnostic checking is mainly limited to the univariate framework. 
As above-mentioned, \cite{LM} and \cite{LL} studied a portmanteau test based on the autocorrelations of the squared residuals.  
\cite{BHK} developed an asymptotic theory of portmanteau tests in the standard GARCH framework,
\cite{LKN} suggested a consistent specification test for GARCH$(1,1)$ model. Recently, \cite{G2020} proposed goodness-of-fit tests for certain parametrizations of conditionally heteroscedastic time series with unobserved components. \cite{FWZ-PortTest} proposed a portmanteau test for the Log-GARCH model and the exponential GARCH (EGARCH) model. For the univariate APGARCH model, a portmanteau test based on the autocovariances of the squared residuals is developed by \cite{CF-PortTest} for the APGARCH model when the power $\underline{\delta}_0$ is known and by \cite{Nous-PortTestUni} when the power $\underline{\delta}_0$ is unknown and is jointly estimated with the other parameters. See also \cite{BF2020} who recently extended the work of \cite{CF-PortTest} (when $\underline{\delta}_0$ is known and when some parameters lie on the boundary) to the class of APGARCH with exogenous covariates (APGARCH-X). In the multivariate analysis, there are a few works.  \cite{LL1997}  proposed portmanteau  statistic for multivariate conditional heteroscedasticity models (see also \cite{DL2003} and  \cite{DL2003cor}).  
%However, \cite{TT1999}
%noted that there is a loss of information in the transformation of residuals, that may induce a
%severe loss of power. 
%Other situations where the standard  tests are not robust can be found for instance in   \cite{BMIA2020,BMS2018,yac}.
\cite{D2004} (see also \cite{D2004cor}) introduced the test which is a direct generalization of the portmanteau test of \cite{LM} to the VEC-GARCH model. \cite{WTR2013} extend Duchesne's approach to the case of multivariate GARCH models with Student$-t$ innovations. Recently, \cite{K2021}  provide a residual-based approach to examine the adequacy of multivariate GARCH models. Other tests for multivariate ARCH models include those developed can be found for instance in \cite{KN1998}, \cite{TT1999} and \cite{W2002}.

Contrary to the univariate APGARCH models, there are no validation tests for the class of the model \eqref{MAPGARCH}. In this paper we generalize the results obtained by \cite{CF-PortTest},  \cite{Nous-PortTestUni} and \cite{LL1997} to the CCC-APGARCH$(p,q)$ models defined in \eqref{MAPGARCH}. This extension raises difficult problems. First, non trivial constraints on the parameters must be imposed for identifiability of the parameters (see \cite{FZ2019}). Secondly, the implementation of standard estimation methods (for instance the Gaussian quasi-maximum likelihood estimation) is not obvious because this requires a constrained high-dimensional optimization (see also \cite{L2005}).
These technical difficulties certainly explain why univariate GARCH models are much more used than MGARCH in applied works.
%This paper is devoted to the problem of the validation step of the CCC-APGARCH$(p,q)$ representation \eqref{MAPGARCH} processes, when the powers $\underline{\delta}_0$ are known and when they are unknown and are jointly estimated with other paramaters. 

The paper is organized as follows. In Section \ref{estimation} we recall the results on the quasi-maximum likelihood estimator (QMLE) and its asymptotic distribution obtained by \cite{Nous-MAPGARCH}  when the power $\underline{\delta}_0$ is known or unknown. 
Section \ref{Port-test} presents our main results which give the asymptotic theory of the sum of squared residuals  autocovariances and autocorrelations  for the wide class of CCC-APGARCH models \eqref{MAPGARCH} when the power $\underline{\delta}_0$ is known (Section \ref{Port-deltaconnu}) and when the power $\underline{\delta}_0$ is unknown and estimated (Section \ref{PortTest-inconnu}). In Section \ref{simulations} we test the null hypothesis of the CCC-APGARCH model for different values of $\underline{\delta}_0$ in both cases. Section \ref{donnees reelles} illustrates the proposed tests for CCC-APGARCH models applied to a bivariate exchange rates. 
%To establish these results, we use a recent work of \cite{Nous-MAPGARCH} and we use the asymptotic properties obtained for the CCC-APGARCH model \eqref{MAPGARCH}.

\section{Quasi-maximum likelihood estimation}\label{estimation}
When the power $\underline{\delta}_0'=(\delta_{0,1},...,\delta_{0,d})$ is known, we write 
\begin{equation*}
\theta := (\underline{\omega}',  {\alpha_{1}^{+}} ', \ldots, {\alpha_q^+}', {\alpha_1^-}', \ldots, {\alpha_q^-}', \beta'_{1}, \ldots, \beta'_{p},\rho')',
\end{equation*}
where $\alpha_i^+$ and $\alpha_i^-$ are defined by $\alpha_i^\pm = \mbox{vec}(A_i^\pm)$ for $i=1, \ldots, q$, $\beta_j = \mbox{vec}(B_j)$ for $j = 1,\ldots, p$,  and $\rho = (\rho_{21},\ldots \rho_{d1}, \rho_{32}, \ldots, \rho_{d2}, \ldots, \rho_{dd-1})'$ such that the $\rho_{ij}$'s are the components of the matrix $R$.  
The  parameter $\theta$ belongs to the parameter space
\begin{equation*}
\Theta \subset ]0,+\infty[^{d} \times [0,\infty[^{d^2(2q+p)} \times ]-1,1[^{d(d-1)/2}.
\end{equation*}
The unknown true parameter value is denoted by
\begin{equation*}
\theta_0 := (\underline{\omega}_0', {\alpha_{01}^+}', \ldots, {\alpha_{0q}^+}', {\alpha_{01}^-}', \ldots, {\alpha_{0q}^-}', {\beta_{01}}',\ldots, {\beta_{0p}}',\rho_0')'.
\end{equation*}
Similarly when the power $\underline{\delta}=(\delta_1,...,\delta_d)'$ is unknown and is jointly estimated  with the  parameter $\theta$ we denote by $\vartheta := (\theta', \underline{\delta}')'$. The parameter $\vartheta$ belongs to the parameter space
\begin{equation*}
\Delta \subset ]0,+\infty[^{d} \times [0,\infty[^{d^2(2q+p)} \times ]-1,1[^{d(d-1)/2}\times]0,+\infty[^{d}.
\end{equation*}  
The unknown true parameter value is denoted by $\vartheta_0 := (\theta_0', \underline{\delta}_0')'$,  where $\underline{\delta}_0=(\delta_{0,1},...,\delta_{0,d})'$.
%
%\begin{equation*}
%\vartheta_0 := (\theta_0', \underline{\delta}_0')', \text{ where }\underline{\delta}_0=(\delta_{0,1},...,\delta_{0,d})'.
%\end{equation*}
%%where $\underline{\delta}_0=(\delta_{0,1},...,\delta_{0,d})'$.
%%
\subsection{Estimation when the power $\underline{\delta}_0$ is known}\label{QMLEdeltaconnu}
The goal is to estimate the $s_0= d + d^2(p+2q) + d(d-1)/{2}$ coefficients of the model \eqref{MAPGARCH}. For all $\theta\in\Theta$ we let $H_t=H_t(\theta)$. 
We assume that $H_t$ is a strictly stationary and non anticipative solution of
\begin{equation}\label{Ht-connu}
\left\{\begin{aligned}
&H_t = D_t R D_t, \qquad D_t = \text{diag}\left(\sqrt{h_{1,t}}, \ldots, \sqrt{h_{d,t}}\right),\qquad R=R(\theta),\\
&\underline{h}_t^{\underline{\delta}_0/2} := \underline{h}_t^{\underline{\delta}_0/2}(\theta) = \underline{\omega} + \sum\limits_{i=1}^q \left\{A_i^+(\underline{\varepsilon}_{t-i}^+)^{\underline{\delta}_0/2} + A_i^-(\underline{\varepsilon}_{t-i}^-)^{\underline{\delta}_0/2}\right\} + \sum\limits_{j=1}^p B_j \underline{h}_{t-j}^{\underline{\delta}_0/2}.
\end{aligned}\right.
\end{equation}
Given a realization $(\underline{\varepsilon}_1, \ldots, \underline{\varepsilon}_n)$ of length $n$ satisfying the representation \eqref{MAPGARCH}, the variable $H_t$ can  be approximated for $t\geq 1$ by $\tilde{H}_t$ defined recursively by
\begin{equation*}
\left\{\begin{aligned}
&\tilde{H}_t = \tilde{D}_t R \tilde{D}_t, \qquad \tilde{D}_t = \text{diag}\left(\sqrt{\tilde{h}_{1,t}}, \ldots, \sqrt{\tilde{h}_{d,t}}\right)\\
&\underline{\tilde{h}}_t^{\underline{\delta}_0/2} := \underline{\tilde{h}}_t^{\underline{\delta}_0/2}(\theta) = \underline{\omega} + \sum\limits_{i=1}^q \left\{A_i^+(\underline{\varepsilon}_{t-i}^+)^{\underline{\delta}_0/2} + A_i^-(\underline{\varepsilon}_{t-i}^-)^{\underline{\delta}_0/2}\right\} + \sum\limits_{j=1}^p B_j \underline{\tilde{h}}_{t-j}^{\underline{\delta}/2},
\end{aligned}\right.
\end{equation*}
conditional to the initial values $\underline{\varepsilon}_0, \ldots, \underline{\varepsilon}_{1-q}, \underline{\tilde{h}}_0,\ldots, \underline{\tilde{h}}_{1-p}$.
The quasi-maximum likelihood (QML) method is particularly efficient for the MGARCH class models because it provides consistent and asymptotically normal estimator for strictly stationary MGARCH processes under mild regularity conditions (but with no moment assumptions on the observed process).
The quasi-maximum likelihood estimator (QMLE) of model \eqref{MAPGARCH} is obtained by the standard estimation procedure for MGARCH class models. Thus the QMLE of $\theta_0$ of model \eqref{MAPGARCH} is defined as any measurable solution $\hat{\theta}_n$ of 
\begin{equation}\label{QMLEconnu}
\hat{\theta}_n = \underset{\theta \in \Theta}{\arg\min} \dfrac{1}{n} \sum\limits_{t=1}^n \tilde{l}_t(\theta), \qquad \tilde{l}_t(\theta) = \underline{\varepsilon}_t' \tilde{H}_t^{-1} \underline{\varepsilon}_t + \log(\det(\tilde{H}_t)).
\end{equation}
To ensure the asymptotic properties of the QMLE for model \eqref{MAPGARCH} obtained by \cite{Nous-MAPGARCH}, we need the following assumptions:
%\indent \textbf{A1:}

\indent \textbf{A1}: $\theta_0 \in \Theta$ and $\Theta$ is compact.

Now, we rewrite the  first equation of \eqref{MAPGARCH} as
\begin{equation*}\label{eta-tilde}
 \underline{\varepsilon}_t = D_t\overline{\eta}_t,\qquad \mbox{ where } \overline{\eta}_t = (\overline{\eta}_{1,t},\ldots, \overline{\eta}_{d,t})' = R_0^{1/2}\eta_t.
\end{equation*}
Using the third equation of model \eqref{MAPGARCH}, we may write
\begin{equation*}\label{upsilon}
(\underline{\varepsilon}_t^\pm)^{\underline{\delta}_0/2} = (\Upsilon_t^{\pm,(\underline{\delta}_0)})\underline{h}_{0t}^{\underline{\delta}_0/2}, \mbox{ with } \Upsilon_t^{\pm,(\underline{\delta}_0)} = \mbox{diag}\left((\pm\overline{\eta}_{1,t}^\pm)^{\delta_{0,1}},\ldots, (\pm\overline{\eta}_{d,t}^\pm)^{\delta_{0,d}}\right).
\end{equation*}
	\indent\textbf{A2}: $\forall \theta \in \Theta, \det(\mathcal{B}_0(z)) = 0 \Rightarrow \vert z\vert > 1$ for $\mathcal{B}_0(z) = I_d - \sum_{j=1}^p B_{0j}z^j$ and $\gamma(\textbf{C}_0) < 0$, where $\gamma(\cdot)$ is the top Lyapunov exponent of the sequence of matrix $C_0 = \{C_{0t}, t \in \mathbb{Z}\}$ with the matrix $C_{0t}$ of size $(p+2q)d \times (p+2q)d$ been defined by %in the appendix (see \eqref{C0t}).\\
%	The probabilistic properties of the model \eqref{MAPGARCH} rely on the sequence of matrices $(C_{0t})$ defined by
\begin{equation*}\label{C0t}
C_{0t} = \left(
\begin{tabular}{ccc}
$\Upsilon_t^{+,(\underline{\delta}_0)}A_{01 : q}^+$ & $\Upsilon_t^{+,(\underline{\delta}_0)}A_{01 : q}^-$ & $\Upsilon_t^{+,(\underline{\delta}_0)}B_{01:p}$\\
$I_{d(q-1)}$ & \multicolumn{2}{c}{$0_{d(q-1) \times d(p+q+1)}$}\\
$\Upsilon_t^{-,(\underline{\delta}_0)}A_{01 : q}^+$ & $\Upsilon_t^{-,(\underline{\delta}_0)}A_{01 : q}^-$ & $\Upsilon_t^{-,(\underline{\delta}_0)}B_{01:p}$\\
$0_{d(q-1)\times dq}$ & $I_{d(q-1)}$ & $0_{d(q-1) \times d(p+1)}$ \\
$A_{01:q}^+$ & $A_{01:q}^-$ & $B_{01:p}$\\
$0_{d(p-1)\times 2dq}$ & $I_{d(p-1)}$ & $0_{d(p-1)\times d}$  \\
\end{tabular}
\right),
\end{equation*}
where the $d\times qd$ matrices $A_{01:q}^+ = (A_{01}^+\ldots A_{0q}^+)$, $A_{01:q}^- = (A_{01}^-\ldots A_{0q}^-)$ and $B_{01:q} = (B_{01}\ldots B_{0q})\in\mathbb{R}^{d\times pd}$. 
%The matrix $C_{0t}$ is a $(p+2q)d \times (p+2q)d$ matrix.

	\indent\textbf{A3}: For $i=1,\ldots, d$ the distribution of $\overline{\eta}_{it}$ is not concentrated on $2$ points and $\mathbb{P}(\overline{\eta}_{it}>0)\in  (0,1)$.
	
%	, where $ \overline{\eta}_t = (\overline{\eta}_{1,t},\ldots, \overline{\eta}_{d,t})' = R_0^{1/2}\eta_t$.\\
	\indent\textbf{A4}: For $\mathcal{A}_0^+(z)=\sum_{i=1}^q A_{0i}^+z^{i}$ and $\mathcal{A}_0^-(z)=\sum_{i=1}^q A_{0i}^-z^{i}$ if $p>0, \mathcal{A}_0^+(1) + \mathcal{A}_0^-(1) \neq 0$,  $\mathcal{A}_0^+(z), \mathcal{A}_0^-(z)$ and $\mathcal{B}_0(z)$ are left-coprime and the matrix 
	%$M(\mathcal{A}_0^+, \mathcal{A}_0^-, \mathcal{B}_0)$ has full rank $d$
	\[M(\mathcal{A}_0^+, \mathcal{A}_0^-, \mathcal{B}_0) = \left[a_{q_1^+}^+(1) \ldots a_{q_d^+}^+(d) a_{q_1^-}^-(1) \ldots a_{q_d^-}^-(d) b_{p_1}(1) \ldots b_{p_d}(d)\right]\]
has full rank $d$, with $q_i^+ = q_i^+(\theta_0), q_i^- = q_i^-(\theta_0)$ and $p_i = p_i(\theta_0)$ for any value of $i$, where 
$q_i^+(\theta_0)$, $q_i^-(\theta_0)$, and $p_i(\theta_0)$ denote the maximal degrees for any column $i$ of the matrix operators $\mathcal{A}_0^+$, $\mathcal{A}_0^-$ and $\mathcal{B}_0$. We also denote by $a_{q_i^+}^+(i)$ the column vector of the coefficients $L^{q_i^+}$ , by $a_{q_i^-}^-(i)$ the column vector of the coefficients $L^{q_i^-}$ in the column $i$ of $\mathcal{A}_0^+$, respectively $\mathcal{A}_0^-$ and by $b_{p_i}(i)$ the column vector of the coefficients $L^{p_i}$ in the column $i$ of $\mathcal{B}_0$.\\
	\indent\textbf{A5}: $R$ is a positive-definite correlation matrix for all $\theta \in \Theta$.\\
	\indent\textbf{A6}: $\theta_0\in \stackrel{\circ}{\Theta}$, where $\stackrel{\circ}{\Theta}$ is the interior of $\Theta$.\\
	\indent\textbf{A7}: $\mathbb{E}\Vert \eta_t\eta_t'\Vert^2 < \infty$.  \\
	
%For the strong consistency of the QMLE, we need to assume the compactness of $\Theta$ (i.e. \textbf{A1}). The assumption \textbf{A2} makes reference to the condition of strict stationarity of the model \eqref{MAPGARCH}. Identifiability assumption is also required. 
Then under Assumptions \textbf{A1}--\textbf{A7}, \cite{Nous-MAPGARCH} showed that $\hat{\theta}_n \to \theta_0$ a.s. when $n$ goes to infinity and $\sqrt{n}(\hat{\theta}_n - \theta_0)$ is asymptotically normal with mean $0$ and covariance matrix $\Omega:=J^{-1}IJ^{-1}$, where $J$ is a positive-definite matrix and $I$ is a positive semi-definite matrix, defined by
\begin{equation*}
I :=I(\theta_0)= \mathbb{E} \left[ \dfrac{\partial l_t(\theta_0)}{\partial \theta} \dfrac{\partial l_t(\theta_0)}{\partial \theta'}\right], \quad 
J:=J(\theta_0) = \mathbb{E} \left[ \dfrac{\partial^2l_t(\theta_0)}{\partial\theta\partial\theta'}\right]\quad 
\text{with } l_t(\theta) = \underline{\varepsilon}_t'H_t^{-1}\underline{\varepsilon}_t + \log(\det(H_t)).
\end{equation*}

\subsection{Estimation when the power $\underline{\delta}_0$ is unknown}\label{QMLEdeltainconnu}
Similar to the previous section we have $s_0= 2d + d^2(p+2q) + d(d-1)/{2}$ coefficients of model \eqref{MAPGARCH} to estimate. In order to ensure that parameter $\underline{\delta}_0$ is identified we need the following additional assumption: \\[2mm]
\indent \textbf{A8}: $\eta_t$ has a positive density on some neighbourhood of zero.\\[2mm]
For all $\vartheta\in\Delta$ we let $\mathcal{H}_t=\mathcal{H}_t(\vartheta)$. 
We assume that $\mathcal{H}_t$ is a strictly stationary and non anticipative solution of
\begin{equation}\label{Ht-inconnu}
\left\{\begin{aligned}
&\mathcal{H}_t = D_t R D_t, \qquad D_t = \text{diag}\left(\sqrt{h_{1,t}}, \ldots, \sqrt{h_{d,t}}\right)\\
&\underline{h}_t := \underline{h}_t(\vartheta) = \left(\underline{\omega} + \sum\limits_{i=1}^q A_i^+(\underline{\varepsilon}_{t-i}^+)^{\underline{\delta}/2} + A_i^-(\underline{\varepsilon}_{t-i}^-)^{\underline{\delta}/2} + \sum\limits_{j=1}^p B_j \underline{h}_{t-j}^{\underline{\delta}/2}\right)^{2/\underline{\delta}}.
\end{aligned}\right.
\end{equation}
Conditionally to the initial values $\underline{\varepsilon}_0, \ldots, \underline{\varepsilon}_{1-q}, \underline{\tilde{h}}_0,\ldots, \underline{\tilde{h}}_{1-p}$, for $t\geq 1$ the variable $\mathcal{H}_t$ can also be approximated recursively by
\begin{equation*}
\left\{\begin{aligned}
&\tilde{\mathcal{H}}_t = \tilde{D}_t R \tilde{D}_t, \qquad \tilde{D}_t = \text{diag}\left(\sqrt{\tilde{h}_{1,t}}, \ldots, \sqrt{\tilde{h}_{d,t}}\right)\\
&\underline{\tilde{h}}_t := \underline{\tilde{h}}_t(\vartheta) = \left(\underline{\omega} + \sum\limits_{i=1}^q A_i^+(\underline{\varepsilon}_{t-i}^+)^{\underline{\delta}/2} + A_i^-(\underline{\varepsilon}_{t-i}^-)^{\underline{\delta}/2} + \sum\limits_{j=1}^p B_j \underline{\tilde{h}}_{t-j}^{\underline{\delta}/2}\right)^{2/\underline{\delta}}.
\end{aligned}\right.
\end{equation*}
The QMLE of $\vartheta_0$ is defined as any measurable solution $\hat{\vartheta}_n$ of 
\begin{equation}\label{QMLEinconnu}
\hat{\vartheta}_n = \underset{\vartheta \in \Delta}{\arg\min} \dfrac{1}{n} \sum\limits_{t=1}^n \tilde{\ell}_t(\vartheta), \qquad \tilde{\ell}_t(\vartheta) = \underline{\varepsilon}_t' \tilde{\mathcal{H}}_t^{-1} \underline{\varepsilon}_t + \log(\det(\tilde{\mathcal{H}}_t)).
\end{equation}
To ensure the asymptotic properties of the QMLE of $\vartheta_0$ for model \eqref{MAPGARCH} obtained by \cite{Nous-MAPGARCH}, we need assumptions similar to those we assumed in the case when the power $\underline{\delta}_0$ is known. We will assume Assumptions \textbf{A1}--\textbf{A6} with parameter $\theta$ replaced by $\vartheta$ and the space parameter $\Theta$ replaced by $\Delta$.

Under Assumptions \textbf{A1}--\textbf{A8} \cite{Nous-MAPGARCH} showed that $\hat{\vartheta}_n \to \vartheta_0$ a.s. when $n$ goes to infinity and $\sqrt{n}(\hat{\vartheta}_n - \vartheta_0)$ is asymptotically normal with mean $0$ and covariance matrix $\Omega:={\cal J}^{-1}{\cal I}{\cal J}^{-1}$, where ${\cal J}$ is a positive-definite matrix and ${\cal I}$ is a positive semi-definite matrix, defined by
\begin{equation*}
{\cal I} :={\cal I}(\vartheta_0)= \mathbb{E} \left[ \dfrac{\partial \ell_t(\vartheta_0)}{\partial \vartheta} \dfrac{\partial \ell_t(\vartheta_0)}{\partial \vartheta'}\right], \quad 
{\cal J}:={\cal J}(\vartheta_0) = \mathbb{E} \left[ \dfrac{\partial^2 \ell_t(\vartheta_0)}{\partial\vartheta\partial\vartheta'}\right]\quad 
\text{with } \ell_t(\vartheta) = \underline{\varepsilon}_t'\mathcal{H}_t^{-1}\underline{\varepsilon}_t + \log(\det(\mathcal{H}_t)).
\end{equation*}
In all the sequel we denote  by $\xrightarrow[]{\mathrm{d}}$
%, respectively by $\xrightarrow[]{\mathbb{P}}$
the convergence in distribution.
%, respectively in probability.
The symbol $\mathrm{o}_{\mathbb P}(1)$ is used for a sequence of random variables that converge to zero in probability.
\section{Diagnostic checking with portmanteau tests}\label{Port-test}
To check the adequacy of a given multivariate time series model, for instance for an estimated VARMA$(p,q)$ model, it is common practice to test the significance of the multivariate residuals autocovariances. In the MGARCH framework this approach is not relevant because the process $\eta_t$
%$\eta_t = H_t^{-1/2}\underline{\varepsilon}_t$ 
is always a white noise (possibly a
martingale difference) even when the volatility is misspecified.
For this reason the following portmanteau test is based on the squared residuals autocovariances. The null hypothesis is 
\begin{equation*}
\boldsymbol{\mathcal{H}}_0 : \text{the process $(\underline{\varepsilon}_t)$ satisfies model \eqref{MAPGARCH}}.
\end{equation*}
%
%Define the autocovariances of the sum of squared residuals at lag $h$, for $h < n$, by
\subsection{Portmanteau test when the power $\underline{\delta}_0$ is known}\label{Port-deltaconnu}
Let $\hat{\eta}_{t}=\tilde{\eta}_{t}(\hat\theta_n)=\tilde{H}_t^{-1/2}(\hat\theta_n)\underline{\varepsilon}_t=\hat{H}_t^{-1/2}\underline{\varepsilon}_t$ be the QMLE residuals when $p+q>0$ and where $\tilde{\eta}_t(\theta) = \tilde{H}_t^{-1/2}({\theta})\underline{\varepsilon}_t$.

We define the autocovariances  of the sum of squared residuals at lag $h>0$, for $h < n$, by
\begin{eqnarray*}
{\hat{r}}_h ={\tilde{r}}_h (\hat{\theta}_n)\qquad\text{where }\quad\tilde{r}_h(\theta)& = &\dfrac{1}{n} \sum\limits_{t = h +1}^n [\text{Tr}(\tilde{s}_t(\theta))][\text{Tr}(\tilde{s}_{t-h}(\theta))]\qquad \text{with } \tilde{s}_t(\theta) = \tilde{\eta}_t(\theta)\tilde{\eta}_t'(\theta) - I_d\\&=& \dfrac{1}{n} \sum\limits_{t = h +1}^n [\tilde{\eta}_t'(\theta)\tilde{\eta}_t(\theta)-d][\tilde{\eta}_{t-h}'(\theta)\tilde{\eta}_{t-h}(\theta)-d]
\\&=& \dfrac{1}{n} \sum\limits_{t = h +1}^n [\underline{\varepsilon}_t'\tilde{H}_t^{-1}({\theta})\underline{\varepsilon}_t-d][\underline{\varepsilon}_{t-h}'\tilde{H}_{t-h}^{-1}({\theta})\underline{\varepsilon}_{t-h}-d].
\end{eqnarray*}
Similarly we define the "empirical" autocovariances  of the sum of squared white noise at lag $h$ by
\begin{eqnarray*}
{{r}}_h ={{r}}_h ({\theta_0})\qquad\text{where }\quad {r}_h(\theta)& = &\dfrac{1}{n} \sum\limits_{t = h +1}^n [\text{Tr}({s}_t(\theta))][\text{Tr}({s}_{t-h}(\theta))],
%\qquad \text{with } {s}_t(\theta) = {\eta}_t(\theta){\eta}_t'(\theta) - I_d.
%\\&=& \dfrac{1}{n} \sum\limits_{t = h +1}^n [{\eta}_t'(\theta){\eta}_t(\theta)-d][{\eta}_{t-h}'(\theta){\eta}_{t-h}(\theta)-d]
%\\&=& \dfrac{1}{n} \sum\limits_{t = h +1}^n [\underline{\varepsilon}_t'{H}_t^{-1}({\theta})\underline{\varepsilon}_t-d][\underline{\varepsilon}_{t-h}'{H}_{t-h}^{-1}({\theta})\underline{\varepsilon}_{t-h}-d].
\end{eqnarray*}
with ${s}_t(\theta) = {\eta}_t(\theta){\eta}_t'(\theta) - I_d$ and ${\eta}_{t}(\theta)={{ H}}_t^{-1/2}(\theta)\underline{\varepsilon}_t$.
It should be noted that ${{r}}_h$ is not a
statistic (unless if $p=q=0$) because it depends on the unobserved
innovations $\eta_t$.

For a fixed integer $m\geq 1$ and in the sequel we will also need these following vectors: 
\begin{equation*}
\mathbf{\hat{r}}_m = \left(\hat{r}_1, \ldots, \hat{r}_m\right)'\quad\text{ and }\quad
{\mathbf{{r}}}_m = \left({r}_1, \ldots, {r}_m\right)',\  \text{such that}\  1\leq m\leq n.
\end{equation*}
To ensure the invertibility of the asymptotic covariance matrix of the vector of the sum of squared residuals autocovariances we need the following technical assumption on the distribution of  $\eta_t$.

\indent \textbf{A9}: For $d\geq2$, $\eta_t$ takes more than $3(d+1)$ positive values and more than $3(d+1)$ negative values.

Let $\mathbb{S}_{t-1:t-m} = (S_{t-1},\ldots, S_{t-m})'$, where $S_t=\eta_t'\eta_t-d$. The following theorem gives the asymptotic distribution  of the vector of the sum of squared residuals autocovariances. 
\begin{thm}\label{thmD-connu}
Under Assumptions \textbf{A1}--\textbf{A7} and \textbf{A9}, if $\underline{\varepsilon}_t$ is the non-anticipative and stationary solution of the CCC-APGARCH$(p,q)$ model defined in \eqref{MAPGARCH}, then  we have
\begin{equation*}
\sqrt{n}{\mathbf{\hat{r}}}_m \xrightarrow[n\to\infty]{\mathrm{d}}\mathcal{N}(0, D),\quad \text{where}\quad D = \left(\mathbb{E}\left[S_{t}^2\right]\right)^2I_m + C_mJ^{-1}IJ^{-1}C_m' + C_m\Sigma_{\hat{\theta}_n,{\mathbf{r}}_m}  + \Sigma_{\hat{\theta}_n,{\mathbf{r}}_m}'C_m'
\end{equation*}
is a non-singular matrix and where  $\Sigma_{\hat{\theta}_n,{\mathbf{r}}_m}=\mathbb{E}\left[J^{-1}\mathbf{h}'_t\text{vec}(\eta_t\eta_t'-I_d)S_t\mathbb{S}_{t-1:t-m}'\right]$ and the matrix $C_m$ is given by \eqref{Cm} in the proof of Theorem \ref{thmD-connu}.
\end{thm}
The proof of Theorem \ref{thmD-connu} is postponed to Section \ref{preuves}. 
\begin{rem}\label{remCon}
When we  assume that: $\mathbb{E}(\eta_{it}^3)=0$, for $i,=1,\ldots,d$; for $i,j\in\{1,\ldots,d\}$ and $i\neq j$,  $\eta_{it}$ and $\eta_{jt}$ are mutually uncorrelated up to the fourth order and $\eta_{it}$'s have the same fourth order moment, we have:
 $\mathbb{E}\left[S_{t}^2\right]=d\left(\mathbb{E}\left[\eta_{it}^4\right]-1\right)$ and $\Sigma_{\hat{\theta}_n,{\mathbf{r}}_m}=-\left(\mathbb{E}\left[\eta_{it}^4\right]-1\right)J^{-1}C'_m$.
Thus we obtain
$$D = d^2\left(\mathbb{E}\left[\eta_{it}^4\right]-1\right)^2I_m+C_m\left(J^{-1}IJ^{-1}  -2\left(\mathbb{E}\left[\eta_{it}^4\right]-1\right)J^{-1}\right)C'_m.$$
Therefore we retrieve the well-known result obtained by \cite{LL1997}.
\end{rem}
The above theorem is useless for practical purpose because it does not involve any observable quantities.
In order to state our second result we need to define a consistent estimator of the asymptotic matrix $D$ (see Theorem \ref{thmD-connu}). 

%We also need further notations. 
In view of \cite{Nous-MAPGARCH} the matrices $I$ and $J$ can be estimated by their empirical or observable counterparts given by
\begin{equation*}
\begin{aligned}
\hat{I}(i,j) &= \dfrac{1}{n} \sum\limits_{t=1}^n\left[  Tr\left(\left(\hat{H}_{t}^{-1}- \hat{H}_{t}^{-1} \underline{\varepsilon}_t\underline{\varepsilon}_t' \hat{H}_{t}^{-1}\right) \dfrac{\partial \hat{H}_{t}}{\partial\theta_i}\right)
 Tr\left(\left(\hat{H}_{t}^{-1}- \hat{H}_{t}^{-1} \underline{\varepsilon}_t\underline{\varepsilon}_t' \hat{H}_{t}^{-1}\right) \dfrac{\partial \hat{H}_{t}}{\partial\theta_j}\right) \right]\\
\text{and }\quad \hat{J}(i,j) &= \dfrac{1}{n} \sum\limits_{t=1}^n\left[Tr\left(\hat{H}_{t}^{-1} \dfrac{\partial \hat{H}_{t}}{\partial\theta_j}\hat{H}_{t}^{-1} \dfrac{\partial \hat{H}_{t}}{\partial\theta_i} \right)\right],%, \quad\text{ where }\quad\hat{H}_t:=\tilde{H}_{t}(\hat\theta_n),
\quad\text{ for }\quad i,j=1,\ldots,s_0 .
\end{aligned}
\end{equation*}
Let $\hat{\Sigma}_{\hat{\theta}_n,{\mathbf{r}}_m}$ and $\hat{C}_m$ be weakly consistent estimators of
$\Sigma_{\hat{\theta}_n,{\mathbf{r}}_m}$ and ${C}_m$ involved in the asymptotic normality of $\sqrt{n}{\mathbf{\hat{r}}}_m$. 
%Let $\hat{S}_t =\underline{\varepsilon}_t'\hat{H}_t^{-1}\underline{\varepsilon}_t-d$ and 
Define the  matrix $\hat{C}_m$ of size $m \times s_0$ whose $(h,i)-$th element is given by
\begin{equation*}
\hat{C}_m(h,i) = -\dfrac{1}{n}\sum\limits_{t=h+1}^n\hat{S}_{t-h}
\text{Tr}\left(\hat{H}_{t}^{-1}\dfrac{\partial \tilde{H}_t(\hat{\theta}_n)}{\partial\theta_i}\right) \text{ for }1\leq h \leq m\text{ and }1\leq i \leq s_0,
\end{equation*}
where $\hat{S}_t =\underline{\varepsilon}_t'\hat{H}_t^{-1}\underline{\varepsilon}_t-d$. 
The matrix $\Sigma_{\hat{\theta}_n,{\mathbf{r}}_m}$ can be estimated by its empirical or observable counterpart given by
\begin{equation*}
\hat{\Sigma}_{\hat{\theta}_n,{\mathbf{r}}_m} = \dfrac{1}{n} \sum\limits_{t=1}^n\hat{J}^{-1}\hat{\mathbf{h}}'_t\text{vec}(\hat{\eta}_t\hat{\eta}_t'-d)\hat{S}_t\hat{\mathbb{S}}_{t-1:t-m}'\quad\text{where}\quad \hat{\mathbf{h}}_t(i) = \text{vec}\left(\hat{H}_{t}^{-1/2}\dfrac{\partial \tilde{H}_t(\hat{\theta}_n)}{\partial\theta_i}\hat{H}_{t}^{-1/2}\right) \text{ for }1\leq i \leq s_0 .
\end{equation*}
Let $\hat{D} = \hat{\kappa}^2I_m + \hat{C}_m\hat{J}^{-1}\hat{I}\hat{J}^{-1}\hat{C}_m' + \hat{C}_m\hat{\Sigma}_{\hat{\theta}_n,{\mathbf{r}}_m}  + \hat{\Sigma}_{\hat{\theta}_n,{\mathbf{r}}_m}'\hat{C}_m'$ and $\hat{\kappa}=n^{-1}\sum_{t=1}^n\hat{S}_{t}^2$   be weakly consistent estimators of the matrix $D$ and $\mathbb{E}\left[S_{t}^2\right]$.

The above quantities are now all observable, we are able to state our second theorem which gives the asymptotic distribution for quadratic forms of the autocovariances.
% of the sum of squared residuals.
%
\begin{thm}\label{thm-Port-test-connu}
Under the assumptions of Theorem \ref{thmD-connu} and $\boldsymbol{\mathcal{H}}_0$, we have
\begin{equation*}
n{\mathbf{\hat{r}}}_m'\hat{D}^{-1}{\mathbf{\hat{r}}}_m \xrightarrow[n\to\infty]{\mathrm{d}} \chi^2_{m}.
\end{equation*}  
\end{thm}
The proof of Theorem \ref{thm-Port-test-connu} is postponed to Section \ref{preuves}. 
\begin{rem}\label{bobo}
If we focuse on the following alternative hypothesis 
\begin{equation*}
\boldsymbol{\mathcal{H}}_1 : \text{the process $(\underline{\varepsilon}_t)$ does not admit the representation  \eqref{MAPGARCH} with parameter $\theta_0$} ,
\end{equation*}
at least one $ r^0_h= \mathbb E [ (\eta_t'\eta_t-d)(\eta'_{t-h}\eta_{t-h}-d)] \neq 0$ under $\boldsymbol{\mathcal{H}}_1$. One may prove that under $\boldsymbol{\mathcal{H}}_1$
\begin{equation*}
 {\mathbf{\hat{r}}}_m'\hat{D}^{-1}{\mathbf{\hat{r}}}_m \xrightarrow[n\to\infty]{\mathbb{P}}  {\mathbf{r}^0_m}'  D^{-1}   \mathbf{r}^0_m 
\end{equation*}
where the vector $\mathbf{r}^0_m = (r^0_1,\ldots, r^0_m)'$. 
Therefore the test statistic $n {\mathbf{\hat{r}}}_m'\hat{D}^{-1}{\mathbf{\hat{r}}}_m $ is consistent in detecting $\boldsymbol{\mathcal{H}}_1$.
\end{rem}
The proof of this remark is also postponed to Section \ref{preuves}.

Consider the vector of the first $m$ autocorrelations  of the sum of squared residuals $$\hat{\rho}_m=(\hat{\rho}(1),\dots,\hat{\rho}(m))'\quad\text{where}
\quad\hat{\rho}(h)=\frac{\hat{r}_h}{\hat{r}_0}.$$
\begin{cor}\label{thm-Port-test-connubis}
Under the assumptions of Theorem \ref{thm-Port-test-connu}, we have
\begin{align}\label{rho1-connu}
\sqrt{n}\hat{\rho}_m& \xrightarrow[n\to\infty]{\mathrm{d}}\mathcal{N}(0, D_{\hat{\rho}})\quad \text{where}\quad D_{\hat{\rho}} = \frac{D}{\left(\mathbb{E}\left[\eta_{t}'\eta_t-d\right]^2\right)^2},
\\ \label{rho2-connu}
n\hat{\rho}_m'\hat{D}_{\hat{\rho}}^{-1}\hat{\rho}_m &\xrightarrow[n\to\infty]{\mathrm{d}} \chi^2_{m}\quad \text{where}\quad \hat{D}_{\hat{\rho}} = \frac{\hat{D}}{\left(\frac{1}{n}\sum_{t=1}^n\left[\hat{\eta}_{t}'\hat{\eta}_t-d\right]^2\right)^2}.
\end{align}  
\end{cor}
The proof of Corollary \ref{thm-Port-test-connubis} is postponed to Section \ref{preuves}. 
\begin{rem}\label{remConbis}
Under the assumptions of Remark \ref{remCon}, we have: 
$\mathbb{E}[\underline{\varepsilon}_t'{H}_t^{-1}\underline{\varepsilon}_t-d]^2=
\left(\mathbb{E}\left[\eta_{it}^4\right]-1\right)d$, so that 
 $D_{\hat{\rho}} = {D}/{\left(\mathbb{E}\left[\eta_{it}^4\right]-1\right)^2d^2}$
and $\hat{D}_{\hat{\rho}} = {\hat{D}}/{\left(\hat{\kappa}_i-1\right)^2d^2},
\quad\text{ for }\quad i=1,\ldots,d$, where $\hat{\kappa}_i=n^{-1}\sum_{t=1}^n\hat{\eta}_{it}^4$.
\end{rem}
The adequacy of the CCC-APGARCH$(p,q)$ model defined in \eqref{MAPGARCH} is then rejected at the asymptotic level $\alpha$ when
\begin{equation*}
n{\mathbf{\hat{r}}}_m'\hat{D}^{-1}{\mathbf{\hat{r}}}_m > \chi^2_{m}(1-\alpha)\quad\text{or}\quad n\hat{\rho}_m'\hat{D}_{\hat{\rho}}^{-1}\hat{\rho}_m > \chi^2_{m}(1-\alpha),
\end{equation*}
where $\chi^2_{m}(1-\alpha)$ represents the $(1-\alpha)$-quantile of the chi-square distribution with $m$ degrees of freedom.
\subsection{Portmanteau test when the power is unknown}\label{PortTest-inconnu}
The results are close to those given in Section \ref{Port-deltaconnu}.
It consists to adapt the notations in Section \ref{Port-deltaconnu} by replacing $\theta_0$ (resp. $\hat{\theta}_n$) by $\vartheta_0$ (resp. $\hat{\vartheta}_n$) and $H_t$ (resp. $\tilde{H}_t$) by ${\cal H}_t$ (resp. $\tilde{{\cal H}}_t$).

To establish the asymptotic distribution of the portmanteau test statistic, when $\underline{\delta}_0$ is unknown, Assumption \textbf{A9} is replacing by

\indent $\textbf{A9'}:$ For $d\geq2$, $\eta_t$ takes more than $11d+1$ positive values and more than $11d+1$ negative values.
%
%
%%
%The following theorem gives the asymptotic distribution  of the vector of the sum of squared residuals autocovariances. 
\begin{thm}\label{thmD-inconnu}
Under Assumptions \textbf{A1}--\textbf{A8} and \textbf{A9'}, if $\underline{\varepsilon}_t$ is the non-anticipative and stationary solution of the CCC-APGARCH$(p,q)$ model defined in \eqref{MAPGARCH}, then  we have
\begin{equation*}
\sqrt{n}{\mathbf{\hat{r}}}_m \xrightarrow[n\to\infty]{\mathrm{d}}\mathcal{N}(0, {\cal D}),\quad \text{where}\quad {\cal D} = \left(\mathbb{E}\left[S_{t}^2\right]\right)^2I_m + {\cal C}_m{\cal J}^{-1}{\cal I}{\cal J}^{-1}{\cal C}_m' + {\cal C}_m\Sigma_{\hat{\vartheta}_n,{\mathbf{r}}_m}  + \Sigma_{\hat{\vartheta}_n,{\mathbf{r}}_m}'{\cal C}_m'
\end{equation*}
is a non-singular matrix and where the matrix ${\cal C}_m$ is given by \eqref{CmInc} in the proof of Theorem \ref{thmD-inconnu}.
\end{thm}
The proof of Theorem \ref{thmD-inconnu} is postponed to Section \ref{preuves}.

In view of \cite{Nous-MAPGARCH} the matrices ${\cal I}$ and ${\cal J}$ can be estimated by 
\begin{equation*}
\begin{aligned}
\hat{{\cal I}}(i,j) &= \dfrac{1}{n} \sum\limits_{t=1}^n\left[  Tr\left(\left(\hat{{\cal H}}_{t}^{-1}- \hat{{\cal H}}_{t}^{-1} \underline{\varepsilon}_t\underline{\varepsilon}_t' \hat{{\cal H}}_{t}^{-1}\right) \dfrac{\partial \hat{{\cal H}}_{t}}{\partial\vartheta_i}\right)
 Tr\left(\left(\hat{{\cal H}}_{t}^{-1}- \hat{{\cal H}}_{t}^{-1} \underline{\varepsilon}_t\underline{\varepsilon}_t' \hat{{\cal H}}_{t}^{-1}\right) \dfrac{\partial \hat{{\cal H}}_{t}}{\partial\vartheta_j}\right) \right]\\
\text{and }\quad \hat{{\cal J}}(i,j) &= \dfrac{1}{n} \sum\limits_{t=1}^n\left[Tr\left(\hat{{\cal H}}_{t}^{-1} \dfrac{\partial \hat{{\cal H}}_{t}}{\partial\vartheta_j}\hat{{\cal H}}_{t}^{-1} \dfrac{\partial \hat{{\cal H}}_{t}}{\partial\vartheta_i} \right)\right], %, \quad\text{ where }\quad\hat{H}_t:=\tilde{H}_{t}(\hat\vartheta_n),\quad\text{ for }\quad i,j=1,\ldots,s_0 .
\quad\text{ for }\quad i,j=1,\ldots,s_0 .
\end{aligned}
\end{equation*}
Let $\hat{S}_t =\underline{\varepsilon}_t'\hat{{\cal H}}_t^{-1}\underline{\varepsilon}_t-d$ and define the matrix $\hat{C}_m$ of size $m \times s_0$ whose $(h,i)-$th element is given by
\begin{equation*}
\hat{{\cal C}}_m(h,i) = -\dfrac{1}{n}\sum\limits_{t=h+1}^n\hat{S}_{t-h}
\text{Tr}\left(\hat{{\cal H}}_{t}^{-1}\dfrac{\partial \tilde{{\cal H}}_t(\hat{\vartheta}_n)}{\partial\vartheta_i}\right) \text{ for }1\leq h \leq m\text{ and }1\leq i \leq s_0.
\end{equation*}
%where $\hat{{\cal C}}_m(h,i)$ denotes the $(h,i)-$th element of the 
The matrix $\hat{{\cal C}}_m$ is  a weakly consistent estimator of ${{\cal C}}_m$. The matrix $\Sigma_{\hat{\vartheta}_n,{\mathbf{r}}_m}$ can also be estimated by 
\begin{equation*}
\hat{\Sigma}_{\hat{\vartheta}_n,{\mathbf{r}}_m} = \dfrac{1}{n} \sum\limits_{t=1}^n\hat{{\cal J}}^{-1}\hat{\mathbf{h}}'_t\text{vec}(\hat{\eta}_t\hat{\eta}_t'-d)\hat{S}_t\hat{\mathbb{S}}_{t-1:t-m}'\quad\text{where}\quad \hat{\mathbf{h}}_t(i) = \text{vec}\left(\hat{{\cal H}}_{t}^{-1/2}\dfrac{\partial \tilde{{\cal H}}_t(\hat{\vartheta}_n)}{\partial\vartheta_i}\hat{{\cal H}}_{t}^{-1/2}\right) \text{ for }1\leq i \leq s_0 .
\end{equation*}
Let $\hat{{\cal D}} = \hat{\kappa}^2I_m + \hat{{\cal C}}_m\hat{{\cal J}}^{-1}\hat{{\cal I}}\hat{{\cal J}}^{-1}\hat{{\cal C}}_m' + \hat{{\cal C}}_m\hat{\Sigma}_{\hat{\vartheta}_n,{\mathbf{r}}_m}  + \hat{\Sigma}_{\hat{\vartheta}_n,{\mathbf{r}}_m}'\hat{{\cal C}}_m'$
be a weakly consistent estimator of ${\cal D}$.

We are able to state the second theorem of this section  which gives the asymptotic distribution for quadratic forms of the sum of squared residuals autocovariances.
\begin{thm}\label{thm-Port-test-inconnu}
Under the assumptions of Theorem \ref{thmD-inconnu} and $\boldsymbol{\mathcal{H}}_0$, we have
\begin{equation*}
n{\mathbf{\hat{r}}}_m'\hat{{\cal D}}^{-1}{\mathbf{\hat{r}}}_m \xrightarrow[n\to\infty]{\mathrm{d}} \chi^2_{m}.
\end{equation*}  
\end{thm}
The proof of Theorem \ref{thm-Port-test-inconnu} is postponed to Section \ref{preuves}. 
\begin{rem}\label{bobo2}
If we focuse on the following alternative hypothesis 
\begin{equation*}
\boldsymbol{\mathcal{H}}_1 : \text{the process $(\underline{\varepsilon}_t)$ does not admit the representation  \eqref{MAPGARCH} with parameter $\vartheta_0$} ,
\end{equation*}
at least one $ r^0_h= \mathbb E [ (\eta_t'\eta_t-d)(\eta'_{t-h}\eta_{t-h}-d)] \neq 0$ under $\boldsymbol{\mathcal{H}}_1$. One may prove that under $\boldsymbol{\mathcal{H}}_1$
\begin{equation*}
 {\mathbf{\hat{r}}}_m'\hat{{\cal D}}^{-1}{\mathbf{\hat{r}}}_m \xrightarrow[n\to\infty]{\mathbb{P}}  {\mathbf{r}^0_m}'  {\cal D}^{-1}   \mathbf{r}^0_m. 
\end{equation*}
%where the vector $\mathbf{r}^0_m = (r^0_1,\ldots, r^0_m)'$. 
Therefore the test statistic $n {\mathbf{\hat{r}}}_m'\hat{{\cal D}}^{-1}{\mathbf{\hat{r}}}_m $ is consistent in detecting $\boldsymbol{\mathcal{H}}_1$.
\end{rem}
The proof of this remark is similar to that of Remark \ref{bobo} and is omitted.

\begin{cor}\label{thm-Port-test-inconnubis}
Under the assumptions of Theorem \ref{thm-Port-test-inconnu}, we have
\begin{align}\label{rho1-inconnu}
\sqrt{n}\hat{\rho}_m& \xrightarrow[n\to\infty]{\mathrm{d}}\mathcal{N}(0,{\cal D}_{\hat{\rho}})\quad \text{where}\quad {\cal D}_{\hat{\rho}} = \frac{{\cal D}}{\left(\mathbb{E}\left[\eta_{t}'\eta_t-d\right]^2\right)^2},
\\ \label{rho2-inconnu}
n\hat{\rho}_m'\hat{{\cal D}}_{\hat{\rho}}^{-1}\hat{\rho}_m &\xrightarrow[n\to\infty]{\mathrm{d}} \chi^2_{m}\quad \text{where}\quad \hat{{\cal D}}_{\hat{\rho}} = \frac{\hat{{\cal D}}}{\left(\frac{1}{n}\sum_{t=1}^n\left[\hat{\eta}_{t}'\hat{\eta}_t-d\right]^2\right)^2}.
\end{align}  
\end{cor}
The proof of Corollary \ref{thm-Port-test-inconnubis} is postponed to Section \ref{preuves}. 

\begin{rem}\label{remInCon}
In view of Remark \ref{remConbis}, we have: 
 ${\cal D}_{\hat{\rho}} = {\cal D}/{\left(\mathbb{E}\left[\eta_{it}^4\right]-1\right)^2d^2}$
and $\hat{{\cal D}}_{\hat{\rho}} = {\hat{{\cal D}}}/{\left(\hat{\kappa}_i-1\right)^2d^2},
\quad\text{ for }\quad i=1,\ldots,d$.
%, where $\hat{\kappa}_i=n^{-1}\sum_{t=1}^n\hat{\eta}_{it}^4$.
\end{rem}

The adequacy of the CCC-APGARCH$(p,q)$ model, define in \eqref{MAPGARCH} is then rejected at the asymptotic level $\alpha$ when
\begin{equation*}
n{\mathbf{\hat{r}}}_m'\hat{{\cal D}}^{-1}{\mathbf{\hat{r}}}_m > \chi^2_{m}(1-\alpha)\quad\text{or}\quad n\hat{\rho}_m'\hat{{\cal D}}_{\hat{\rho}}^{-1}\hat{\rho}_m > \chi^2_{m}(1-\alpha).
\end{equation*}
In view of Corollary \ref{thm-Port-test-connubis} (resp. Corollary \ref{thm-Port-test-inconnubis}) 
for any $1\leq h\leq m$, a $100(1-\alpha)\%$ confidence region for $\rho(h)$ is given
by 
$$-u_\alpha\frac{\hat{{ D}}_{\hat{\rho}}(h,h)}{\sqrt{n}}\leq\hat{\rho}(h)\leq u_\alpha\frac{\hat{{ D}}_{\hat{\rho}}(h,h)}{\sqrt{n}}\quad
\left(\text{resp. } -u_\alpha\frac{\hat{{\cal D}}_{\hat{\rho}}(h,h)}{\sqrt{n}}\leq\hat{\rho}(h)\leq u_\alpha\frac{\hat{{\cal D}}_{\hat{\rho}}(h,h)}{\sqrt{n}}\right)$$
where $u_\alpha$ denotes the quantile of order $1-\alpha$ of the $\mathcal{N}(0,1)$ distribution.
\section{Numerical illustration}\label{simulations}
By means of Monte Carlo experiments we investigate the finite sample properties of the tests introduced in this paper.
The numerical illustrations are made with the free statistical software RStudio (see https://www.rstudio.com). 

We generate a bivariate CCC-APGARCH$(0,1)$ model (Model \eqref{MAPGARCH} with $p=0$ and $q=1$)
\begin{equation}\label{arch}
\left\{
\begin{aligned}
&\underline{\varepsilon}_t = H_{0t}^{1/2}\eta_t, \\
&H_{0t} = D_{0t} R_0 D_{0t},\qquad D_{0t} = \mbox{diag}(\sqrt{h_{1,0t}}, \sqrt{h_{2,0t}}),\\
&\underline{h}_{0t}^{\underline{\delta}_0/2} = \underline{\omega}_0 +  A_{01}^+(\underline{\varepsilon}_{t-1}^+)^{\underline{\delta}_0/2} + A_{01}^-(\underline{\varepsilon}_{t-1}^-)^{\underline{\delta}_0/2},
\end{aligned}
\right.	
\end{equation}
for different values of $\underline{\delta}_0\in\{(1,1),(0.8,1.5),(2,2),(3,2.5)\}$ and 
where $\underline{h}_{0t} = (h_{1,0t}, h_{2,0t})'$, 
$\underline{\varepsilon}_t^+ = \left(\{\varepsilon_{1,t}^+\}^2, \{\varepsilon_{2,t}^+\}^2\right)'$ and $ \underline{\varepsilon}_t^- = \left(\{\varepsilon_{1,t}^-\}^2, \{\varepsilon_{2,t}^-\}^2\right)'.$
The innovation process $(\eta_t)$ is
defined by
\begin{equation*}% \label{bruitfort}
\left(\begin{array}{c}\eta_{1,t}\\\eta_{2,t}\end{array}\right)\sim\mathrm{IID}\,{\cal
N}(0,I_2).
\end{equation*}
Considering other distributions for $\eta_t$ does not affect the conclusion. 

The coefficients of the data generating process (DGP for short) in \eqref{arch} are chosen such that Assumption \textbf{A2} holds. The coefficient $\underline{\omega}_0$ is a vector of size $2 \times 1$, $A_{01}^+, A_{01}^-$ and $R_0$ are matrices of size $2\times 2$ taken as:
%\begin{eqnarray}
%\label{coef}
%\underline{\omega}_0=\left(\begin{array}{c}0.2\\0.3\end{array}\right),\quad
%A_{01}^+=\left(\begin{array}{cc}0.25&0.10\\0.10&0.15\end{array}\right)
%,\quad A_{01}^-= \left(\begin{array}{cc}0.45&0.25\\0.25&0.35\end{array}\right)
%\quad\text{and}\quad
%R_{0}=\left(\begin{array}{cc}1&0.7\\0.7&1\end{array}\right).
%\end{eqnarray} 
\begin{eqnarray}
\label{coef}
\underline{\omega}_0=\left(\begin{array}{c}0.2\\0.3\end{array}\right),\;
R_{0}=\left(\begin{array}{cc}1&0.7\\0.7&1\end{array}\right)
\;\text{and}\;
\left\{
\begin{aligned}
&A_{01}^+=\left(\begin{array}{cc}0.25&0.10\\0.10&0.15\end{array}\right)
,\; A_{01}^-= \left(\begin{array}{cc}0.45&0.25\\0.25&0.35\end{array}\right), \;\text{if}\; A_{01}^-\neq A_{01}^+\\
&A_{01}^+=A_{01}^-= \left(\begin{array}{cc}0.45&0.25\\0.25&0.35\end{array}\right), \quad\text{if}\quad A_{01}^-= A_{01}^+.
\end{aligned}
\right.	
\end{eqnarray} 
We simulated $N=1,000$ independent replications of size $n=250$, $n=500$ and $n=2,000$ of Model \eqref{arch} with coefficients \eqref{coef}. 

For each of these $N$ replications of model \eqref{arch}, we use the QMLE method to estimate the  coefficient $\theta_{0}\in\mathbb R^{11}$ when the power is known (resp. $\vartheta_{0}\in\mathbb R^{13}$ when the power is unknown). After estimating the model considered we apply portmanteau test to the sum of squared residuals for different values of $m\in\{1,\dots,12\}$,  where $m$ is the number of autocorrelations used in the portmanteau test statistic.

We use in the following tables 3 nominal levels $\alpha=1\%$, 5\% and 10\%. For these nominal levels, the empirical  relative frequency of rejection size over the $N$ independent replications should vary  respectively within the confidence intervals $[0.3\%,1.7\%]$, $[3.6\%, 6.4\%]$  and $[8.1\%,11.9\%]$ with probability 95\% and  $[0.3\%, 1.9\%]$, $[3.3\%, 6.9\%]$ and $[7.6\%, 12.5\%]$ with probability 99\% under the assumption that the true probabilities of rejection are respectively $\alpha=1\%$, $\alpha=5\%$ and $\alpha=10\%$.

We repeat the same experiments to examine the empirical power of the proposed test for the null hypothesis of a bivariate CCC-APGARCH$(0,1)$ model of the form \eqref{arch} against the following bivariate CCC-APGARCH$(1,1)$ alternative  defined by
\begin{equation}\label{model-puiss}
\left\{
\begin{aligned}
&\underline{\varepsilon}_t = H_{0t}^{1/2}\eta_t, \\
&H_{0t} = D_{0t} R_0 D_{0t},\qquad D_{0t} = \mbox{diag}(\sqrt{h_{1,0t}}, \sqrt{h_{2,0t}}),\\
&\underline{h}_{0t}^{\underline{\delta}_0/2} = \underline{\omega}_0 +  A_{01}^+(\underline{\varepsilon}_{t-1}^+)^{\underline{\delta}_0/2} + A_{01}^-(\underline{\varepsilon}_{t-1}^-)^{\underline{\delta}_0/2}+ B_{01}\underline{h}_{0t-1}^{\underline{\delta}_0/2},
\end{aligned}
\right.	
\end{equation}
where the matrices $\underline{\omega}_0$, $A_{01}^+$, $A_{01}^-$ and $R_0$ are given by \eqref{coef} and 
$B_{01}= \left(\begin{array}{cc}0.43&0.1\\0.1&0.42\end{array}\right).$
%\begin{eqnarray}\label{ARCH1}
%\left(\begin{array}{c}\epsilon_{1,t}\\\epsilon_{2,t}\end{array}\right)
%&=&\left(\begin{array}{cc}h_{11,t}&0\\0&h_{22,t}\end{array}\right)
%\left(\begin{array}{c}\eta_{1,t}\\\eta_{2,t}\end{array}\right)
%\end{eqnarray}
%where 
%\begin{eqnarray*}
%\left(\begin{array}{c}h_{11,t}^{2}\\h_{22,t}^{2}\end{array}\right)=
%\left(\begin{array}{c}0.3\\0.2\end{array}\right)+
%\left(\begin{array}{cc}0.45&0.00\\0.40&0.25\end{array}\right)
%\left(\begin{array}{c}\epsilon_{1,t-1}^{2}\\
%\epsilon_{2,t-1}^{2}\end{array}\right).
%\end{eqnarray*}
\subsection{When the power is known}
Table \ref{archsym} (resp. Table \ref{archasym}) displays the empirical relative frequencies of rejection over the $N$ independent replications for the
3 nominal levels $\alpha=1\%$, 5\% and 10\%  when the DGP is the APGARCH$(0,1)$ model \eqref{arch}--\eqref{coef} with $A_{01}^+=A_{01}^-$ (resp. with $A_{01}^+\neq A_{01}^-$).

As expected, Tables \ref{archsym} and \ref{archasym} show that the percentages of rejection belong to the confidence interval with probability $95\%$ and $99\%$. Thus the type I error is better controlled.

In term of power performance, we investigate two experiments given in the following tables:

Table \ref{puiss-archasym} displays (in \%) the empirical power of the proposed test for the null hypothesis of the CCC-APGARCH$(0,1)$ model defined by \eqref{arch}--\eqref{coef} with $\underline{\delta}_0=(1,1)$ against the alternative given by \eqref{arch}--\eqref{coef} when $\underline{\delta}_0\neq (1,1)$.

Table \ref{puiss-archasymbis} displays also  (in \%) the empirical power of the proposed test for the null hypothesis of a bivariate CCC-APGARCH$(0,1)$ model of the form \eqref{arch} against the bivariate CCC-APGARCH$(1,1)$ alternative given by \eqref{model-puiss} when $\underline{\delta}_0$ is known.

We draw the conclusion that: 
\begin{itemize}
\item[$a)$]in the first experiment given in Table \ref{puiss-archasym}, the portmanteau tests are more disappointing  since they fail to detect some alternatives of the form $\underline{\delta}_0\neq (1,1)$ when the null is $\underline{\delta}_0= (1,1)$, except for $\underline{\delta}_0\geq (2,2)$  when $n$ increases.
\item[$b)$] Whereas the second experiment given in Table \ref{puiss-archasymbis} reveals that the portmanteau tests are much more powerful to detect wrong values of the order $(p,q)$ even when $n$ is small.
\end{itemize}

\subsection{When the power is unknown}
In this case, the power $\underline{\delta}_0$ is jointly estimated with the parameter $\theta_0$. As in the case where $\underline{\delta}_0$ is known,  Table \ref{archsymInc} (resp. Table \ref{archasymInc}) displays the empirical relative frequencies of rejection over the $N$ independent replications for the
3 nominal levels $\alpha=1\%$, 5\% and 10\%  when the DGP is the APGARCH$(0,1)$ model \eqref{arch}--\eqref{coef} with $A_{01}^+=A_{01}^-$ (resp. with $A_{01}^+\neq A_{01}^-$). Even in this case, Tables \ref{archsymInc} and \ref{archasymInc} show that the percentages of rejection belong to the confidence interval with probability $95\%$ and $99\%$. Thus the type I error is better controlled.
In term of power performance, Table \ref{puiss-archasymInc} shows that the powers of the test are quite satisfactory even when $n$ is small.
%the portmanteau tests are much more powerful to detect wrong values of the order $(p,q)$.

%\subsection{Remark}
%Considering other distributions for $\eta_t$ does not affect the conclusion. Other simulations have been performed with 
%\begin{itemize}
%\item a standardized Student's distribution with $\nu$ degrees of freedom ($t_{\nu}$) where $\nu=8$, such that $\mathbb{E}(\eta_{it}^2)=1$,
%\item and a centered and standardized two-components Gaussian mixture distribution ($\eta_{it}\sim 0.1\mathcal{N}(-2,2)+0.9\mathcal{N}(2,0.16)$) to obtain $\mathbb{E}(\eta_{it}^2)=1$, which is highly leptokurtic since $\kappa_{\eta_i}=10.53$.
%\begin{itemize}

\section{Adequacy of CCC-APGARCH models for real datasets}\label{donnees reelles} 
We consider the daily return of two exchange rates EUR/USD (Euros Dollar) and  EUR/JPY (Euros Yen). %EUR/GBP (Euros Pounds) and EUR/CAD (Euros Canadian dollar). 
The observations covered the period from January 4, 1999 to March 9, 2021 which correspond to $n = 5,679$ observations. The data were obtained from the website of the European Central Bank (http://www.ecb.int/stats/exchange/eurofxref/html/index.en.html). On these series, several CCC-APGARCH$(p,q)$ models of the form  \eqref{MAPGARCH} have been estimated by QML. For each estimated model, we apply the  portmanteau tests proposed in Section \ref{Port-test} to the sum of squared residuals for different values of $m\in\{1,\dots,12\}$ to test the adequacy of CCC-APGARCH models.

Table \ref{tabexchange} displays the $p-$values for adequacy of the CCC-APGARCH$(p,q)$ models for daily returns of exchange rates  based on $m$ squared residuals autocovariances, as well as the true and estimated powers (denoted $\underline{\delta}_0$ and $\hat{\underline{\delta}}$) and the likelihood. When $\underline{\delta}_0$ is known, the two corresponding models with $\underline{\delta}_0=(1,1)$ and $\underline{\delta}_0=(2,2)$ have the same number of parameters so it makes sense to prefer the model with the higher likelihood (the likelihood is given in the last column
of Table  \ref{tabexchange}). According to this criterion, the Log$-$likelihood of the preferred model is given in bold face (see Table  \ref{tabexchange}).

Table \ref{tabexchange} shows that the  CCC-APGARCH$(0,q)$ models (for $q=1,2,3$) are generally rejected whereas the CCC-APGARCH$(p,q)$ models   are not generally rejected and seem more appropriate. 
% for these series $\underline{\varepsilon}_t=(\text{USD}_t,\text{JPY}_t)'$. 
When $\underline{\delta}_0$ is known, the CCC-APGARCH$(1,1)$ and CCC-APGARCH$(2,1)$ models seem to be relevant for $\underline{\varepsilon}_t=(\text{USD}_t,\text{JPY}_t)'$. In contrast, when $\underline{\delta}_0$ is estimated, the CCC-APGARCH$(2,1)$ and CCC-APGARCH$(2,2)$ models seem to be relevant for $\underline{\varepsilon}_t$.

From the second last column of Table \ref{tabexchange}, we can also see that the estimated power $\hat\delta$ is not necessary equal to 1 or 2 and is different for each model. 

The portmanteau test is thus an important tool in the validation process.
From the empirical results and the simulation experiments, we draw the conclusion that the proposed portmanteau tests based on the sum of squared residuals of a CCC-APGARCH$(p,q)$ controls well the error of first kind at different asymptotic level $\alpha$ and is efficient to detect a misspecification of the order $(p,q)$.

\begin{table}[H]
 \caption{\small{Portmanteau test $p-$values for adequacy of the CCC-APGARCH$(p,q)$ models  for daily returns of exchange rates of the (Dollar,Yen), based on $m$ of the sum of squared residuals autocovariances.}}
{\scriptsize
\begin{center}
\begin{tabular}{p{1.9cm}c
@{\hskip0.2cm}c@{\hskip0.2cm}c@{\hskip0.2cm}c@{\hskip0.2cm}c@{\hskip0.2cm}c@{\hskip0.2cm}c@{\hskip0.2cm}c@{\hskip0.2cm}c@{\hskip0.2cm}c@{\hskip0.2cm}c@{\hskip0.2cm}c@{\hskip0.2cm}c@{\hskip0.2cm}c@{\hskip0.2cm}c}
\hline\hline
& \multicolumn{12}{c}{Lag $m$} & \multicolumn{1}{c}{\multirow{2}{*}{$\underline{\delta}_0$ or $\hat{\underline{\delta}}$}}& \multicolumn{1}{c}{\multirow{2}{*}{Log-lik}}\\
\cline{2-13}
Currency & 1 & 2 & 3 & 4 & 5 & 6 & 7 & 8 & 9 & 10 & 11 & 12 &\\
\hline
\\

\multicolumn{13}{l}{Portmanteau tests for adequacy of the CCC-APGARCH$(0,1)$ when $\underline{\delta}_0$ is known}\\
(USD,JPY) & 0.880 &0.000& 0.000& 0.000 &0.000 &0.000 &0.000 &0.000 &0.000 &0.000& 0.000& 0.000 & $(1,1)$ &{\bf -0.1295}\\
(USD,JPY) & 0.056 &0.000 &0.000 &0.000 &0.000 &0.000 &0.000 &0.000 &0.000 &0.000 &0.000 &0.000& $(2,2)$ & -0.1291\\
\\
\multicolumn{13}{l}{Portmanteau tests for adequacy of the CCC-APGARCH$(0,2)$ when $\underline{\delta}_0$ is known}\\
(USD,JPY) & 0.401& 0.520& 0.003& 0.000 &0.000 &0.000 &0.000 &0.000 &0.000 &0.000& 0.000& 0.000 & $(1,1)$ & -0.1827\\
(USD,JPY) & 0.492& 0.045& 0.001 &0.000 &0.000 &0.000 &0.000& 0.000 &0.000 &0.000& 0.000 &0.000& $(2,2)$ &{\bf -0.1844}\\
\\
\multicolumn{13}{l}{Portmanteau tests for adequacy of the CCC-APGARCH$(0,3)$ when $\underline{\delta}_0$ is known}\\
(USD,JPY) & 0.496 &0.683 &0.372 &0.000 &0.000 &0.000 &0.000 &0.000& 0.000 &0.000 &0.000 &0.000 & $(1,1)$ & -0.2002\\
(USD,JPY) & 0.600 &0.114 &0.010 &0.000 &0.000 &0.000 &0.000 &0.000 &0.000& 0.000 &0.000 &0.000& $(2,2)$ & {\bf -0.2010}\\
\\

\multicolumn{13}{l}{Portmanteau tests for adequacy of the CCC-APGARCH$(1,1)$ when $\underline{\delta}_0$ is known}\\
(USD,JPY) & 0.118& 0.214& 0.362& 0.074& 0.022& 0.027 &0.039 &0.063& 0.064& 0.039& 0.035& 0.042 & $(1,1)$ &{\bf -0.3410}\\
(USD,JPY) & 0.280& 0.479& 0.677& 0.232& 0.110& 0.146 &0.216& 0.254 &0.207& 0.159& 0.159 &0.140 & $(2,2)$&-0.3406 \\
\\

\multicolumn{13}{l}{Portmanteau tests for adequacy of the CCC-APGARCH$(2,1)$ when $\underline{\delta}_0$ is known}\\
(USD,JPY) & 0.164& 0.279& 0.322 &0.160& 0.092 &0.067 &0.072 &0.099& 0.047 &0.020& 0.009 &0.014 & $(1,1)$ &-0.2492\\
(USD,JPY) & 0.337& 0.595 &0.545& 0.260 &0.092 &0.073 &0.114 &0.145 &0.135 &0.124& 0.128& 0.169 & $(2,2)$&{\bf -0.2858} \\
\\

\multicolumn{13}{l}{Portmanteau tests for adequacy of the CCC-APGARCH$(1,2)$ when $\underline{\delta}_0$ is known}\\
(USD,JPY) & 0.402 &0.082 &0.151& 0.180& 0.081& 0.064& 0.091& 0.139 &0.086 &0.062& 0.030 &0.046 & $(1,1)$ & -0.2402\\
(USD,JPY) & 0.610 &0.082 &0.102 &0.025& 0.013 &0.004 &0.007 &0.008& 0.006 &0.009& 0.013 &0.007 & $(2,2)$ & {\bf -0.2937}\\
\\

\multicolumn{13}{l}{Portmanteau tests for adequacy of the CCC-APGARCH$(2,2)$ when $\underline{\delta}_0$ is known}\\
(USD,JPY) & 0.166 &0.191 &0.206 &0.102 &0.079 &0.060 &0.083 &0.116 &0.076 &0.045 &0.032 &0.047 & $(1,1)$ &-0.2559 \\
(USD,JPY) & 0.152 &0.205 &0.333 &0.077 &0.040& 0.020 &0.027 &0.038& 0.012 &0.012 &0.012& 0.012 & $(2,2)$ & {\bf -0.3062}\\
\\
\hline\\
\multicolumn{13}{l}{Portmanteau tests for adequacy of the CCC-APGARCH$(0,1)$ when $\underline{\delta}_0$ is unknown}\\
(USD,JPY) & 0.518 &0.000& 0.000 &0.000 &0.000 &0.000 &0.000 &0.000 &0.000 &0.000 &0.000 &0.000& $(4.595 , 1.201)$&-0.1321 \\
\\
\multicolumn{13}{l}{Portmanteau tests for adequacy of the CCC-APGARCH$(0,2)$ when $\underline{\delta}_0$ is unknown}\\
(USD,JPY) & 0.779& 0.102 &0.002 &0.000 &0.000 &0.000& 0.000 &0.000 &0.000 &0.000 &0.000 &0.000& $(1.743 ,1.496)$ &-0.1856 \\
\\
\multicolumn{13}{l}{Portmanteau tests for adequacy of the CCC-APGARCH$(0,3)$ when $\underline{\delta}_0$ is unknown}\\
(USD,JPY) & 0.952 &0.361 &0.064& 0.000& 0.000 &0.000 &0.000 &0.000 &0.000 &0.000& 0.000& 0.000& $(1.245, 1.518)$ &-0.2031\\
\\
\multicolumn{13}{l}{Portmanteau tests for adequacy of the CCC-APGARCH$(1,1)$ when $\underline{\delta}_0$ is unknown}\\
(USD,JPY) & 0.044 &0.084 &0.167& 0.040 &0.031& 0.046& 0.077& 0.111& 0.118& 0.102& 0.129 &0.087 & $(2.149, 1.568)$ &-0.3520\\
\\
\multicolumn{13}{l}{Portmanteau tests for adequacy of the CCC-APGARCH$(2,1)$ when $\underline{\delta}_0$ is unknown}\\
(USD,JPY) & 0.115 &0.262 &0.253 &0.134& 0.101& 0.151 &0.207 &0.219& 0.194 &0.234& 0.260 &0.269 & $(2.229 ,1.578)$ &-0.2756\\
\\
\multicolumn{13}{l}{Portmanteau tests for adequacy of the CCC-APGARCH$(1,2)$ when $\underline{\delta}_0$ is unknown}\\
(USD,JPY) & 0.681& 0.288& 0.336& 0.038 &0.018& 0.027& 0.042 &0.035& 0.028 &0.025 &0.032 &0.036& $(1.618, 4.547)$ &-0.2805 \\
\\
\multicolumn{13}{l}{Portmanteau tests for adequacy of the CCC-APGARCH$(2,2)$ when $\underline{\delta}_0$ is unknown}\\
(USD,JPY) & 0.237 &0.412 &0.407 &0.188& 0.120 &0.176 &0.213 &0.106 &0.070 &0.081 &0.103 &0.097& $(0.921, 2.033)$ & -0.2783 \\
\hline\hline
\end{tabular}
\end{center}
}
\label{tabexchange}
\end{table}

\begin{table}[H]
 \caption{\small{Empirical size of the proposed test: relative frequencies (in \%)
 of rejection of an APGARCH$(0,1)$. }}
  {\scriptsize 
\begin{center}
\begin{tabular}{ccc cccccccccccc}
\hline\hline
$\delta_0$& Length $n$ & Level $\alpha$ & \multicolumn{12}{c}{Lag $m$}\\
\hline
&&&1&2&3&4&5&6&7&8&9&10&11&12\\
&& $1\%$&0.7 & 0.7 & 0.7 & 1.3 & 1.1 & 0.7 & 0.7 & 0.5 & 0.6&  0.5 & 0.4 & 0.4\\
$(1,1)$&$250$& $5\%$& 3.6 & 4.8 & 5.3 & 6.3 & 6.6 & 5.2 & 5.5 & 4.9 & 4.6 & 4.6 & 4.4 & 4.0\\
 && $10\%$&  8.6 & 9.8 &11.4& 11.7 &11.9& 11.2& 11.4& 10.6& 10.0 & 9.7 & 9.3 & 9.7\\
\cline{3-15}
&& $1\%$&0.8 & 0.7 & 0.7 & 0.8&  1.2 & 1.4 & 0.9 & 1.1 & 0.9 & 0.9 & 0.8 & 0.8\\
$(1,1)$&$500$& $5\%$& 4.4 & 4.2 & 4.3 & 4.3 & 4.7 & 4.3  &3.8  &3.4  &3.3  &3.9  &3.6 & 3.6\\
 && $10\%$& 8.5&  9.1& 10.3 & 8.9 & 9.8 &10.2 &10.4&  9.4 & 9.3 & 8.2 & 8.5 & 9.0\\
\cline{3-15}
&& $1\%$&1.2 & 1.6 & 1.4 & 0.9 & 0.9 & 0.7 & 1.3 & 1.2 & 1.2 & 1.4 & 1.4 & 1.4\\
$(1,1)$&$2,000$& $5\%$&4.1 & 4.7 & 5.3 & 5.5 & 6.0 & 5.6 & 6.0 & 5.0 & 5.1&  5.8 & 5.7 & 5.3 \\
 && $10\%$& 8.9 & 9.4 & 9.7& 10.4& 10.8& 11.1 &12.6& 10.9& 12.1& 11.5& 11.4 &10.6\\
 \hline
&&&1&2&3&4&5&6&7&8&9&10&11&12\\
&& $1\%$& 0.6 & 0.9 & 0.5 & 1.2 & 0.9 & 0.8 & 0.4 & 0.3 & 0.5 & 0.4 & 0.2 & 0.3\\
$(0.8,1.5)$&$250$& $5\%$& 2.4 & 3.9 & 4.8 & 5.8 & 5.0 & 5.1 & 5.4 & 4.3 & 4.5 & 3.9 & 4.0 & 3.8\\
 && $10\%$& 8.6 & 7.7 &10.7& 10.3& 11.1& 11.0 &10.7 & 9.6 & 9.4 & 9.1 & 8.4 & 8.4\\
\cline{3-15}
&& $1\%$& 0.6 & 0.5 & 0.9 & 0.7 & 1.1 & 1.2 & 1.1 & 1.0 & 0.7 & 0.7 & 0.8 & 0.9\\
$(0.8,1.5)$&$500$& $5\%$& 4.2 & 3.9 & 3.9 & 3.9 & 4.4 & 4.1 & 4.3 & 3.9 & 3.5 & 4.7 & 3.7 & 4.5\\
 && $10\%$& 8.4 & 8.1 & 9.6 & 8.3 & 9.8 & 9.7 & 9.7&  9.1 & 9.1 & 8.5 & 7.7 & 8.7\\
\cline{3-15}
&& $1\%$&1.1 & 1.7 & 1.3 & 1.4 & 0.9 & 0.9 & 1.5 & 1.3 & 1.2 & 1.5 & 1.6 & 1.2\\
$(0.8,1.5)$&$2,000$& $5\%$& 4.6 & 4.9 & 5.4 & 5.7 & 5.8 & 5.6 & 5.8 & 5.6 & 5.5 & 5.9 & 6.2 & 5.2\\
 && $10\%$& 9.0& 10.0 & 9.0 & 9.7 &10.6 &11.2& 11.8& 11.5& 11.2& 11.3 &10.7 &11.1\\
\hline
&&&1&2&3&4&5&6&7&8&9&10&11&12\\
&& $1\%$& 0.8 & 0.9 & 1.2 & 1.7 & 1.2 & 1.7 & 1.0 & 0.9&  0.8 & 0.8 & 0.5 & 0.6\\
$(2,2)$&$250$& $5\%$& 4.2 & 6.0 & 6.7 & 6.9 & 6.4 & 6.7 & 6.4 & 5.9 & 5.7 & 5.1 & 4.9 & 4.6\\
 && $10\%$& 9.1 &11.8& 12.3& 12.4 &12.4& 12.2& 11.8 &11.0& 10.5& 10.6& 10.6 &10.0\\
\cline{3-15}
&& $1\%$& 0.6 & 0.8 & 0.8 & 1.0 & 1.3 & 1.1 & 1.2 & 1.2 & 1.1 & 0.6 & 0.7 & 0.8\\
$(2,2)$&$500$& $5\%$& 5.1 & 5.3 & 4.9 & 4.4 & 5.5 & 4.6 & 4.5 & 3.9 & 3.8 & 4.3 & 4.6 & 4.4\\
 && $10\%$& 9.9& 10.2& 10.2 & 9.7 & 9.4 &10.2 & 9.6 & 9.2 & 9.2& 10.0 & 8.7 & 9.7\\
\cline{3-15}

&& $1\%$& 1.2 & 1.7&  1.3 & 1.2 & 1.2 & 1.1 & 1.2 & 1.2 & 1.4 & 1.7 & 1.5 & 1.6\\
$(2,2)$&$2,000$& $5\%$& 4.5 & 5.5 & 5.7 & 6.1 & 6.3 & 6.5 & 6.3 & 5.0&  5.9 & 5.7 & 5.8 & 5.4\\
 && $10\%$& 8.5&  9.4& 10.9& 10.8& 10.5& 11.2& 12.5& 11.7 &12.2 &11.8 &11.4& 10.5\\
\hline
&&&1&2&3&4&5&6&7&8&9&10&11&12\\
&& $1\%$& 0.9 & 0.9 & 1.2 & 1.7 & 1.4 & 1.8 & 1.3 & 1.3 & 0.9 & 0.8 & 0.7 & 0.6\\
$(3,2.5)$&$250$& $5\%$& 4.7&  5.7 & 6.7 & 6.9 & 6.2 & 6.6 & 5.9 & 5.7 & 5.6 & 4.9 & 4.7 & 4.7\\
 && $10\%$& 9.6 &12.0& 11.7 &12.6& 12.3& 12.3 &11.2 &11.2 & 9.3 &10.5 &10.3 & 8.9\\
\cline{3-15}
&& $1\%$&  0.6 & 0.9 & 0.9 & 1.0 & 1.3 & 1.0 & 1.1 & 1.3 & 1.2 & 0.6 & 0.7 & 0.7\\
$(3,2.5)$&$500$& $5\%$& 5.3 & 5.3 & 5.6 & 4.7 & 5.7 & 4.6 & 4.5 & 4.1 & 4.2 & 4.3 & 4.6 & 4.6\\
 && $10\%$& 11.1& 10.4 & 9.9& 10.1& 10.0& 10.7 & 9.9 &10.1 & 9.2 & 9.5 & 8.9 & 9.6\\
\cline{3-15}

&& $1\%$& 1.2 & 1.6 & 1.3 & 1.3 & 1.1 & 1.0 & 1.0 & 1.3 & 1.4 & 1.7 & 1.6 & 1.6\\
$(3,2.5)$&$2,000$& $5\%$& 4.6 & 4.9 & 5.8 & 6.2&  6.0&  6.4 & 6.5 & 5.3 & 5.5 & 5.7 & 5.8 & 5.2\\
 && $10\%$& 9.2 & 9.2& 10.5 &10.8& 10.2& 11.1& 12.1 &12.4 &12.2 &11.9 &11.3 &10.8\\
\hline\hline 
\multicolumn{14}{l}{ Model \eqref{arch}--\eqref{coef} with $A_{01}^+=A_{01}^-$ when $\underline{\delta}_0$ is known.}
\\
\end{tabular}
\end{center}
}
\label{archsym}
\end{table}

\begin{table}[H]
 \caption{\small{Empirical size of the proposed test: relative frequencies (in \%)
 of rejection of an APGARCH$(0,1)$.}}
  {\scriptsize 
\begin{center}
\begin{tabular}{ccc cccccccccccc}
\hline\hline
$\delta_0$& Length $n$ & Level $\alpha$ & \multicolumn{12}{c}{Lag $m$}\\
\hline
&&&1&2&3&4&5&6&7&8&9&10&11&12\\
&& $1\%$&0.8&  1.2 & 0.9 & 1.6 & 1.6 & 0.9 & 0.8 & 0.9 & 0.7 & 0.5  &0.4 & 0.5\\
$(1,1)$&$250$& $5\%$& 4.1 & 4.9 & 6.1 & 5.9 & 6.4 & 5.7 & 5.7 & 5.3 & 4.7 & 4.2 & 4.6&  4.4\\
 && $10\%$& 8.7 &10.7& 12.1 &12.5 &12.4 &11.7 &11.4 &10.5 &11.4 &10.2 & 9.3 & 9.5\\
\cline{3-15}
&& $1\%$&1.1 & 1.1 & 1.0 & 1.0 & 1.4 & 1.6 & 1.1 & 1.2 & 1.0 & 0.8 & 0.9 & 1.0\\
$(1,1)$&$500$& $5\%$& 4.8 & 5.5 & 5.1 & 4.6 & 5.2&  4.6 & 4.5 & 3.9 & 4.0 & 4.3 & 4.3 & 4.2\\
 && $10\%$& 9.6 &10.6 &11.7 &10.9  &9.6 &10.9 &10.6 & 9.6 & 9.3 & 8.9 & 8.5 & 9.9\\
\cline{3-15}
&& $1\%$&1.2 & 1.7 & 1.6 & 1.4 & 1.1 & 1.3 & 1.5 & 1.2 & 1.2 & 1.5 & 1.7 & 1.4\\
$(1,1)$&$2,000$& $5\%$&4.4 & 5.2 & 5.8 & 5.1 & 5.6 & 5.6 & 6.4 & 5.0 & 6.0 & 6.3 & 6.1 & 5.2\\
 && $10\%$& 8.7 & 9.2 &10.6 & 9.7 &11.3 &11.1& 12.2& 10.7& 11.7 &11.9& 11.5& 11.1\\
 \hline
&&&1&2&3&4&5&6&7&8&9&10&11&12\\
&& $1\%$& 0.5 & 0.9 & 1.3 & 1.3 & 1.4  &1.3 & 1.0 & 0.8  &0.8 & 0.5 & 0.3 & 0.6\\
$(0.8,1.5)$&$250$& $5\%$& 3.7 & 4.2 & 5.1 & 5.9 & 6.1 & 5.5 & 5.7 & 5.1 & 4.6 & 4.5 & 4.0 & 3.9\\
 && $10\%$& 9.3 & 9.9 &11.6& 11.3 &11.4& 11.6 &10.3 &10.5 &10.5  &9.9 & 9.1  &9.0\\
\cline{3-15}
&& $1\%$& 0.9 & 0.9 & 1.2 & 1.4 & 1.7 & 1.5 & 1.4 & 1.2 & 0.9 & 0.8 & 0.9  &0.9\\
$(0.8,1.5)$&$500$& $5\%$& 4.4 & 5.2 & 5.2 & 4.3 & 5.1 & 4.8 & 4.6 & 4.3 & 4.2 & 4.7 & 4.8 & 4.6\\
 && $10\%$& 9.2& 10.6& 11.0& 10.6& 10.2 & 9.8 & 9.7 & 8.8 & 9.0 & 9.2  &8.3 &10.3\\
\cline{3-15}
&& $1\%$&1.1 & 1.5 & 1.2 & 1.3 & 1.7 & 1.2 & 1.5 & 1.3 & 1.1 & 1.5 & 1.8 & 1.5\\
$(0.8,1.5)$&$2,000$& $5\%$& 4.9 & 5.2 & 5.6 & 5.3 & 5.6 & 6.1 & 6.1 & 5.3 & 6.2 & 6.3 & 6.0 & 5.3\\
 && $10\%$& 8.5& 10.3 &10.9 &10.1& 10.9& 10.6 &12.1 &11.6 &11.4& 11.7& 11.8 &10.6\\
\hline
&&&1&2&3&4&5&6&7&8&9&10&11&12\\
&& $1\%$& 1.2 & 1.3 & 1.1 & 1.7 & 1.5 & 1.6 & 1.2 & 1.0 & 0.9 & 0.7 & 0.3 & 0.5\\
$(2,2)$&$250$& $5\%$& 5.6 & 6.5 & 7.2 & 7.5 & 6.7 & 6.6 & 6.4 & 6.2 & 5.6 & 4.9 & 4.6 & 4.6\\
 && $10\%$& 10.2& 11.8& 13.4& 12.7 &13.4& 13.0& 12.3& 11.3& 11.1 &11.6 &10.0 &10.3\\
\cline{3-15}
&& $1\%$& 0.7 & 1.1 & 0.9 & 1.2 & 1.5 & 1.5 & 1.3  &1.2  &1.0 & 1.0 & 1.0 & 1.1\\
$(2,2)$&$500$& $5\%$& 5.8 & 5.8 & 5.8 & 5.0 & 5.6 & 5.2 & 5.2  &4.4 & 3.9&  4.7 & 4.7 & 5.0\\
 && $10\%$& 10.5 &11.5 &11.4& 11.4& 10.7& 11.2 &10.5& 10.2& 10.6& 10.0 & 8.2 & 9.6\\
\cline{3-15}

&& $1\%$& 1.3 & 1.4 & 1.1 & 1.2 & 1.4 & 1.3 & 1.3 & 1.0 & 1.3 & 1.6 & 1.9  &1.7\\
$(2,2)$&$2,000$& $5\%$& 4.5 & 5.1 & 5.3 & 5.0 & 5.9 & 6.2 & 5.9 & 5.2 & 6.0 & 6.1 & 5.5 & 5.2\\
 && $10\%$& 9.3 & 9.5 &10.4 & 9.9 &10.8& 10.9 &11.6 &11.4& 11.8 &12.1& 11.0 &10.7\\
\hline
&&&1&2&3&4&5&6&7&8&9&10&11&12\\
&& $1\%$& 1.4 & 1.2 & 1.4&  1.7 & 1.8 & 1.8 & 1.6 & 1.1 & 1.0 & 0.8 & 0.4 & 0.4\\
$(3,2.5)$&$250$& $5\%$& 5.0 & 6.7 & 7.3 & 7.7 & 6.5 & 7.1 & 6.7 & 5.9 & 5.6 & 5.1 & 4.9 & 4.3\\
 && $10\%$& 9.7 &11.8 &13.6 &12.9& 13.1 &12.8 &11.9 &11.5 &11.4 &11.3 &10.5 & 9.2\\
\cline{3-15}
&& $1\%$& 1.0 & 1.1&  1.0 & 1.2 & 1.4 & 1.3 & 1.2 & 1.2 & 1.3 & 1.0 & 1.0 & 1.0\\
$(3,2.5)$&$500$& $5\%$& 5.9 & 6.0 & 5.5 & 4.9 & 5.3 & 5.0 & 4.8 & 4.2 & 3.8 & 4.4 & 4.3 & 4.6\\
 && $10\%$& 11.0& 10.7& 10.9& 11.2& 10.5& 11.3 &10.9 & 9.7& 10.2 & 9.4 & 8.5 & 9.6\\
\cline{3-15}

&& $1\%$& 1.0 & 1.5 & 0.9 & 1.2 & 1.3 & 1.2 & 1.5 & 1.0 & 1.3 & 1.6 & 1.7 & 1.8\\
$(3,2.5)$&$2,000$& $5\%$& 5.0 & 5.0 & 5.5 & 5.6 & 5.9 & 6.1 & 6.6 & 4.8 & 5.8 & 5.8 & 5.3 & 5.2\\
 && $10\%$& 9.9 & 9.2& 10.5& 10.2 &10.9 &10.8 &11.2& 11.7& 12.1& 11.8& 11.2 &10.1\\
\hline\hline 
\multicolumn{14}{l}{ Model \eqref{arch}--\eqref{coef} with $A_{01}^+\neq A_{01}^-$ when $\underline{\delta}_0$ is known.}
\\
\end{tabular}
\end{center}
}
\label{archasym}
\end{table}

\begin{table}[H]
 \caption{\small{Empirical power of the proposed test for the null hypothesis of the CCC-APGARCH$(0,1)$ model defined by \eqref{arch} with $\underline{\delta}_0=(1,1)$ against the alternative given by \eqref{arch} when $\underline{\delta}_0\neq (1,1)$.}}
  {\scriptsize 
\begin{center}
\begin{tabular}{ccc cccccccccccc}
\hline\hline
$\delta_0$& Length $n$ & Level $\alpha$ & \multicolumn{12}{c}{Lag $m$}\\
\hline
&&&1&2&3&4&5&6&7&8&9&10&11&12\\
&& $1\%$&0.8&  1.2 & 0.9 & 1.6 & 1.6 & 0.9 & 0.8 & 0.9 & 0.7 & 0.5  &0.4 & 0.5\\
$(1,1)$&$250$& $5\%$& 4.1 & 4.9 & 6.1 & 5.9 & 6.4 & 5.7 & 5.7 & 5.3 & 4.7 & 4.2 & 4.6&  4.4\\
 && $10\%$& 8.7 &10.7& 12.1 &12.5 &12.4 &11.7 &11.4 &10.5 &11.4 &10.2 & 9.3 & 9.5\\
\cline{3-15}
&& $1\%$&1.1 & 1.1 & 1.0 & 1.0 & 1.4 & 1.6 & 1.1 & 1.2 & 1.0 & 0.8 & 0.9 & 1.0\\
$(1,1)$&$500$& $5\%$& 4.8 & 5.5 & 5.1 & 4.6 & 5.2&  4.6 & 4.5 & 3.9 & 4.0 & 4.3 & 4.3 & 4.2\\
 && $10\%$& 9.6 &10.6 &11.7 &10.9  &9.6 &10.9 &10.6 & 9.6 & 9.3 & 8.9 & 8.5 & 9.9\\
\cline{3-15}
&& $1\%$&1.2 & 1.7 & 1.6 & 1.4 & 1.1 & 1.3 & 1.5 & 1.2 & 1.2 & 1.5 & 1.7 & 1.4\\
$(1,1)$&$2,000$& $5\%$&4.4 & 5.2 & 5.8 & 5.1 & 5.6 & 5.6 & 6.4 & 5.0 & 6.0 & 6.3 & 6.1 & 5.2\\
 && $10\%$& 8.7 & 9.2 &10.6 & 9.7 &11.3 &11.1& 12.2& 10.7& 11.7 &11.9& 11.5& 11.1\\
 \hline
&&&1&2&3&4&5&6&7&8&9&10&11&12\\
&& $1\%$& 0.5 & 1.2 & 1.0 & 1.6 & 1.4 & 1.3 & 0.8 & 0.9 & 0.7 & 0.7 & 0.5 & 0.4\\
$(0.8,1.5)$&$250$& $5\%$& 4.4 & 5.2 & 6.3 & 5.9 & 5.3 & 6.4 & 5.9 & 5.8 & 5.3 & 4.9 & 4.5 & 4.7\\
 && $10\%$& 9.5& 10.7 &11.4 &12.4& 12.7& 12.4 &11.6& 11.2& 10.3& 10.3 &10.1& 10.0 \\
\cline{3-15}
&& $1\%$& 1.0 & 1.0 & 1.1 & 1.0 & 1.4 & 1.5 & 1.3 & 1.2 & 1.0 & 1.0 & 0.9 & 1.1\\
$(0.8,1.5)$&$500$& $5\%$& 5.3 & 5.2 & 5.4 & 5.0 & 5.3 & 4.5 & 4.4 & 4.1 & 4.2 & 4.5 & 4.5 & 4.4\\
 && $10\%$& 10.3 &10.5& 11.4& 11.0& 10.8& 10.1& 10.0 & 9.2 & 9.2& 10.0 & 9.2 & 9.2\\
\cline{3-15}
&& $1\%$&1.4 & 1.5 & 1.2 & 1.2 & 1.7 & 1.3 & 1.4 & 1.4 & 1.4 & 1.2 & 1.4 & 1.3\\
$(0.8,1.5)$&$2,000$& $5\%$& 4.8 & 5.3 & 5.7 & 5.5 & 5.9 & 6.4 & 5.8 & 5.8 & 6.7 & 6.3 & 6.2 & 6.0\\
 && $10\%$& 9.9 &10.1 &11.0 &10.5 &11.6& 11.6& 12.1& 11.6 &11.8 &12.9 &12.2 &11.4\\
\hline
&&&1&2&3&4&5&6&7&8&9&10&11&12\\
&& $1\%$& 1.8 & 2.4 & 1.7 & 1.9 & 2.0 & 1.9 & 1.5 & 1.1 & 0.7 & 0.6 & 0.5 & 0.7\\
$(2,2)$&$250$& $5\%$& 8.4 & 9.1 & 9.8 & 8.6 & 8.5 & 7.8 & 8.0 & 7.0 & 5.8 & 5.8 & 5.0 & 5.1\\
 && $10\%$& 14.8& 16.2& 16.1& 16.1& 15.7& 16.0& 13.7& 12.4& 13.2& 11.8& 11.2& 10.7\\
\cline{3-15}
&& $1\%$& 3.1 & 3.4 & 3.2 & 2.7 & 2.3 & 2.6 & 2.8 & 2.1 & 1.9 & 1.9 & 2.0 & 1.9\\
$(2,2)$&$500$& $5\%$& 11.1 &13.6& 11.7& 10.8 & 9.4 & 9.5 & 8.1 & 7.6 & 7.6 & 7.6 & 6.3 & 7.1\\
 && $10\%$& 18.8& 19.9 &19.8& 18.6 &17.6& 16.2& 15.9& 15.0& 14.7& 14.0& 13.0& 14.0\\
\cline{3-15}

&& $1\%$& 10.3 &11.0 & 9.8 & 8.3&  8.1 & 6.9 & 6.4 & 5.7&  5.9 & 5.5 & 5.6 & 4.8\\
$(2,2)$&$2,000$& $5\%$& 25.4 &29.1 &25.1 &23.0& 21.2& 19.3& 17.9& 18.0 &17.4& 16.6& 15.7 &14.7\\
 && $10\%$& 35.8 &39.5 &36.5 &34.5& 31.9 &30.1 &29.7 &28.5 &27.0 &26.2 &25.9 &26.0\\
\hline
&&&1&2&3&4&5&6&7&8&9&10&11&12\\
&& $1\%$& 1.9 & 2.4 & 1.8 & 1.7 & 2.1 & 2.0 & 1.6 & 0.9 & 0.9 & 0.6 & 0.4 & 0.7\\
$(3,2.5)$&$250$& $5\%$& 8.2 & 9.7 & 9.5 & 9.2 & 8.2 & 7.8 & 7.1 & 7.1 & 5.8 & 5.3 & 5.3 & 5.1\\
 && $10\%$& 15.2 &16.7& 16.1 &16.7 &16.0 &15.8 &14.4& 12.7& 12.7 &12.4 &10.6 &11.2\\
\cline{3-15}
&& $1\%$& 3.7 & 3.6 & 3.7 & 2.7 & 2.5 & 2.6 & 2.6 & 2.5 & 1.8 & 2.0 & 2.2 & 2.2\\
$(3,2.5)$&$500$& $5\%$& 11.5 &13.1 &12.3& 11.5 & 9.9 & 9.8 & 9.0 & 8.6 & 8.4 & 7.9 & 7.2 & 8.1\\
 && $10\%$& 19.4 &21.3& 20.3& 19.7 &17.9 &17.1 &17.0 &16.2 &14.9 &14.3 &14.4& 14.5\\
\cline{3-15}

&& $1\%$& 11.0 &12.8 &10.9 &10.0 & 9.0 & 8.2 & 7.5 & 7.0 & 6.8 & 6.4 & 6.0 & 5.4\\
$(3,2.5)$&$2,000$& $5\%$& 26.3 &31.7& 28.4& 25.7& 22.7& 21.2& 20.5& 20.8& 19.4& 18.3& 17.5& 17.5\\
 && $10\%$& 38.7 &44.3 &40.5& 38.0& 35.4& 33.7& 32.6 &31.0 &29.8 &29.4& 28.6 &27.4\\
\hline\hline 
\multicolumn{14}{l}{ Model \eqref{arch}--\eqref{coef} with $A_{01}^+\neq A_{01}^-$ when $\underline{\delta}_0$ is known.}
\\
\end{tabular}
\end{center}
}
\label{puiss-archasym}
\end{table}

\begin{table}[H]
 \caption{\small{Empirical power of the proposed test for the null hypothesis of a bivariate CCC-APGARCH$(0,1)$ model of the form \eqref{arch} against the bivariate CCC-APGARCH$(1,1)$ alternative given by \eqref{model-puiss} when $\underline{\delta}_0$ is known.}}
  {\scriptsize 
\begin{center}
\begin{tabular}{ccc cccccccccccc}
\hline\hline
$\delta_0$& Length $n$ & Level $\alpha$ & \multicolumn{12}{c}{Lag $m$}\\
\hline
&&&1&2&3&4&5&6&7&8&9&10&11&12\\
&& $1\%$&11.1 &12.7 &21.2 &22.6& 20.2& 16.5 &13.5& 10.7 & 9.4 & 7.1 & 6.5 & 5.5\\
$(1,1)$&$250$& $5\%$& 25.0& 36.0 &49.7& 50.5& 48.0 &44.3& 40.3 &35.6& 32.0 &28.9& 26.8 &24.8\\
 && $10\%$& 35.1 &52.2 &66.6 &68.9 &65.3 &61.4& 56.7& 52.0& 48.9 &45.9& 42.0& 39.9\\
\cline{3-15}
&& $1\%$&27.6& 42.6 &65.1 &72.4 &70.6 &66.8 &63.7 &57.7& 52.3 &48.8& 44.0 &39.6\\
$(1,1)$&$500$& $5\%$& 46.1 &71.7& 89.4& 92.2& 91.1& 90.4& 86.7 &83.3 &80.7 &77.0 &75.7 &72.3\\
 && $10\%$& 57.8 &84.1 &94.9& 96.3 &96.2 &95.1 &94.0 &92.6 &89.9 &87.6 &85.2& 83.7 \\
\cline{3-15}
&& $1\%$&78.0& 97.5 &99.3 &99.6 &99.8& 99.8& 99.8 &99.8& 99.8 &99.7& 99.7& 99.8\\
$(1,1)$&$2,000$& $5\%$& 87.9 &99.3 &99.7 &99.9 &99.9 &99.9& 99.9 &99.9 &99.8 &99.8 &99.8 &99.8\\
 && $10\%$& 92.3 &99.6 &99.9 &99.9 &99.9 &99.9 &99.9 &99.9 &99.9& 99.9 &99.9& 99.9\\
 \hline
&&&1&2&3&4&5&6&7&8&9&10&11&12\\
&& $1\%$& 12.5& 12.6 &18.8& 20.4& 16.3 &13.9 &11.6 & 9.6 & 7.2 & 6.2 & 5.7 & 5.0\\
$(0.8,1.5)$&$250$& $5\%$& 24.4& 34.2& 47.1& 48.3& 46.3 &41.2 &37.4& 32.5& 28.0 &25.0 &24.6 &21.7\\
 && $10\%$& 34.2& 49.9& 64.1& 64.7 &62.6 &58.6& 53.6& 50.6 &46.6 &44.3 &40.2 &37.2 \\
\cline{3-15}
&& $1\%$& 26.8& 37.7 &55.9 &62.7& 61.3& 58.3& 53.4 &47.8& 42.4& 38.9 &34.6& 29.1\\
$(0.8,1.5)$&$500$& $5\%$& 45.3 &65.9& 80.8 &86.6& 85.3& 83.2& 80.8 &77.6& 73.6 &70.9 &68.7 &62.8\\
 && $10\%$& 55.7& 79.4& 91.1& 93.1& 93.0& 91.8 &89.9 &87.8 &85.6 &83.0& 79.9 &77.8\\
\cline{3-15}
&& $1\%$&71.9 &90.9 &95.7& 97.8 &98.2 &98.0 &98.0 &98.1 &97.6 &97.1 &96.6 &96.9\\
$(0.8,1.5)$&$2,000$& $5\%$& 84.4 &96.7& 98.3 &98.8& 99.3 &99.2& 99.2& 99.3& 99.1 &98.7& 98.8 &98.8\\
 && $10\%$& 87.4 &98.3 &98.7 &99.3 &99.3 &99.3& 99.3 &99.3 &99.3 &99.3 &99.1 &99.3\\
\hline
&&&1&2&3&4&5&6&7&8&9&10&11&12\\
&& $1\%$& 14.6& 13.3 &21.8 &23.6& 22.3 &19.0 &16.3& 12.8& 10.7 &10.4 & 9.1 & 7.4\\
$(2,2)$&$250$& $5\%$& 30.4 &38.3& 53.0 &57.0& 53.8& 49.5 &45.5 &39.8 &36.1 &33.2 &30.1 &27.8\\
 && $10\%$& 41.7& 53.3 &68.8 &73.2 &71.5& 67.0 &63.4 &57.7 &54.8 &51.6& 46.9 &44.3\\
\cline{3-15}
&& $1\%$& 34.9 &44.7 &66.3 &74.0& 74.5 &72.6 &68.2& 64.4 &59.6 &55.9& 50.4 &45.8\\
$(2,2)$&$500$& $5\%$& 55.1 &70.9& 88.6 &92.3 &91.5& 90.8& 89.1& 87.4 &85.4& 81.9 &79.2 &77.5\\
 && $10\%$& 64.4 &82.7 &94.6 &96.2 &97.0& 96.0& 95.3 &94.3& 92.8 &91.1 &89.6 &87.7\\
\cline{3-15}

&& $1\%$& 84.6 &96.2 &98.7 &99.4 &99.6 &99.8 &99.8 &99.8 &99.7 &99.7 &99.4 &99.3 \\
$(2,2)$&$2,000$& $5\%$& 92.9 & 99.1  &99.5  &99.8& 100.0 & 99.9 & 99.8 & 99.9  &99.9  &99.9  &99.9  &99.9\\
 && $10\%$& 95.0 & 99.5  &99.8& 100.0& 100.0 &100.0 &100.0& 100.0 & 99.9 & 99.9& 100.0 &100.0\\
\hline
&&&1&2&3&4&5&6&7&8&9&10&11&12\\
&& $1\%$& 20.1 &17.3 &23.9 &26.3 &25.3& 22.9 &20.2 &18.5 &16.2 &12.7 &11.6 & 9.5\\
$(3,2.5)$&$250$& $5\%$& 41.7 &42.2 &54.8 &57.6 &55.7 &52.2 &49.0 &46.3& 42.3 &39.0& 37.8 &34.7\\
 && $10\%$& 51.6& 59.0 &69.6& 73.1& 73.1& 69.5 &67.1 &62.9& 60.1 &57.5& 54.0 &50.6\\
\cline{3-15}
&& $1\%$& 45.8 &48.3 &59.7& 66.0& 68.8 &67.7& 66.0 &66.0& 61.4& 58.0& 54.7 &50.4\\
$(3,2.5)$&$500$& $5\%$& 67.0 &71.7 &83.4 &86.2& 87.8& 87.7 &86.0& 86.2 &85.0& 82.3 &80.3& 77.5\\
 && $10\%$& 76.5& 82.2 &91.4 &92.2 &93.5 &93.7 &93.3 &92.6 &91.9& 90.0 &89.7 &87.9\\
\cline{3-15}

&& $1\%$& 86.2 &91.5 &95.4 &96.1 &96.9 &97.7 &97.9 &97.9 &98.1 &98.0& 97.9& 98.2\\
$(3,2.5)$&$2,000$& $5\%$& 93.2& 96.9 &98.1& 98.3 &98.9 &98.8 &99.0 &99.2 &99.1 &98.8 &98.9& 98.9\\
 && $10\%$& 94.6& 97.6& 98.7& 99.0& 99.2 &99.4& 99.5& 99.5& 99.5& 99.4& 99.5& 99.4\\
\hline\hline 
\multicolumn{14}{l}{ Model \eqref{model-puiss} with $A_{01}^+\neq A_{01}^-$ when $\underline{\delta}_0$ is known.}
\\
\end{tabular}
\end{center}
}
\label{puiss-archasymbis}
\end{table}

\begin{table}[H]
 \caption{\small{Empirical size of the proposed test: relative frequencies (in \%)
 of rejection of an APGARCH$(0,1)$.}}
  {\scriptsize 
\begin{center}
\begin{tabular}{ccc cccccccccccc}
\hline\hline
$\delta_0$& Length $n$ & Level $\alpha$ & \multicolumn{12}{c}{Lag $m$}\\
\hline
&&&1&2&3&4&5&6&7&8&9&10&11&12\\
&& $1\%$&0.6 & 1.1 & 1.0 & 1.6 & 1.3 & 1.0 & 1.0 & 1.0 & 0.6 & 0.3 & 0.3 & 0.5\\
$(1,1)$&$250$& $5\%$& 5.1 & 5.4 & 6.1 & 7.2 & 6.4 & 6.4 & 6.1 & 5.6 & 5.5 & 5.2 & 4.8 & 4.8\\
 && $10\%$&  10.4 &10.3 &11.4& 11.4 &12.5& 12.3 &11.8 &10.5& 10.1 &11.1 &10.7& 10.5\\
\cline{3-15}
&& $1\%$&0.5 & 0.9 & 1.1 & 1.1 & 1.3 & 1.0 & 1.0 & 1.0 & 0.9 & 1.1 & 0.8 & 0.7\\
$(1,1)$&$500$& $5\%$& 4.3 & 4.9 & 4.1 & 4.8 & 4.3 & 4.8 & 4.8 & 4.3&  4.5 & 4.5 & 4.3 & 4.3\\
 && $10\%$& 9.1& 10.1 &10.2 & 8.9 & 9.8 & 9.5 & 9.8 &10.0 & 8.9 & 9.0 & 9.2 & 9.9\\
\cline{3-15}
&& $1\%$&0.8 & 1.3 & 1.1 & 1.0 & 1.4 & 0.8 & 1.2 & 1.2 & 1.4  &1.4 & 1.5 & 1.1\\
$(1,1)$&$2,000$& $5\%$&4.2 & 4.6 & 5.9 & 5.4 & 5.9 & 5.8 & 6.0 & 5.8 & 5.3 & 5.3 & 5.4  &4.7\\
 && $10\%$& 8.6 & 9.4& 10.7& 10.2& 10.0& 10.6 &12.0& 11.3 &12.2 &12.5& 11.4& 10.5\\
 \hline
&&&1&2&3&4&5&6&7&8&9&10&11&12\\
&& $1\%$&  0.5 & 0.8 & 1.3 & 1.2&  0.9 & 1.0 & 0.8  &0.5 & 0.6 & 0.5 & 0.5 & 0.5\\
$(0.8,1.5)$&$250$& $5\%$& 5.0 & 4.9 & 6.4 & 7.1 & 6.8 & 5.8 & 5.9 & 5.4 & 5.7 & 4.7 & 4.8 & 5.0\\
 && $10\%$& 9.4& 10.4& 12.0& 13.1 &13.0 &12.9& 12.2 &10.5 &11.0 &10.2 &10.0 &10.7 \\
\cline{3-15}
&& $1\%$& 0.9 & 1.1 & 1.0 & 1.3 & 1.5 & 1.1 & 1.2 & 1.2 & 1.2 & 0.9 & 0.8 & 1.0\\
$(0.8,1.5)$&$500$& $5\%$& 4.4 & 4.2 & 4.4 & 5.2 & 4.9 & 4.9 & 4.8&  4.6 & 4.3 & 4.7 & 4.3 & 4.8\\
 && $10\%$& 8.2 & 9.3  &9.8 & 9.9 &10.3& 10.6 &10.1 & 9.8 & 9.3 & 9.5 & 9.8 &10.1\\
\cline{3-15}
&& $1\%$&0.9 & 1.4 & 1.2 & 1.2 & 1.3 & 0.9 & 1.3 & 1.2 & 1.3 & 1.6 & 1.5 & 1.4\\
$(0.8,1.5)$&$2,000$& $5\%$& 4.4 & 5.4 & 5.6 & 5.6 & 6.1 & 5.4 & 5.9 & 6.0 & 6.3 & 5.7  &6.0 & 4.9\\
 && $10\%$& 9.4  &9.6 &11.2&10.5 &10.1& 11.7 &12.1 &11.3& 12.0 &11.9 &11.4 &10.8\\
\hline
&&&1&2&3&4&5&6&7&8&9&10&11&12\\
&& $1\%$& 1.0 & 1.6 & 1.6 & 1.8 & 2.0 & 1.6 & 1.2 & 1.2 & 0.8 & 1.2 & 1.4 & 1.0\\
$(2,2)$&$250$& $5\%$& 6.2 & 8.2 & 7.0 & 6.4 & 4.8 & 5.6 & 5.8 & 5.6 & 6.2 & 5.0 & 5.4 & 5.8\\
 && $10\%$& 12.0 &13.4& 14.4 &12.4 &10.8 &10.8 &11.0 &10.0 & 9.6& 11.8 &11.4 &10.6\\
\cline{3-15}
&& $1\%$& 0.9 & 0.9 & 0.9 & 1.3 & 1.2 & 1.1 & 1.2 & 1.2 & 1.0 & 0.8 & 0.7 & 0.9\\
$(2,2)$&$500$& $5\%$& 5.2 & 5.2 & 4.9 & 5.3 & 5.6 & 4.5 & 4.8 & 4.6 & 4.2 & 4.6 & 4.8 & 4.8\\
 && $10\%$& 10.5& 10.8 &11.0& 11.0& 10.3 &10.8 &10.3 &10.3 & 9.7& 10.1 & 9.1 &10.9\\
\cline{3-15}

&& $1\%$& 1.2 & 1.4 & 1.1 & 0.9 & 1.4 & 1.1 & 1.3 & 1.3 & 1.4 & 1.4 & 1.5 & 1.6\\
$(2,2)$&$2,000$& $5\%$& 5.4  &5.0 & 5.2 & 5.7 & 6.1  &6.4 & 6.8 & 5.5 & 5.8  &5.8 & 5.5  &5.9\\
 && $10\%$& 8.8 &10.5 &11.1& 10.8 &10.9 &11.2 &12.1 &11.9& 12.8& 12.4 &11.2& 10.8\\
\hline
&&&1&2&3&4&5&6&7&8&9&10&11&12\\
&& $1\%$& 1.0 & 2.6 & 2.0 & 1.6 & 1.6 & 1.8 & 1.4 & 1.4 & 1.2 & 1.6 & 1.4 & 1.0\\
$(3,2.5)$&$250$& $5\%$& 7.4 & 7.4 & 8.8 & 7.4 & 6.4 & 6.0 & 5.6 & 5.8 & 5.4 & 5.2 & 5.2 & 5.2\\
 && $10\%$& 11.8& 14.8& 13.0& 14.4& 11.4& 12.6& 11.6& 11.4& 10.4& 11.4& 11.6& 10.0\\
\cline{3-15}
&& $1\%$&  0.8 & 0.9 & 0.7 & 1.1 & 1.4 & 1.5 & 1.4 & 1.5 & 1.1 & 0.9 & 0.9 & 1.0\\
$(3,2.5)$&$500$& $5\%$& 4.9 & 4.9 & 5.6 & 5.2 & 5.8  &4.3 & 4.1 & 4.7 & 4.8 & 4.8 & 4.8 & 5.1\\
 && $10\%$& 11.0& 11.0& 10.5 &10.7 &10.9 &11.0& 11.2& 10.1 &10.2 & 9.6 & 9.5& 11.1\\
\cline{3-15}

&& $1\%$& 1.0  &1.4 & 1.2 & 1.4 & 1.0 & 1.2 & 1.1 & 1.4 & 1.5 & 1.4 & 1.6 & 1.9\\
$(3,2.5)$&$2,000$& $5\%$& 5.2 & 4.9  &5.5 & 5.9 & 6.5 & 6.4 & 6.9  &5.8 & 5.9  &6.0 & 6.0  &5.9\\
 && $10\%$& 9.7 &10.8 &11.4 &10.3& 10.5 &10.7 &12.2& 12.3 &12.9 &12.1 &11.7 &10.8\\
\hline\hline 
\multicolumn{14}{l}{ Model \eqref{arch}--\eqref{coef} with $A_{01}^+=A_{01}^-$ when $\underline{\delta}_0$ is unknown.}
\\
\end{tabular}
\end{center}
}
\label{archsymInc}
\end{table}

\begin{table}[H]
 \caption{\small{Empirical size of the proposed test: relative frequencies (in \%)
 of rejection of an APGARCH$(0,1)$.}}
  {\scriptsize 
\begin{center}
\begin{tabular}{ccc cccccccccccc}
\hline\hline
$\delta_0$& Length $n$ & Level $\alpha$ & \multicolumn{12}{c}{Lag $m$}\\
\hline
&&&1&2&3&4&5&6&7&8&9&10&11&12\\
&& $1\%$&0.7 & 0.9 & 1.4 & 1.9 & 1.1 & 1.2 & 0.8 & 0.5 & 0.7 & 0.4 & 0.3 & 0.6\\
$(1,1)$&$250$& $5\%$& 4.2 & 5.3 & 5.8 & 6.6 & 6.2 & 5.5 & 5.9 & 5.6&  5.2 & 4.7 & 4.8 & 4.4\\
 && $10\%$& 9.1 & 9.6& 12.3& 11.3& 12.6 &13.2& 12.0& 10.6 &11.0 &10.5 & 9.7 & 9.7\\
\cline{3-15}
&& $1\%$&0.5 & 0.9 & 1.1 & 1.0 & 1.5 & 1.4  &1.3 & 1.2 & 1.1 & 1.1 & 0.9 & 0.8\\
$(1,1)$&$500$& $5\%$& 4.6&  5.0 & 4.8 & 4.8 & 5.5 & 5.5 & 4.8 & 4.0 & 4.2 & 4.5 & 4.5 & 4.9\\
 && $10\%$& 9.7 &10.0 &10.4& 10.4 &10.1 &10.7& 11.0 &10.0 & 9.3 & 9.6 & 8.3 & 9.8\\
\cline{3-15}
&& $1\%$&1.5 & 1.4 & 1.4 & 1.2 & 1.5 & 1.3 & 1.6 & 1.0 & 1.3 & 1.8 & 1.8 & 1.2\\
$(1,1)$&$2,000$& $5\%$&4.5 & 4.6 & 5.1 & 4.9 & 5.4 & 5.8 & 6.1 & 5.4 & 5.8 & 6.1 & 6.0  &5.4\\
 && $10\%$& 9.4 & 9.2& 10.6 &10.4 &11.0& 10.6 &10.9& 11.4& 11.2 &12.0 &10.9 &10.5\\
 \hline
&&&1&2&3&4&5&6&7&8&9&10&11&12\\
&& $1\%$& 0.5 & 1.1 & 1.3 & 1.7 & 1.5 & 1.4 & 1.1 & 0.6 & 0.7 & 0.8 & 0.5 & 0.5\\
$(0.8,1.5)$&$250$& $5\%$& 4.5 & 5.1 & 6.2 & 6.7 & 5.5 & 5.7 & 6.0 & 5.1 & 5.0 & 5.1 & 4.3 & 4.4\\
 && $10\%$& 9.5& 10.0& 12.1& 11.9 &12.4& 12.7& 11.0& 10.7& 10.0& 10.0 & 9.6& 10.0\\
\cline{3-15}
&& $1\%$& 1.2 & 1.2&  1.1&  1.2 & 1.6 & 1.4 & 1.2 & 1.0 & 1.1 & 1.3 & 1.2 & 1.0\\
$(0.8,1.5)$&$500$& $5\%$& 5.4 & 5.7 & 4.8 & 5.4 & 5.7 & 4.9 & 5.1 & 4.4 & 4.6 & 4.7 & 4.9 & 5.3\\
 && $10\%$& 10.0 &10.3& 11.2 &10.3 &10.8 &10.6 &11.1 &10.0 & 9.7 & 9.4 & 9.1 &10.2\\
\cline{3-15}
&& $1\%$&1.4 & 1.3 & 1.0 & 1.2 & 1.8 & 1.4 & 1.4 & 1.3 & 1.4 & 1.7 & 1.9 & 1.4\\
$(0.8,1.5)$&$2,000$& $5\%$& 4.5 & 5.5 & 5.7 & 5.8& 5.5 & 5.8 & 6.9 & 5.6 & 6.3 & 6.0  &6.4 & 5.4\\
 && $10\%$& 9.2 &10.2 &11.3 &10.6& 10.1 &11.2 &11.3 &11.6 &12.1 &12.2 &12.2 &10.7\\
\hline
&&&1&2&3&4&5&6&7&8&9&10&11&12\\
&& $1\%$& 0.8 & 1.2 & 1.8 & 2.0 & 1.2 & 1.0 & 1.0 & 0.6 & 0.6 & 0.6 & 1.0 & 0.8\\
$(2,2)$&$250$& $5\%$& 5.0 & 6.4 & 6.4 & 6.4 & 4.8 & 4.8 & 5.4 & 5.2 & 5.8 & 5.4 & 5.0 & 5.4\\
 && $10\%$& 10.6& 11.2 &13.4& 12.4& 10.6 &10.6 &10.4 &10.2 &10.6 &11.6 &11.4 &10.0\\
\cline{3-15}
&& $1\%$& 1.0 & 1.2 & 1.1 & 1.2 & 1.5 & 1.6 & 1.8 & 1.4 & 1.2 & 1.3 & 1.0 & 0.9\\
$(2,2)$&$500$& $5\%$& 5.1 & 6.1 & 5.7 & 5.8 & 5.7 & 5.8 & 5.8 & 4.2 & 4.7 & 4.6 & 4.8 & 5.1\\
 && $10\%$& 9.8 &11.9 &12.1 &12.4& 11.6 &11.6 &11.1 &10.8 &10.1& 10.0 & 9.0 & 9.5\\
\cline{3-15}

&& $1\%$& 1.1 & 1.3 & 1.1  &1.2  &1.7 & 1.3 & 1.7&  1.2 & 1.3 & 1.8 & 1.7  &1.6\\
$(2,2)$&$2,000$& $5\%$& 4.8 & 5.0 & 5.8 & 5.5 & 6.1 & 5.9 & 6.1 & 5.4 & 5.7 & 5.5 & 5.5  &5.5\\
 && $10\%$& 10.2 &10.2 &10.1& 10.6 &10.8& 11.4& 11.3& 11.6 &12.4 &12.3 &11.2 &10.1\\
\hline
&&&1&2&3&4&5&6&7&8&9&10&11&12\\
&& $1\%$& 1.0 & 1.4 & 1.8 & 2.0 & 1.0 & 1.0 & 0.8 & 1.0 & 0.8 & 0.8 & 0.8 & 1.0\\
$(3,2.5)$&$250$& $5\%$& 6.0 & 6.4 & 7.4 & 6.2 & 4.2 & 5.2 & 6.2 & 4.2 & 5.2 & 5.4 & 5.6 & 4.6\\
 && $10\%$& 11.6 &11.8& 13.6 &12.2 &10.2 &10.2 &10.4 &10.6 & 9.8& 10.4& 10.2& 10.6\\
\cline{3-15}
&& $1\%$& 1.0 & 1.1 & 1.2 & 1.3 & 1.5 & 1.4 & 1.8 & 1.8 & 1.3 & 1.0 & 1.1&  1.0\\
$(3,2.5)$&$500$& $5\%$& 6.0 & 6.3 & 5.4 & 5.7 & 5.7 & 5.3 & 4.9 & 4.3 & 4.0 & 4.5&  4.8 & 5.1\\
 && $10\%$& 11.8& 12.9 &11.6& 12.0& 12.2 &12.2& 11.5& 11.0 &10.6 & 9.7 & 9.5 &10.3\\
\cline{3-15}

&& $1\%$& 1.2 & 1.5 & 1.1 & 1.1 & 1.7 & 1.4 & 1.8 & 1.5 & 1.4 & 1.9 & 1.7 & 1.7\\
$(3,2.5)$&$2,000$& $5\%$& 5.1 & 4.9 & 5.9 & 6.0 & 5.7&  6.3 & 6.6 & 5.7 & 5.6 & 5.9 & 5.6  &5.5\\
 && $10\%$& 10.8 &10.1& 10.6 &11.0& 11.3& 10.9 &11.4 &12.0& 12.6 &12.0 &11.0 &10.2\\
\hline\hline 
\multicolumn{14}{l}{ Model \eqref{arch}--\eqref{coef} with $A_{01}^+\neq A_{01}^-$ when $\underline{\delta}_0$ is unknown.}
\\
\end{tabular}
\end{center}
}
\label{archasymInc}
\end{table}

\begin{table}[H]
 \caption{\small{Empirical power of the proposed test for the null hypothesis of a bivariate CCC-APGARCH$(0,1)$ model of the form \eqref{arch} against the bivariate CCC-APGARCH$(1,1)$ alternative given by \eqref{model-puiss} when $\underline{\delta}_0$ is unknown.}}
  {\scriptsize 
\begin{center}
\begin{tabular}{ccc cccccccccccc}
\hline\hline
$\delta_0$& Length $n$ & Level $\alpha$ & \multicolumn{12}{c}{Lag $m$}\\
\hline
&&&1&2&3&4&5&6&7&8&9&10&11&12\\
&& $1\%$&21.6& 19.6& 23.8& 27.8 &23.0 &20.4 &16.6 &13.0& 11.2& 10.0 & 9.8 & 8.2\\
$(1,1)$&$250$& $5\%$& 37.0 &43.6 &54.6& 55.0& 50.0& 46.4& 43.0& 39.0& 34.8& 31.4& 30.8& 27.8\\
 && $10\%$& 47.2 &56.4& 67.2& 69.8 &68.0& 62.4& 58.4 &54.8& 52.2& 49.6& 45.2& 42.2\\
\cline{3-15}
&& $1\%$&46.8& 53.8 &69.2 &72.8 &71.8 &70.2& 65.6& 60.0& 56.2 &53.4 &46.8& 43.4\\
$(1,1)$&$500$& $5\%$& 68.6 &77.0& 88.6 &91.0 &91.0 &88.8 &86.2 &83.4 &82.0& 78.2& 75.4& 73.2\\
 && $10\%$& 76.6& 87.4 &95.2& 95.4& 96.0 &95.4 &93.6 &92.0& 89.4& 88.0& 86.2& 85.8\\
\cline{3-15}
&& $1\%$&96.7 & 99.2 & 99.6 & 99.7 & 99.9 & 99.9 & 99.9 & 99.9 & 99.9 & 99.9 & 99.9 &100.0\\
$(1,1)$&$2,000$& $5\%$&98.2 & 99.8 & 99.9 & 99.9 & 99.9  &99.9  &99.9 & 99.9& 100.0 & 99.9 &100.0 &100.0\\
 && $10\%$& 99.1 &100.0  &99.9  &99.9 &100.0& 100.0& 100.0 & 99.9& 100.0& 100.0 &100.0 &100.0\\
 \hline
&&&1&2&3&4&5&6&7&8&9&10&11&12\\
&& $1\%$& 20.8& 17.4 &20.8& 20.2& 18.2& 13.6 &11.2& 10.2 & 8.0 & 6.8 & 6.0 & 5.0\\
$(0.8,1.5)$&$250$& $5\%$& 38.2 &40.6 &47.8& 48.8& 46.8& 43.0& 37.0& 33.6 &29.4 &24.6 &24.4 &23.2\\
 && $10\%$& 49.6 &56.2& 66.0 &66.6 &63.2 &59.2 &53.8& 51.0& 46.4& 43.4 &42.0& 39.4 \\
\cline{3-15}
&& $1\%$&47.5 &49.2 &61.7 &65.6 &63.2& 60.9& 55.4& 50.1& 46.2& 41.0 &38.7 &32.9\\
$(0.8,1.5)$&$500$& $5\%$& 67.5& 74.5& 84.4& 87.5& 86.0& 83.6 &82.4& 78.7& 75.3 &71.7 &69.1& 65.0\\
 && $10\%$& 75.5 &84.5& 92.1& 94.0 &93.2 &92.4 &90.4 &88.2 &86.8& 83.4 &81.6 &79.5\\
\cline{3-15}
&& $1\%$&92.8& 96.2 &98.1& 98.7 &98.8 &98.7 &98.5 &98.4& 98.2 &97.8 &97.5 &97.8\\
$(0.8,1.5)$&$2,000$& $5\%$& 96.0 &98.3 &99.0 &99.3 &99.4 &99.4 &99.4 &99.3 &99.2& 98.9 &99.1 &99.2\\
 && $10\%$& 97.2 &99.3& 99.3 &99.6 &99.6 &99.5 &99.5 &99.5 &99.5& 99.5& 99.5 &99.5\\
\hline
&&&1&2&3&4&5&6&7&8&9&10&11&12\\
&& $1\%$& 14.4& 17.2& 23.2 &23.0& 22.0& 19.0& 17.8& 15.0& 11.6& 11.2 & 8.2 & 7.2\\
$(2,2)$&$250$& $5\%$& 29.0& 38.8& 56.8& 59.0& 55.2 &50.4 &45.2& 43.2& 38.0& 35.0& 32.0& 29.0\\
 && $10\%$& 37.8& 57.8& 72.2& 72.4& 71.8& 69.2& 65.8 &59.0& 53.8& 52.2& 48.4 &46.4\\
\cline{3-15}
&& $1\%$& 40.4& 51.8 &69.8 &76.6& 76.6& 73.6& 70.6& 65.8 &61.6& 58.2& 52.8 &47.6\\
$(2,2)$&$500$& $5\%$& 59.6& 76.6& 90.0& 92.0 &91.4& 91.8& 90.2 &87.8& 86.0& 82.6 &79.8 &77.6\\
 && $10\%$& 69.4& 86.2& 94.8 &97.0& 95.6& 95.8& 95.4& 94.8& 94.2& 93.0& 91.2& 89.4\\
\cline{3-15}

&& $1\%$& 91.8 &98.2 &99.0 &99.4 &99.7 &99.7 &99.8& 99.8 &99.8 &99.8 &99.6& 99.5\\
$(2,2)$&$2,000$& $5\%$& 96.0 & 99.4  &99.7& 100.0 &100.0 & 99.9 & 99.8 & 99.9 & 99.9  &99.9 & 99.9 & 99.9\\
 && $10\%$& 97.0 & 99.6 & 99.9 &100.0& 100.0 &100.0 &100.0& 100.0 & 99.9 &100.0 &100.0  &99.9 \\
\hline
&&&1&2&3&4&5&6&7&8&9&10&11&12\\
&& $1\%$& 24.6& 24.8 &33.2& 36.8& 36.8 &33.2& 29.6& 26.6 &23.2 &21.8& 19.2& 16.4\\
$(3.,2.5)$&$250$& $5\%$& 41.0& 52.0& 62.0 &64.4 &64.8 &62.4& 57.0& 54.8& 50.8& 46.8& 43.6& 43.2\\
 && $10\%$& 48.8 &64.4& 72.4 &76.2& 75.0& 75.4& 72.8& 69.2& 68.4& 65.2 &62.2 &60.4\\
\cline{3-15}
&& $1\%$& 52.8& 56.6& 67.0& 71.6 &72.4 &70.4& 67.4& 65.8& 63.4 &59.2 &56.6& 55.6\\
$(3.,2.5)$&$500$& $5\%$& 66.6& 77.8& 85.8& 87.2& 88.6& 89.4& 88.0& 88.2& 87.4& 84.6& 81.6& 80.0\\
 && $10\%$& 74.6& 84.6& 92.0 &92.0& 93.6& 93.6& 93.4& 92.8& 91.4& 91.0& 91.2& 90.4\\
\cline{3-15}

&& $1\%$& 90.4 &93.7& 95.6 &96.7 &97.7 &97.6& 97.8 &98.0 &97.7& 98.0 &97.8 &97.9\\
$(3.,2.5)$&$2,000$& $5\%$& 94.4 &96.9& 98.3 &98.8 &98.9& 99.0& 99.2 &99.3& 99.1& 98.8 &99.0 &99.0\\
 && $10\%$& 96.4 &98.3 &98.9 &99.3 &99.3 &99.4& 99.5& 99.6 &99.3& 99.3 &99.5 &99.5\\
\hline\hline 
\multicolumn{14}{l}{ Model \eqref{model-puiss} with $A_{01}^+\neq A_{01}^-$ when $\underline{\delta}_0$ is unknown.}
\\\end{tabular}
\end{center}
}
\label{puiss-archasymInc}
\end{table}

\section{Appendix : Proofs of the mains results}\label{preuves}
To prove the main results we need some tools from \cite{Nous-MAPGARCH} summarized in the following lemma.
\subsection{Preliminaries}
%\subsubsection{When the power $\underline{\delta}_0$ is known }
For all $\theta \in \Theta$, recall that $\underline{\tilde{h}}^{\underline{\delta}_0/2}_t(\theta)$ is the strictly stationary and non-anticipative solution of \eqref{Ht-connu}.
\begin{lem}\  {(\cite{Nous-MAPGARCH})}

Under Assumptions \textbf{A1}-\textbf{A7} and for $s\in ]0, 1[$, we have 
\begin{equation}\label{cor2}
\mathbb{E}\Vert \underline{\varepsilon}_t^{\underline{\delta}_0/2} \Vert^s < \infty, \qquad \mathbb{E}\sup\limits_{\theta\in \Theta} \Vert \underline{h}_{0t}^{\underline{\delta}_0/2} \Vert^s < \infty, \qquad \mathbb{E}\sup\limits_{\theta\in \Theta} \Vert \underline{\tilde{h}}_{0t}^{\underline{\delta}_0/2} \Vert^s < \infty.
\end{equation}
Moreover, there exists $K$ a random constant that depends on the past values of $\{\varepsilon_t, t\le 0\}$ and $0<\rho<1$ such that
\begin{equation}\label{Majh.delta}
\begin{aligned}
\sup\limits_{\theta\in\Theta} \left\Vert \underline{ h}_{t}^{\underline{\delta}_0/2}(\theta) - \underline{\tilde h}_{t}^{\underline{\delta}_0/2}(\theta)\right\Vert
\leq K\rho^t \ .
\end{aligned}
\end{equation}
Thus, for $i_1=1,\dots,d$, since $\min\left( { h}_{i_1,t}^{\delta_{0,i_1}/2}(\theta),{\tilde h}_{i_1,t}^{{\delta_{0,i_1}/2}}(\theta)\right)\geq\omega=\inf\limits_{1\leq i\leq d}\underline{\omega}(i)$, the mean-value theorem implies that
\begin{align}\label{MajH2}
\sup\limits_{\theta\in\Theta} \left| { h}_{i_1,t}(\theta)-{\tilde h}_{i_1,t}(\theta)\right| &\leq\frac{2}{\delta_{0,i_1}} \sup\limits_{\theta\in\Theta}\max\left( { h}_{i_1,t}^{1-\delta_{0,i_1}/2}(\theta),{\tilde h}_{i_1,t}^{1-{\delta_{0,i_1}/2}}(\theta)\right) \sup\limits_{\theta\in\Theta}
\left| { h}_{i_1,t}^{\delta_{0,i_1}/2}(\theta)-{\tilde h}_{i_1,t}^{{\delta_0,{i_1}/2}}(\theta)\right| \nonumber \\
&\leq\frac{2K}{\delta_{0,i_1}} \left(\sup\limits_{\theta\in\Theta}\frac{1}{\omega}\right)\sup\limits_{\theta\in\Theta}\max\left( { h}_{i_1,t}(\theta),{\tilde h}_{i_1,t}(\theta)\right) \rho^t
\leq K\rho^t \ ,
\end{align}
and similarly
\begin{align}\label{MajH2.1}
\sup\limits_{\theta\in\Theta} \left| { h}_{i_1,t}^{1/2}(\theta)-{\tilde h}_{i_1,t}^{1/2}(\theta)\right| &\leq\frac{1}{\delta_{0,i_1}} \sup\limits_{\theta\in\Theta}\max\left( { h}_{i_1,t}^{(1-\delta_{0,i_1})/2}(\theta),{\tilde h}_{i_1,t}^{(1-{\delta_{0,i_1}})/2}(\theta)\right) \sup\limits_{\theta\in\Theta}
\left| { h}_{i_1,t}^{\delta_{0,i_1}/2}(\theta)-{\tilde h}_{i_1,t}^{{\delta_{0,i_1}/2}}(\theta)\right|  \nonumber \\
&\leq\frac{K}{\delta_{0,i_1}} \left(\sup\limits_{\theta\in\Theta}\frac{1}{\omega}\right)\sup\limits_{\theta\in\Theta}\max\left( { h}_{i_1,t}^{1/2}(\theta),{\tilde h}_{i_1,t}^{1/2}(\theta)\right) \rho^t
\leq K\rho^t \ .
\end{align}
From \eqref{MajH2} we can deduce that, almost surely, we have
\begin{equation}\label{MajH-H1}
\begin{aligned}
\sup\limits_{\theta\in\Theta} \left\Vert H_t({\theta}) - \tilde{H}_t({\theta})\right\Vert \leq K\rho^t,\qquad \forall t.
\end{aligned}
\end{equation}
Since $\Vert R^{-1} \Vert$ is the inverse of the eigenvalue of smaller module of $R$ and $\Vert \tilde{D}_t^{-1}\Vert = [\min_i(h_{ii,t}^{1/2})]^{-1}$ for $i=1,\dots,d$, we have
\begin{align}\label{MajDR1}
\sup\limits_{\theta\in\Theta} \Vert \tilde{H}_t^{-1}({\theta})\Vert&  \leq \sup\limits_{\theta\in\Theta}\Vert \tilde{D}_t^{-1}\Vert^2\Vert R^{-1}\Vert \leq \sup\limits_{\theta\in\Theta}\left [\min_i(\underline{\omega}(i))\right ]^{-1}\Vert R^{-1}\Vert \leq K,
\end{align}
by using the fact that $R$ is a positive-definite matrix (see Assumption \textbf{A5}), the compactness of $\Theta$ and the strict positivity of the components of $\underline{\omega}$. Similarly, we have 
\begin{equation}\label{MajH1.1}
\sup\limits_{\theta\in\Theta} \Vert H_t^{-1}({\theta})\Vert \leq K.
\end{equation}
There exists a neighborhood $V(\theta_0)$ of $\theta_0\in\overset{\circ}{\Theta}$ such that: for all $r_0\geq 1$, $i_1=1,\dots,d$ and all $i,j=1,\dots,s_1$ we have
\begin{align}
\label{unif1}
\mathbb{E}\sup_{\theta\in V(\theta_0)}\left\vert \dfrac{1}{{h}_{i_1,t}^{{\delta}_{0,i_1}/2}} \dfrac{\partial {h}_{i_1,t}^{{\delta}_{0,i_1}/2}}{\partial\theta_i}(\theta)\right\vert^{r_0}  & < \infty\qquad\text{and}\qquad
\mathbb{E}\sup_{\theta\in V(\theta_0)}\left\vert \dfrac{1}{{h}_{i_1,t}^{{\delta}_{0,i_1}/2}} \dfrac{\partial^2 {h}_{i_1,t}^{{\delta}_{0,i_1}/2}}{\partial\theta_i\partial\theta_j}(\theta)\right\vert^{r_0} &< \infty \ .
\end{align}

%\begin{equation}
%\mathbb{E} \sup\limits_{\theta\in V(\theta_0)}\left\vert \dfrac{h_{i_1,t}^{\delta_{i_1}/2}(\theta_0)}{h_{i_1,t}^{\delta_{i_1}/2}(\theta)}\right\vert  < \infty. 
%\end{equation}
%The matrix $J$ is reversal and 
%\begin{equation}
%\sqrt{n}(\hat{\theta}_n - \theta_0) = - J^{-1}\dfrac{1}{\sqrt{n}}\sum\limits_{t=1}^n\dfrac{\partial l_t(\theta_0)}{\partial \theta_i} + o_{\mathbb{P}}(1).
%\end{equation}
\end{lem}
In the case where the power is unknown, the vector of parameter becomes $\vartheta$ and we replace $H_t$ by $\mathcal{H}_t$ and $\theta$ by $\vartheta$. The previous results must be adapted in consequence. 
\subsection{Proof of Theorem \ref{thmD-connu}}
We decomposed this proof in following steps.

\begin{enumerate}[$\qquad (i)$]
	\item Asymptotic impact of the unknown initial values on the statistic ${\mathbf{\hat{r}}}_m$.
	\item Asymptotic distribution of $\sqrt{n}{\mathbf{\hat{r}}}_m$.
	\item Invertibility of the matrix $D$.
\end{enumerate}
Recall that
\begin{eqnarray*}
 {r}_h(\theta)& = &\dfrac{1}{n} \sum\limits_{t = h +1}^n [\text{Tr}({s}_t(\theta))][\text{Tr}({s}_{t-h}(\theta))]\qquad \text{with } {s}_t(\theta) = {\eta}_t(\theta){\eta}_t'(\theta) - I_d.
\\&=& \dfrac{1}{n} \sum\limits_{t = h +1}^n [{\eta}_t'(\theta){\eta}_t(\theta)-d][{\eta}_{t-h}'(\theta){\eta}_{t-h}(\theta)-d]
\\&=& \dfrac{1}{n} \sum\limits_{t = h +1}^n [\underline{\varepsilon}_t'{H}_t^{-1}({\theta})\underline{\varepsilon}_t-d][\underline{\varepsilon}_{t-h}'{H}_{t-h}^{-1}({\theta})\underline{\varepsilon}_{t-h}-d].
\end{eqnarray*}

\textbf{$(i)$ Asymptotic impact of the unknown initial values on the statistic ${\mathbf{\hat{r}}}_m$}

Let ${S}_t(\theta) =\underline{\varepsilon}_t'{H}_t^{-1}({\theta})\underline{\varepsilon}_t-d$ and $\tilde{S}_t(\theta) =\underline{\varepsilon}_t'\tilde{H}_t^{-1}({\theta})\underline{\varepsilon}_t-d$.
We observe that 
\begin{eqnarray*}
 {r}_h(\theta)-\tilde{r}_h(\theta)& = &\dfrac{1}{n} \sum\limits_{t = h +1}^n  (a_t + b_t),
\end{eqnarray*}
where $a_t = S_{t-h}(\theta)(S_t(\theta) - \tilde{S}_t(\theta))$ and $b_t = (\tilde{S}_{t-h}(\theta)-S_{t-h}(\theta)) \tilde{S}_t(\theta)$.
We obtain
\begin{equation*}
\begin{aligned}
|a_t| &= \left\vert \underline{\varepsilon}_{t-h}'{H}_{t-h}^{-1}({\theta})\underline{\varepsilon}_{t-h}-d \right\vert \left\vert \underline{\varepsilon}_t' H_t^{-1}({\theta})(H_t({\theta}) - \tilde{H}_{t}({\theta}))\tilde{H}_t^{-1}({\theta})\underline{\varepsilon}_t \right\vert \\ %
 &=  \left\vert\text{Tr}({H}_{t-h}^{-1}({\theta})\underline{\varepsilon}_{t-h}\underline{\varepsilon}_{t-h}'-I_d) \right\vert  \left\vert \text{Tr}(H_t^{-1}({\theta})(H_t({\theta}) - \tilde{H}_{t}({\theta}))\tilde{H}_t^{-1}({\theta})\underline{\varepsilon}_t \underline{\varepsilon}_t')\right\vert \\
 &\leq\sup\limits_{\theta\in\Theta}\left(\Vert H_t^{-1}({\theta})\Vert\Vert\underline{\varepsilon}_{t-h} \underline{\varepsilon}_{t-h}'\Vert+\Vert I_d\Vert\right)
 \left(\Vert H_t^{-1}({\theta})\Vert\Vert H_t({\theta}) - \tilde{H}_{t}({\theta})\Vert\Vert\tilde{H}_t^{-1}({\theta})\Vert\Vert\underline{\varepsilon}_t \underline{\varepsilon}_t'\Vert\right). 
\end{aligned}
\end{equation*}
Now using \eqref{MajH-H1}, \eqref{MajDR1} and \eqref{MajH1.1}, we have
\begin{align*}
|a_t|
& \leq K\rho^t(\underline{\varepsilon}_{t-h}'\underline{\varepsilon}_{t-h} + d)\underline{\varepsilon}_t'\underline{\varepsilon}_t
\end{align*}
We have the same bound  for $|b_t|$.  
Using the inequality $(a+b)^s \leq a^s + b^s$, for $a,b \geq 0$ and $s\in ]0,1[$, \eqref{cor2} and H{\"o}lder's inequality, 
we have for some $s^\ast \in ]0,1[$ sufficiently small
\begin{align*}
\mathbb{E}\left\vert \dfrac{1}{\sqrt{n}} \sum\limits_{t=1}^n \sup\limits_{\theta\in\Theta}\vert a_t \vert \right\vert^{s\ast} &\leq \mathbb{E}\left\vert \dfrac{1}{\sqrt{n}}\sum\limits_{t=1}^n \sup\limits_{\theta\in\Theta}\Vert K\rho^t(\underline{\varepsilon}_{t-h}\underline{\varepsilon}_{t-h}' + I_d)\underline{\varepsilon}_t\underline{\varepsilon}_t'\Vert \right\vert^{s\ast}\\
%&\leq K\left(\dfrac{1}{\sqrt{n}}\right)^{s\ast} \mathbb{E}\left\Vert \sum\limits_{t=1}^n\rho^t(\underline{\varepsilon}_{t-h}\underline{\varepsilon}_{t-h}' + I_d)\underline{\varepsilon}_t\underline{\varepsilon}_t'\Vert \right\Vert^{s\ast}\\
%&\leq K\left(\dfrac{1}{\sqrt{n}}\right)^{s\ast} \mathbb{E}\left[ \sum\limits_{t=1}^n\rho^{ts\ast}\right]\mathbb{E}\left[(\underline{\varepsilon}_{t-h}\underline{\varepsilon}_{t-h}' + I_d)\underline{\varepsilon}_t\underline{\varepsilon}_t'^{s\ast}\right]\\
%&\leq K\left(\dfrac{1}{\sqrt{n}}\right)^{s\ast} \sum\limits_{t=1}^n\rho^{ts\ast}\mathbb{E}\left[(\underline{\varepsilon}_{t-h}\underline{\varepsilon}_{t-h}' + I_d)^{s\ast}\right]\mathbb{E}\left[\underline{\varepsilon}_t\underline{\varepsilon}_t'^{s\ast} \right]\\
%&\leq K\left(\dfrac{1}{\sqrt{n}}\right)^{s\ast} \sum\limits_{t=1}^n\rho^{ts\ast}\mathbb{E}\left[(\underline{\varepsilon}_{t-h}\underline{\varepsilon}_{t-h}'^{s\ast} + I_d^{s\ast})\right]\mathbb{E}\left[\underline{\varepsilon}_t\underline{\varepsilon}_t'^{s\ast} \right]\\
&\leq K\left(\dfrac{1}{\sqrt{n}}\right)^{s\ast} \sum\limits_{t=1}^n\rho^{ts\ast} \underset{n\to\infty}{\longrightarrow} 0
\end{align*}
We deduce 
\begin{equation*}
\dfrac{1}{\sqrt{n}} \sum\limits_{t=1}^n \sup\limits_{\theta \in \Theta} \vert a_t \vert = 
\mathrm{o}_{\mathbb{P}}(1).
\end{equation*}
We have the same convergence for $b_t$, and for the derivatives of $a_t$ and $b_t$.
Consequently we obtain
\begin{equation}\label{Impact_as_connu}
\sqrt{n}\Vert {\mathbf{r}}_m(\theta_0) - {\mathbf{\tilde{r}}}_m(\theta_0) \Vert = \mathrm{o}_{\mathbb{P}}(1), \qquad \sup\limits_{\theta\in \Theta}\left\Vert \dfrac{\partial {\mathbf{r}}(\theta)}{\partial \theta'} - \dfrac{\partial {\mathbf{\tilde{r}}}(\theta)}{\partial \theta'}\right\Vert = \mathrm{o}_{\mathbb{P}}(1).
\end{equation}
The unknown initial values have no asymptotic impact on the statistic ${\mathbf{\hat{r}}}_m$.\\

\textbf{$(ii)$ Asymptotic distribution of $\sqrt{n}{\mathbf{\hat{r}}}_m$}

We now show that the asymptotic distribution of $\sqrt{n}{\mathbf{\hat{r}}}_m$ is deduced from the joint distribution of $\sqrt{n}{\mathbf{r}}_m$ and the QMLE.

Using \eqref{Impact_as_connu} and a Taylor expansion of ${\mathbf{r}}_m(\cdot)$ around $\hat{\theta}_n$ and $\theta_0$, we obtain
\begin{align*}
\sqrt{n}{\mathbf{\hat{r}}}_m &= \sqrt{n}{\mathbf{\tilde{r}}}_m(\theta_0) + \dfrac{\partial {\mathbf{\tilde{r}}}_m(\theta^\ast)}{\partial \theta'}\sqrt{n}(\hat{\theta}_n - \theta_0)\\
&= \sqrt{n}{\mathbf{r}}_m(\theta_0) + \dfrac{\partial {\mathbf{{r}}}_m(\theta^\ast)}{\partial \theta'}\sqrt{n}(\hat{\theta}_n - \theta_0) + \mathrm{o}_{\mathbb{P}}(1),
\end{align*}
for some $\theta^\ast$ between $\hat{\theta}_n$ and $\theta_0$.
For $i,j=1,\dots,s_0$ the first and the second  derivatives of $S_t(\theta)$ give
%
%\[\dfrac{\partial S_t(\theta)}{\partial \theta_i} = -\text{Tr}\left[H_t^{-1}(\theta)\underline{\varepsilon}_t\underline{\varepsilon}_t'H_t^{-1}(\theta)\dfrac{\partial H_t(\theta)}{\partial\theta_i}\right],\]
%%
%and the second derivative
\begin{equation*}
\begin{aligned}
\dfrac{\partial S_t(\theta)}{\partial \theta_i}& = -\text{Tr}\left[H_t^{-1}(\theta)\underline{\varepsilon}_t\underline{\varepsilon}_t'H_t^{-1}(\theta)\dfrac{\partial H_t(\theta)}{\partial\theta_i}\right],
\\
\dfrac{\partial^2 S_t(\theta)}{\partial \theta_i\partial\theta_j} &= \text{Tr} \left[H_t^{-1}(\theta)\dfrac{\partial H_t(\theta)}{\partial\theta_j}H_t^{-1}(\theta)\underline{\varepsilon}_t\underline{\varepsilon}_t'H_t^{-1}(\theta)\dfrac{\partial H_t(\theta)}{\partial \theta_i} -H_t^{-1}(\theta)\underline{\varepsilon}_t\underline{\varepsilon}_t'H_t^{-1}(\theta)\dfrac{\partial^2H_t(\theta)}{\partial\theta_i\partial\theta_j} \right.\\
&\qquad \left. \qquad\qquad+ H_t^{-1}(\theta)\underline{\varepsilon}_t\underline{\varepsilon}_t'H_t^{-1}(\theta)\dfrac{\partial H_t(\theta)}{\partial \theta_j}H_t^{-1}(\theta)\dfrac{\partial H_t(\theta)}{\partial\theta_i}  \right].
\end{aligned}
\end{equation*}
In view of \eqref{unif1}, there exists a neighborhood ${V}(\theta_0)$ of $\theta_0$ such that
\begin{equation*}
\mathbb{E} \sup\limits_{\theta\in{V}(\theta_0)}\left\| \dfrac{\partial^2 S_{t-h}(\theta)S_t(\theta)}{\partial \theta\partial \theta'}\right\| < \infty.
\end{equation*}
%Using a second Taylor expansion, we obtain 
For $i=1,\dots,s_0$, let $\mathbf{h}_t(i) = \left[\text{vec}\left(H_{0t}^{-1/2}({\partial H_t(\theta_0)}/{\partial \theta_i})H_{0t}^{-1/2}\right)\right]$ and we define the matrix of size $d^2\times s_0$,  $\mathbf{h}_t = (\mathbf{h}_t(1) \vert \ldots \vert \mathbf{h}_t({s_0}))$.
For a fixed $r_h$, using the previous  inequality, Assumption \textbf{A7},
the almost sure convergence of $\theta^\ast$ to $\theta_0$, a second Taylor expansion and the ergodic theorem, we obtain
\begin{align*}
\dfrac{\partial {{r}}_h(\theta^\ast)}{\partial \theta_i} = \dfrac{\partial 
{{r}}_h(\theta_0)}{\partial \theta_i}+ \mathrm{o}_{\mathbb P}(1) \underset{n\to \infty}{\longrightarrow}C(h,i) &:= \mathbb{E}\left[S_{t-h} \dfrac{\partial S_t}{\partial \theta_i}\right] = -\mathbb{E}\left[S_{t-h}
\text{Tr}\left(H_{0t}^{-1}\dfrac{\partial H_t(\theta_0)}{\partial\theta_i}\right)\right]
\\&
= - \mathbb{E}\left[S_{t-h}\mathbf{h}'_t(i)\text{vec}(I_d)\right],
\end{align*}
by the fact $\mathbb{E}[S_t\partial S_{t-h}(\theta_0)/\partial \theta] = 0$ and using the property
$\text{Tr}(A'B)=(\text{vec}(A))'\text{vec}(B)$. 
Note that $C(h,i)$ is the $(h,i)$-th element of the matrix $C_m$. Consequently we have
\begin{equation}\label{Cm}
\dfrac{\partial {\mathbf{r}}_m(\theta_0)}{\partial \theta'} \underset{n\to\infty}{\longrightarrow} C_m := [C(h,i)]_{1\leq h\leq m,1\leq i\leq s_0}
= - \mathbb{E}\left[\left(\mathbb{S}_{t-1:t-m}\right)\left(\mathbf{h}'_t\text{vec}(I_d)\right)'\right],
\end{equation}
where $\mathbb{S}_{t-1:t-m} = (S_{t-1},\ldots, S_{t-m})'$. It follows that
\begin{equation}\label{DLTrm}
\sqrt{n}{\mathbf{\hat{r}}}_m = \sqrt{n}{\mathbf{r}}_m + C_m\sqrt{n}(\hat{\theta}_n - \theta_0) + \mathrm{o}_{\mathbb{P}}(1).
\end{equation}
From \eqref{DLTrm} it is clear that the asymptotic distribution of $\sqrt{n}{\mathbf{\hat{r}}}_m$ is related to the asymptotic behavior of $\sqrt{n}( \hat{\theta}'_n - \theta'_0,{\mathbf{r}}'_m)'$.
We note that
\begin{equation*}
\sqrt{n}(\hat{\theta}_n - \theta_0) = -J^{-1}\left(\dfrac{1}{\sqrt{n}}\sum\limits_{t=1}^n\dfrac{\partial l_t(\theta_0)}{\partial \theta}\right),
\end{equation*}
with $l_t(\theta) = \underline{\varepsilon}_t'H_t^{-1}(\theta)\underline{\varepsilon}_t + \log(\det(H_t(\theta)))$. The derivatives are recursively calculated with respect to $H_t(\theta)$ for a fixed $i=1,\dots,s_0$ 
\begin{align*}
\dfrac{\partial l_t(\theta_0)}{\partial \theta_i}= \text{Tr}\left[(H_{0t}^{-1} - H_{0t}^{-1}\underline{\varepsilon}_t\underline{\varepsilon}_t'H_{0t}^{-1})\dfrac{\partial H_t(\theta_0)}{\partial\theta_i}\right]
 = - \left[\text{vec}\left(H_{0t}^{-1/2}\dfrac{\partial H_t(\theta_0)}{\partial \theta_i}H_{0t}^{-1/2}\right)\right]'\text{vec}(s_t) = - \mathbf{h}'_t(i)\text{vec}(s_t).
\end{align*}
%with $\mathbf{h}_t(i) = \left[\text{vec}\left(H_{0t}^{-1/2}\dfrac{\partial H_t(\theta_0)}{\partial \theta_i}H_{0t}^{-1/2}\right)\right]$. 
%We can also rewrite $\mathbf{h}_t(i)$ as follow, using the relation $\text{vec}(ABC) = (C'\otimes A)\text{vec}(B)$
%\begin{align*}
%\mathbf{h}_t(i) = \left[\text{vec}\left(H_{0t}^{-1/2}\dfrac{\partial H_t(\theta_0)}{\partial \theta_i}H_{0t}^{-1/2}\right)\right] = \left[H_{0t}^{-1/2}\otimes H_{0t}^{-1/2}\right] \text{vec}\left(\dfrac{\partial H_t(\theta_0)}{\partial \theta_i}\right) = \mathbf{H}_t\mathbf{d}_t(i),
%\end{align*}
%where $\mathbf{H}_t =\left[H_{0t}^{-1/2}\otimes H_{0t}^{-1/2}\right] $ and $\mathbf{d}_t(i) = \text{vec}\left({\partial H_t(\theta_0)}/{\partial \theta_i}\right)$. 
%%Finally, we define the matrices of size $d^2\times s_0$, $\mathbf{h}_t = (\mathbf{h}_t(1) \vert \ldots \vert \mathbf{h}_t({s_0}))$ and $\mathbf{d}_t = (\mathbf{d}_t(1) \vert \ldots \vert \mathbf{d}_t({s_0}))$ such that $\mathbf{h}_t = \mathbf{H}_t\mathbf{d}_t$.
%%%
%We define the derivative by 
We then deduce that
\begin{equation*}
\dfrac{\partial l_t(\theta_0)}{\partial \theta} = \left(\dfrac{\partial l_t(\theta_0)}{\partial \theta_1}, \ldots, \dfrac{\partial l_t(\theta_0)}{\partial \theta_{s_0}}\right)'= - \mathbf{h}'_t\text{vec}(s_t).
\end{equation*}
Observe that $\sqrt{n}\underline{\mathbf{r}}_m = n^{-1} \sum\limits_{t=1}^n\mathbb{S}_{t-1:t-m} S_t$. Now we can obtained the asymptotic distribution of $\sqrt{n}(\hat{\theta}'_n - \theta'_0, {\mathbf{r}}'_m)'$ by applying the central limit theorem to the multivariate martingale difference
\begin{equation*}
\left\{\Upsilon_t = \left(\left\{ J^{-1}\mathbf{h}'_t\text{vec}(s_t)\right\}', \mathbb{S}_{t-1:t-m} S_t\right)'; \mathcal{F}_{t-1}^\eta := \sigma(\eta_u, u\leq t )\right\}.
\end{equation*}
The expectation of the distribution is given by
\begin{equation*}
\mathbb{E}\left[\Upsilon_t \vert \mathcal{F}_{t-1}^\eta\right] = \mathbb{E}\left[\left.\begin{pmatrix} J^{-1}\mathbf{h}'_t\text{vec}(s_t) \\ 
S_{t-1}  S_t \\ \vdots \\ S_{t-m}  S_t \end{pmatrix} \right\vert \mathcal{F}_{t-1}^\eta\right] = \begin{pmatrix} J^{-1}\mathbf{h}'_t \mathbb{E}[\text{vec}(s_t)\vert \mathcal{F}_{t-1}^\eta] \\ S_{t-1}  \mathbb{E}[S_t\vert \mathcal{F}_{t-1}^\eta] \\ \vdots \\ S_{t-m}  \mathbb{E}[S_t\vert \mathcal{F}_{t-1}^\eta] \end{pmatrix}  =0,
\end{equation*}
because $(S_{t-i})_{i\geq 1}$ is measurable with respect to the $\sigma$-field $\mathcal{F}_{t-1}^\eta$ and $\mathbb{E}[S_t]=\mathbb{E}[\eta_t'\eta_t]-d = 0$ and $\mathbb{E}[\text{vec}(s_t)]=\text{vec}[\mathbb{E}(\eta_t\eta'_t)-I_d] = 0$. For $i=1,\dots,d$, the variance is given by 
\begin{align}\nonumber
\Xi &:= \mathbb{E}\left[ \Upsilon_t \Upsilon_t'\right]= \begin{pmatrix} \Sigma_{\hat{\theta}_n} & \Sigma_{\hat{\theta}_n,{\mathbf{r}}_m} \\ \Sigma_{\hat{\theta}_n,{\mathbf{r}}_m}' & \Sigma_{{\mathbf{r}}_m} \end{pmatrix}\\& = \begin{pmatrix} J^{-1}IJ^{-1} & \mathbb{E}\left[J^{-1}\mathbf{h}'_t\text{vec}(s_t)S_t\mathbb{S}_{t-1:t-m}'\right] \\ \mathbb{E}\left[\mathbb{S}_{t-1:t-m}\left(\text{vec}(s_t)\right)' S_t\mathbf{h}_tJ^{-1}\right] &\left(\mathbb{E}\left[S_{t}^2\right]\right)^2I_m\end{pmatrix},\label{Loi_as_connue}
\end{align}
which leads to 
\begin{equation*}
\dfrac{1}{\sqrt{n}} \sum\limits_{t=1}^n \Upsilon_t \xrightarrow[n\to\infty]{\mathrm{d}} \mathcal{N}(0, \Xi).
\end{equation*}
Using \eqref{Impact_as_connu} and \eqref{Loi_as_connue}, the asymptotic distribution of $\sqrt{n}{\mathbf{\hat{r}}}_m$ gives 
\begin{equation*}
\sqrt{n} {\mathbf{\hat{r}}}_m \xrightarrow[n\to\infty]{\mathrm{d}} \mathcal{N}(0, D),
\end{equation*}
where $D$ is a matrix defined as follows 
%\small{
\begin{align*}
D := \lim\limits_{n \to \infty} \text{Var}(\sqrt{n}{\mathbf{\hat{r}}}_m) &= \lim\limits_{n \to\infty} \text{Var}(\sqrt{n}{\mathbf{r}}_m) + C_m \left[\lim\limits_{n\to\infty} \text{Var}(\sqrt{n}(\hat{\theta}_n - \theta_0))\right]C_m' \\ 
& \quad + C_m\left[\lim\limits_{n \to\infty} \text{Cov}(\sqrt{n}(\hat{\theta}_n - \theta_0),\sqrt{n}{\mathbf{r}}_m)\right] \\ 
& \quad  + \left[\lim\limits_{n \to\infty} \text{Cov}(\sqrt{n}(\hat{\theta}_n - \theta_0),\sqrt{n}{\mathbf{r}}_m)\right]C_m'\\
& = \Sigma_{{\mathbf{r}}_m} + C_mJ^{-1}IJ^{-1}C_m' + C_m\Sigma_{\hat{\theta}_n,{\mathbf{r}}_m}  + \Sigma_{\hat{\theta}_n,{\mathbf{r}}_m}'C_m'.
%\\&=d^2\left(\mathbb{E}\left[\eta_{it}^4\right]-1\right)^2I_m+C_m\left(J^{-1}IJ^{-1}  -2\left(\mathbb{E}\left[\eta_{it}^4\right]-1\right)J^{-1}\right)C'_m.
\end{align*}
%}
%\normalsize

\textbf{$(iii)$ Invertibility of the matrix $D$}

%with $\mathbf{h}_t(i) = \left[\text{vec}\left(H_{0t}^{-1/2}\dfrac{\partial H_t(\theta_0)}{\partial \theta_i}H_{0t}^{-1/2}\right)\right]$. 
Note that by using the relation $\text{vec}(ABC) = (C'\otimes A)\text{vec}(B)$ we can also rewrite $\mathbf{h}_t(i)$ as follows
\begin{align*}
\mathbf{h}_t(i) = \left[\text{vec}\left(H_{0t}^{-1/2}\dfrac{\partial H_t(\theta_0)}{\partial \theta_i}H_{0t}^{-1/2}\right)\right] = \left[H_{0t}^{-1/2}\otimes H_{0t}^{-1/2}\right] \text{vec}\left(\dfrac{\partial H_t(\theta_0)}{\partial \theta_i}\right) = \mathbf{H}_t\mathbf{d}_t(i),
\end{align*}
where $\mathbf{H}_t =\left[H_{0t}^{-1/2}\otimes H_{0t}^{-1/2}\right] $ and $\mathbf{d}_t(i) = \text{vec}\left({\partial H_t(\theta_0)}/{\partial \theta_i}\right)$. 
Thus we define the matrix of size $d^2\times s_0$,  $\mathbf{d}_t = (\mathbf{d}_t(1) \vert \ldots \vert \mathbf{d}_t({s_0}))$ such that $\mathbf{h}_t = \mathbf{H}_t\mathbf{d}_t$.

To study the invertibility of the matrix $D$ we write $V = \mathbb{S}_{t-1:t-m}  S_t - C_mJ^{-1}\dfrac{\partial l_t(\theta_0)}{\partial \theta}$ such that
\begin{align*}
\mathbb{E}[VV'] &= \mathbb{E}\left[(S_t^2)\mathbb{S}_{t-1:t-m} \mathbb{S}'_{t-1:t-m} \right] - \mathbb{E}\left[( S_t)\mathbb{S}_{t-1:t-m}\dfrac{\partial l_t(\theta_0)}{\partial \theta}\right]J^{-1}C_m' \\
&\qquad  - C_mJ^{-1}\mathbb{E} \left[\dfrac{\partial l_t(\theta_0)}{\partial \theta}(\mathbb{S}_{t-1:t-m}' S_t)\right] + C_mJ^{-1}\mathbb{E}\left[\dfrac{\partial l_t(\theta_0)}{\partial \theta}\dfrac{\partial l_t(\theta_0)}{\partial \theta'}\right]J^{-1}C_m'\\
&= \Sigma_{{\mathbf{r}}_m} + \Sigma_{\hat{\theta}_n,{\mathbf{r}}_m}'C_m' + C_m\Sigma_{\hat{\theta}_n,{\mathbf{r}}_m} + C_mJ^{-1}IJ^{-1}C_m.
%\\&=d^2\left(\mathbb{E}\left[\eta_{it}^4\right]-1\right)^2I_m+C_m\left(J^{-1}IJ^{-1}  -2\left(\mathbb{E}\left[\eta_{it}^4\right]-1\right)J^{-1}\right)C'_m.
\end{align*}
We can rewrite the vector $V$ as
\begin{equation*}
V = \mathbb{S}_{t-1:t-m}  S_t + C_mJ^{-1}\mathbf{d}_t'\mathbf{H}_t'\text{vec}(s_t).
\end{equation*}
If the matrix $\mathbb{E}[VV']$ is singular, then there exist a vector $\lambda = (\lambda_1,\ldots, \lambda_{m})'$ not equal to zero such that 
\begin{equation}\label{lambdaV}
\lambda'V = \lambda'\mathbb{S}_{t-1:t-m}  S_t + \mu  \mathbf{d}_t'\mathbf{H}_t'\text{vec}(s_t) = 0, \quad \mbox{a.s.},
\end{equation}
with $\mu = \lambda'C_mJ^{-1}$. We have $\mu \neq 0$, else $\lambda'\mathbb{S}_{t-1:t-m} S_t = 0$ almost surely, that implies there exists $j\in \{1,\ldots, m\}$ such that $S_{t-j}$ be mesurable respect to the $\sigma$-field $\{S_r, t-1\leq r\leq t-m\}$ with $r\neq t-j$. That is impossible because the $S_t$ are independent and not degenerated. 
Consequently \eqref{lambdaV} becomes 
\begin{align}\label{Raisonnement}
\mu'\mathbf{d}_t' &= \sum\limits_{i=1}^{s_0} \mu_i \mathbf{d}_t(i) = \sum\limits_{i=1}^{s_0}\mu_i    \text{vec}\left(\dfrac{\partial H_t(\theta_0)}{\partial \theta_i}\right) = \sum\limits_{i=1}^{s_0}\mu_i\dfrac{\partial}{\partial \theta_i}\left[(D_{0t}\otimes D_{0t})\text{vec}(R_0)\right]=0\quad \mbox{a.s.},\nonumber\\
&=\sum\limits_{i=1}^{s_1}\mu_i\dfrac{\partial (D_{0t}\otimes D_{0t})}{\partial \theta_i}\text{vec}(R_0) + \sum\limits_{i = s_1 + 1}^{s_0}\mu_i (D_{0t} \otimes D_{0t})\dfrac{\partial \text{vec}(R_0)}{\partial \theta_i}=0\quad \mbox{a.s.}.
\end{align}
Since the vectors $\partial \text{vec}(R_0)/\partial \theta_i$, $i = s_1+1,\ldots, s_0$ are linearly independent, the vector  $(\mu_{s_1+1}, \ldots, \mu_{s_0})'$ is null and thus Equation \eqref{Raisonnement} yields
\begin{equation}\label{DD}
\sum\limits_{i=1}^{s_1} \mu_i \dfrac{\partial (D_{0t}\otimes D_{0t})}{\partial \theta_i}\text{vec}(R_0)=0,\quad \mbox{a.s.}.
\end{equation}
The rows $1,d+1,\ldots,d^2$ of the Equation \eqref{DD} yield
\begin{align}\label{derht}
\sum\limits_{i=1}^{s_1}\mu_i\dfrac{\partial\underline{h}_t(\theta_0)}{\partial \theta_i} = 0,\quad \text{a.s.}
\end{align}
We have for $i_1 = 1,\ldots, d$ and $i=1,\ldots,s_1$
\begin{align}\label{derD}
\dfrac{\partial h_{i_1,t}(\theta_0)}{\partial  \theta_{i}} =\dfrac{\partial \left(h_{i_1,t}^{\delta_{0,i_1}/2}\right)^{2/\delta_{0,i_1}}}{\partial \theta_{i}}(\theta_0) = \dfrac{2}{\delta_{0,i_1}} h_{i_1,0t} \times \dfrac{1}{h_{i_1,0t}^{\delta_{0,i_1}/2}} \dfrac{\partial h_{i_1,t}^{\delta_{0,i_1}/2}}{\partial\theta_i}(\theta_0),
\end{align}
where the derivatives involved in \eqref{derD} are defined for all $\theta\in\Theta$ recursively by
\begin{equation*}
\dfrac{\partial h_{i_1,t}^{\delta_{0,i_1}/2}(\theta)}{\partial \theta} = c_t(\theta) + \sum\limits_{i_2=1}^d\sum\limits_{i=1}^p B_i(i_1,i_2)\dfrac{\partial h_{i_2,t-i}^{\delta_{0,i_2}/2}}{\partial \theta},
\end{equation*}
with
\begin{equation}\label{AnDctilde}
\begin{aligned}
c_t(\theta) &= \left(0,\dots,1,0,\dots, \left(\varepsilon_{i_1,t-1}^+\right)^{\delta_{0,i_1}},0, \dots,\left(\varepsilon_{i_d,t-1}^+\right)^{\delta_{0,i_d}},0,\dots, \left(\varepsilon_{i_1,t-q}^+\right)^{\delta_{0,i_1}},0, \dots,\left(\varepsilon_{i_d,t-q}^+\right)^{\delta_{0,i_d}},\right.\\ &\qquad\left. 0,
,\dots, \left(\varepsilon_{i_1,t-1}^-\right)^{\delta_{0,i_1}},0, \dots,\left(\varepsilon_{i_d,t-1}^-\right)^{\delta_{0,i_d}},0,\dots, \left(\varepsilon_{i_1,t-q}^-\right)^{\delta_{0,i_1}},0, \dots,\left(\varepsilon_{i_d,t-q}^-\right)^{\delta_{0,i_d}},\right.\\ &\qquad\left. 0,
,\dots, h_ {i_1,t-1}^{\delta_{0,i_1}/2},0, \dots,h_ {i_d,t-1}^{\delta_{0,i_d}/2},0,\dots, h_ {i_1,t-p}^{\delta_{0,i_1}/2},0, \dots,h_ {i_d,t-p}^{\delta_{0,i_d}/2},\dots, 0\right)'.
\end{aligned}
\end{equation}
The distribution of $\eta_t$ is non-degenerated, so Equation \eqref{lambdaV} becomes 
\begin{equation*}
\lambda'V = \lambda'\mathbb{S}_{t-1:t-m} + \mu'\mathbf{d}'\mathbf{H}'\mathds{1} = 0,\quad \mbox{a.s.} 
\end{equation*}
where  $\mathds{1}$ represents a vector composed by $1$ of size $d^2\times 1$.
%Since
%\begin{align}\label{derD}
%\dfrac{\partial D_{0t}(i_1,i_1)}{\partial \theta_{i}} =\dfrac{\partial \left(h_{i_1,t}^{\delta_{0,i_1}/2}\right)^{1/\delta_{0,i_1}}}{\partial \theta_{i}}(\theta_0) = \dfrac{1}{\delta_{0,i_1}} h_{i_1,0t}^{1/2} \times \dfrac{1}{h_{i_1,0t}^{\delta_{0,i_1}/2}} \dfrac{\partial h_{i_1,t}^{\delta_{0,i_1}/2}}{\partial\theta_i}(\theta_0).
%\end{align}
%
For $i_1=1,\dots,d$ and in view of \eqref{derD} we can write
\begin{equation}\label{lambdaV2}
\lambda' V = \lambda' \mathbb{S}_{t-1:t-m} h_{i_1,t}^{\delta_{0,i_1}/2}(\theta_0) + \sum\limits_{i=1}^{s_1}\mu_{i}^\ast \dfrac{\partial h_{i_1,t}^{\delta_{0,i_1}/2}(\theta_0)}{\partial \theta_i} = 0, \quad \mbox{a.s.}
\end{equation}
where $\mu_i^\ast = 2\mu_i / \delta_{0,i_1}$.

Denote by $R_t$ a random variable measurable with respect to $\sigma\{\eta_u, u\leq t\}$ whose value will be modified along the proof. Thus we write
\begin{equation*}
h_{i_1,t}^{\delta_{0,i_1}/2} (\theta_0)= \sum\limits_{i_2=1}^d\left[A_{01}^+(i_1,i_2)(\varepsilon_{i_2,t-1}^+)^{\delta_{0,i_2}} + A_{01}^-(i_1,i_2)(\varepsilon_{i_2,t-1}^-)^{\delta_{0,i_2}}\right] + R_{t-2}.
\end{equation*}
We remind that $\underline{\varepsilon}_t^+ = H_{0t}^{1/2}\eta_t^+$ and $\underline{\varepsilon}_t^- = H_{0t}^{1/2}\eta_t^-$. We decompose Equation \eqref{lambdaV2} in two terms. The first one of \eqref{lambdaV2} can be rewritten 
\begin{align}\label{lambdaVh}
\lambda'\mathbb{S}_{t-1:t-m}h_{i_1,t}^{\delta_{0,i_1}/2}(\theta_0) &= \left\{\sum\limits_{i_2=1}^d \left[A_{01}^+(i_1,i_2)\left(\sum\limits_{j_1=1}^dH_{0,t-1}^{1/2}(i_2,j_1)\eta_{j_1,t-1}^+\right)^{\delta_{0,i_2}} \right.\right.\nonumber\\
&\qquad\qquad \left.\left. + A_{01}^-(i_1,i_2)\left(\sum\limits_{j_1=1}^dH_{0,t-1}^{1/2}(i_2,j_1)\eta_{j_1,t-1}^-\right)^{\delta_{0,i_2}}\right]\right\}R_{t-2}\nonumber\\
&\qquad + \left\{\sum\limits_{i_2=1}^d \left[A_{01}^+(i_1,i_2)\left(\sum\limits_{j_1=1}^dH_{0,t-1}^{1/2}(i_2,j_1)\eta_{j_1,t-1}^+\right)^{\delta_{0,i_2}} \right.\right.\nonumber\\
&\qquad\qquad \left.\left. + A_{01}^-(i_1,i_2)\left(\sum\limits_{j_1=1}^dH_{0,t-1}^{1/2}(i_2,j_1)\eta_{j_1,t-1}^-\right)^{\delta_{0,i_2}}\right]\right\}\nonumber\\
&\qquad \times \left(\lambda_{1}\sum\limits_{i=1}^d \eta_{i,t-1}^2\right) + \left(\lambda_{1}\sum\limits_{i=1}^d \eta_{i,t-1}^2\right)R_{t-2} + R_{t-2},
\end{align}
by using the fact that  
\begin{align*}
\lambda'\mathbb{S}_{t-1:t-m}  &=  \lambda_1S_{t-1} + R_{t-2}
= \lambda_{1}\sum\limits_{i=1}^d \eta_{i,t-1}^2 + R_{t-2}.
\end{align*}
The second term of \eqref{lambdaV2} can also be rewritten as
%\small{
\begin{align}\label{hi1t}
\mu^{\ast'}\dfrac{\partial h_{i_1,t}^{\delta_{0,i_1}/2}(\theta_0)}{\partial \theta} &= \sum\limits_{i_2=1}^d \left[\mu_{i_1+i_2d}^\ast(\varepsilon_{i_2,t-1}^+)^{\delta_{0,i_2}} + \mu_{i_1+(i_2+q)d^2}^\ast(\varepsilon_{i_2,t-1}^-)^{\delta_{0,i_2}}\right] + R_{t-2}\nonumber\\
&= \sum\limits_{i_2=1}^d \left[\mu_{i_1+i_2d}^\ast \left(\sum\limits_{j_1=1}^d H_{0,t-1}^{1/2}(i_2,j_1)\eta_{j_1,t-1}^+\right)^{\delta_{0,i_2}} \right.\nonumber\\
&\qquad\qquad \left.  + \mu_{i_1 + (i_2+q)d^2}^\ast\left(\sum\limits_{j_1=1}^d H_{0,t-1}^{1/2}(i_2,j_1)\eta_{j_1,t-1}^-\right)^{\delta_{0,i_2}}\right]+R_{t-2},
\end{align}
%}
%\normalsize
where the vector $\mu^\ast=(\mu^\ast_{1},\dots,\mu^\ast_{s_1})'$.

Combining the expressions \eqref{lambdaVh} and \eqref{hi1t}, Equation \eqref{lambdaV} comes down almost surely to
\begin{align*}
\lambda'V &= \left\{\sum\limits_{i_2=1}^d \left[A_{01}^+(i_1,i_2)\left(\sum\limits_{j_1=1}^dH_{0,t-1}^{1/2}(i_2,j_1)\eta_{j_1,t-1}^+\right)^{\delta_{0,i_2}}
 \right.\right.\nonumber\\
&\qquad\qquad \left.\left. 
+ A_{01}^-(i_1,i_2)\left(\sum\limits_{j_1=1}^dH_{0,t-1}^{1/2}(i_2,j_1)\eta_{j_1,t-1}^-\right)^{\delta_{0,i_2}}\right]\right\}R_{t-2}\nonumber\\
&\qquad\qquad + \left\{\sum\limits_{i_2=1}^d \left[A_{01}^+(i_1,i_2)\left(\sum\limits_{j_1=1}^dH_{0,t-1}^{1/2}(i_2,j_1)\eta_{j_1,t-1}^+\right)^{\delta_{0,i_2}} \right.\right.\nonumber\\
&\qquad \left.\left. + A_{01}^-(i_1,i_2)\left(\sum\limits_{j_1=1}^dH_{0,t-1}^{1/2}(i_2,j_1)\eta_{j_1,t-1}^-\right)^{\delta_{0,i_2}}\right]\right\}
%\nonumber\\
%&\qquad \times 
\left(\lambda_{1}\sum\limits_{i=1}^d \eta_{i,t-1}^2\right) + \left(\sum\limits_{i=1}^d \eta_{i,t-1}^2\right)R_{t-2}\\
&\qquad +
R_{t-2}\sum\limits_{i_2=1}^d \left[ \left(\sum\limits_{j_1=1}^d H_{0,t-1}^{1/2}(i_2,j_1)\eta_{j_1,t-1}^+\right)^{\delta_{0,i_2}} 
  + \left(\sum\limits_{j_1=1}^d H_{0,t-1}^{1/2}(i_2,j_1)\eta_{j_1,t-1}^-\right)^{\delta_{0,i_2}}\right]+R_{t-2} = 0,
\end{align*}
or equivalent to the two equations
\begin{align}\label{EqPos}
&\left[\sum\limits_{i_2=1}^d A_{01}^+(i_1,i_2)\left(\sum\limits_{j_1=1}^dH_{0,t-1}^{1/2}(i_2,j_1)\eta_{j_1,t-1}^+\right)^{\delta_{0,i_2}}\right]\left[\lambda_{1}\sum\limits_{i_2 = 1}^d\left(\eta_{i_2,t-1}^+\right)^2+ R_{t-2}\right]\nonumber\\
&\quad + R_{t-2}\sum\limits_{i_2=1}^d  \left(\sum\limits_{j_1=1}^d H_{0,t-1}^{1/2}(i_2,j_1)\eta_{j_1,t-1}^+\right)^{\delta_{0,i_2}} + R_{t-2}\sum\limits_{i_2 = 1}^d\left(\eta_{i_2,t-1}^+\right)^2 + R_{t-2} = 0 \quad \mbox{a.s.}
\end{align}
\begin{align}\label{EqNeg}
&\left[\sum\limits_{i_2=1}^d A_{01}^-(i_1,i_2)\left(\sum\limits_{j_1=1}^dH_{0,t-1}^{1/2}(i_2,j_1)\eta_{j_1,t-1}^-\right)^{\delta_{0,i_2}}\right]\left[\lambda_{1}\sum\limits_{i_2 = 1}^d\left(\eta_{i_2,t-1}^-\right)^2 + R_{t-2}\right]\nonumber\\
&\quad + R_{t-2}\sum\limits_{i_2=1}^d \left(\sum\limits_{j_1=1}^d H_{0,t-1}^{1/2}(i_2,j_1)\eta_{j_1,t-1}^-\right)^{\delta_{0,i_2}} + R_{t-2}\sum\limits_{i_2 = 1}^d\left(\eta_{i_2,t-1}^-\right)^2 + R_{t-2} = 0 \quad \mbox{a.s.}.
\end{align}
When $d=1$, from \eqref{EqPos} and \eqref{EqNeg} we retrieve an equation of the following form obtained by \cite{CF-PortTest}
\begin{equation*}
f(y) = a\vert y \vert^{\underline{\delta}_0+2} + b \vert y \vert^{\underline{\delta}_0} + cy^2 +d = 0, 
\end{equation*}
which cannot have more than $3$ positive roots or more than $3$ negative roots, except if $a = b = c = d = 0$.

When $d\geq2$ and also from \eqref{EqPos} and \eqref{EqNeg}, for a fixed component, we obtain an equation of the form
\begin{equation*}
f(y) = \sum\limits_{i=1}^d a_i\vert y \vert^{\delta_{0,i}+2} + \sum\limits_{i=1}^d b_i \vert y \vert^{\delta_{0,i} + 1} + \sum\limits_{i=1}^d c_i \vert y \vert^{\delta_{0,i}} + ay^2 + b\vert y \vert + c = 0.
\end{equation*}
Note that an equation of this form can not have more than  $3(d+1)$ non negative roots or more than $3(d+1)$ non positive roots for $d\geq 2$, except if  $a_i = b_i = c_i = a = b = c = 0$.

By Assumption \textbf{A9}, Equations \eqref{EqPos} and \eqref{EqNeg} imply that $\lambda_{1}\left[\sum_{i_2=1}^d A_{01}^+(i_1,i_2) + A_{01}^-(i_1,i_2)\right] = 0$. But under the assumption \textbf{A4}, if $p>0$, $\mathcal{A}_0(1)^+ + \mathcal{A}_0^- \neq 0$. It is impossible to have $A_{01}^+(i_1,i) = A_{01}^+(i_1,i) = 0$, for all $i=1,\ldots,d$. Then, there exists an $i_0$ such that $A_1(i_1,i_0)^+ + A_1(i_1,i_0)^- \neq 0$ and we then have  $\lambda_{1} = 0$.

In the general case, Equation \eqref{Raisonnement} necessarily leads
\begin{equation*}
A_{01}^+(i_1,i_0) + A_{01}^-(i_1,i_0) = \cdots = A_{0q}^+(i_1,i_0) + A_{0q}^-(i_1,i_0) = 0, \qquad \forall i_0,i_1 = 1,\ldots, d,
\end{equation*}
that is impossible under the assumption $\mathbf{A4}$ and $\lambda = 0$. This is in contradiction with $\lambda'V = 0$, almost surely, that leads that the assumption of non invertibility of matrix $D$ is absurd.
\zak
\subsection{Proof of Theorem \ref{thm-Port-test-connu}}
The almost sure convergence of $\hat{D}$ to $D$ as $n$ goes to infinity is easy to show using the consistency result. We remind the expression of the matrix $D$ 
\begin{equation*}
D = \Sigma_{{\mathbf{r}}_m} + C_mJ^{-1}IJ^{-1}C_m' + C_m\Sigma_{\hat{\theta}_n,{\mathbf{r}}_m} + \Sigma_{\hat{\theta}_n,{\mathbf{r}}_m}'C_m'.
\end{equation*}
%where $\Sigma_{{\mathbf{r}}_m}=d^2\left(\mathbb{E}\left[\eta_{it}^4\right]-1\right)^2I_m$ and $
%\Sigma_{\hat{\theta}_n,{\mathbf{r}}_m}=-\left(\mathbb{E}\left[\eta_{it}^4\right]-1\right)C_mJ^{-1}$ (see Equation \eqref{Loi_as_connue}). 
%%and is omitted.
%
The matrix $D$ can be rewritten as
\begin{equation*}
D = \Sigma_{{\mathbf{r}}_m} + A + B + B',
\end{equation*}
where the matrices $A$ and $B$ are given by
\begin{align*}
A &= (C_m - \hat{C}_m)J^{-1}IJ^{-1}C_m' + \hat{C}_m(J^{-1}-\hat{J}^{-1})IJ^{-1}C_m' + \hat{C}_m\hat{J}^{-1}(I - \hat{I})J^{-1}C_m' \\
&\qquad + \hat{C}_m\hat{J}^{-1}\hat{I}(J^{-1} - \hat{J}^{-1})C_m' + \hat{C}_m\hat{J}^{-1}\hat{I}\hat{J}^{-1}(C_m' - \hat{C}_m') + \hat{A},\\
B &= (C_m - \hat{C}_m)\Sigma_{\hat{\theta}_n,{\mathbf{r}}_m} + \hat{C}_m(\Sigma_{\hat{\theta}_n,{\mathbf{r}}_m} - \hat{\Sigma}_{\hat{\theta}_n,{\mathbf{r}}_m}) + \hat{B},
\end{align*}
with
$\hat{A} = \hat{C}_m\hat{J}^{-1}\hat{I}\hat{J}^{-1}\hat{C}_m'$ and $\hat{B} = \hat{C}_m\hat{\Sigma}_{\hat{\theta}_n,{\mathbf{r}}_m}$ where $
\hat\Sigma_{\hat{\theta}_n,{\mathbf{r}}_m}=-\left(\hat{\kappa}_i-1\right)\hat{C}_m\hat{J}^{-1}$. 
Finally we have
\begin{align*}
D-\hat{D} &= (\Sigma_{{\mathbf{r}}_m} - \hat{\Sigma}_{{\mathbf{r}}_m}) + (A - \hat{A}) + (B-\hat{B}) + (B' - \hat{B}').
\end{align*}
%where $\hat{\Sigma}_{{\mathbf{r}}_m}=d^2(\hat\kappa_{i}-1)^2I_m$.
%
For any multiplicative norm we have
\begin{align*}
\Vert D - \hat{D} \Vert &\leq \Vert \Sigma_{{\mathbf{r}}_m} - \hat{\Sigma}_{{\mathbf{r}}_m} \Vert + \Vert A - \hat{A} \Vert + \Vert B - \hat{B} \Vert + \Vert B' - \hat{B}' \Vert.
\end{align*}
Observe that 
\begin{align}\nonumber
\Vert A - \hat{A} \Vert &\leq \Vert C_m - \hat{C}_m \Vert \Vert J^{-1} \Vert \Vert I \Vert \Vert J^{-1} \Vert  \Vert C_m' \Vert + \Vert \hat{C}_m \Vert \Vert J^{-1} - \hat{J}^{-1} \Vert \Vert I \Vert \Vert J^{-1} \Vert \Vert C_m' \Vert\\ 
&\qquad  \nonumber+ \Vert \hat{C}_m \Vert \Vert \hat{J}^{-1} \Vert \Vert I - \hat{I} \Vert \Vert J^{-1} \Vert \Vert C_m'\Vert +  \Vert \hat{C}_m \Vert \Vert \hat{J}^{-1} \Vert \Vert \hat{I} \Vert \Vert J^{-1} - \hat{J}^{-1} \Vert \Vert C_m' \Vert\\
&\qquad \nonumber+ \Vert \hat{C}_m \Vert \Vert \hat{J}^{-1} \Vert \Vert \hat{I} \Vert \Vert \hat{J}^{-1} \Vert \Vert C_m - \hat{C}_m \Vert, \\
&\leq \nonumber\Vert C_m - \hat{C}_m \Vert \Vert J^{-1} \Vert \Vert I \Vert \Vert J^{-1} \Vert  \Vert C_m' \Vert + \Vert \hat{C}_m \Vert \Vert J^{-1} \Vert \Vert \hat{J} - J \Vert \Vert \hat{J} \Vert \Vert I \Vert \Vert J^{-1} \Vert \Vert C_m' \Vert\\ 
&\qquad \nonumber + \Vert \hat{C}_m \Vert \Vert \hat{J}^{-1} \Vert \Vert I - \hat{I} \Vert \Vert J^{-1} \Vert \Vert C_m'\Vert +  \Vert \hat{C}_m \Vert \Vert \hat{J}^{-1} \Vert \Vert \hat{I} \Vert \Vert J^{-1} \Vert \Vert \hat{J} - J \Vert \Vert \hat{J} \Vert \Vert C_m' \Vert\\
&\qquad + \Vert \hat{C}_m \Vert \Vert \hat{J}^{-1} \Vert \Vert \hat{I} \Vert \Vert \hat{J}^{-1} \Vert \Vert C_m - \hat{C}_m \Vert,\label{normeA}\\
%&\underset{n \to + \infty}{\longrightarrow} 0, \mbox{ p.s.}\\
\Vert B - \hat{B} \Vert &\leq \Vert C_m - \hat{C}_m \Vert \Vert \Sigma_{\hat{\theta}_n,{\mathbf{r}}_m} \Vert  + \Vert \hat{C}_m \Vert \Vert \Sigma_{\hat{\theta}_n,{\mathbf{r}}_m} - \hat{\Sigma}_{\hat{\theta}_n,{\mathbf{r}}_m} \Vert\label{normeB}. 
\end{align}
%To complete the proof of the convergence of the matrix $D$, it remains to show $\Vert \Sigma_{\hat{\theta}_n,{\mathbf{r}}_m} - \hat{\Sigma}_{\hat{\theta}_n,{\mathbf{r}}_m} \Vert $ and $\Vert\Sigma_{{\mathbf{r}}_m} - \hat{\Sigma}_{{\mathbf{r}}_m}\Vert$ are controlled.
%
%We conclude that $\hat{D} \to D$ almost surely when $n \to +\infty$.
%
In view of \eqref{MajH1.1} and \eqref{unif1}, we have $\Vert C_m \Vert<\infty$. We also have $\Vert I \Vert<\infty$. Because the matrix $J$ is nonsingular, we have $\Vert J^{-1} \Vert<\infty$ and $$\Vert \hat{J}^{-1}-J^{-1} \Vert\underset{n\to\infty}{\longrightarrow} 0,\quad a.s.$$
by consistency of $\hat\theta_n$. Under Assumption \textbf{A7}, we have $\vert \mathbb{E}\left[\eta_{t}'\eta_t-d\right] ^2\vert \leq K$. Using the previous arguments and also the strong consistency of  $\hat\theta_n$, we have
$$\vert \mathbb{E}\left[\eta_{t}'\eta_t-d\right] ^2  - \hat{\kappa}\vert\underset{n\to\infty}{\longrightarrow} 0,\quad a.s.\text{ and } \Vert C_m - \hat{C}_m \Vert\underset{n\to\infty}{\longrightarrow} 0,\quad a.s.$$ We then deduce that Equations \eqref{normeA} and \eqref{normeB} converge almost surely to 0 when $n\to\infty$ and the conclusion follows. Thus $\hat{D} \underset{n\to\infty}{\longrightarrow} D$ almost surely.%, when $n \to \infty$.

To conclude the proof of Theorem \ref{thm-Port-test-connu}, it suffices to use Theorem \ref{thmD-connu} and the following result:
if $\sqrt{n} \mathbf{\hat{r}}_m \xrightarrow[n\to\infty]{\mathrm{d}} \mathcal{N}\left(0,D\right)$, with $D$ nonsingular, and if $\hat{D}\underset{n\to\infty}{\longrightarrow} D$ in probability, then $
n \mathbf{\hat{r}}_m'\hat{D}^{-1}\mathbf{\hat{r}}_m \xrightarrow[n\to\infty]{\mathrm{d}} \chi^2_m.$
\zak
%%%%%%%%%%%%%%%%%%%%%%%%
\subsection{Proof of Remark~\ref{bobo}}
We suppose that $\boldsymbol{\mathcal{H}}_1$ holds true. 
One may rewrite the above arguments in order to prove that there exists a nonsingular matrix $D^\ast$ such that 
\begin{equation}\label{tt}
\sqrt{n} ( \mathbf{\hat{r}}_m - \mathbf{{r}}_m^0 )  \xrightarrow[n\to\infty]{\mathrm{d}} \mathcal{N}\left(0,D^\ast \right) \ .
\end{equation}
The matrix $D^\ast$ is given by ${D^\ast}= \Sigma_{{\mathbf{r}}_m^0} + {C_m^\ast}J^{-1}IJ^{-1}{C_m^\ast}' +{C_m^\ast}\Sigma_{\hat{\theta}_n,{\mathbf{r}}_m^0} + \Sigma_{\hat{\theta}_n,{\mathbf{r}}_m^0}'{C_m^\ast}'$, where the matrices $\Sigma_{{\mathbf{r}}_m^0}$ and  $\Sigma_{\hat{\theta}_n,{\mathbf{r}}_m^0}$ are obtained from the asymptotic distribution of 
$\sqrt{n}(\hat{\theta}'_n - \theta'_0, {\mathbf{r}}'_m-{\mathbf{r}^0}'_m)'$.
The $(h,i)$-th element of the matrix $C_m^\ast$ is geven by
\begin{align*}
C^\ast(h,i) := \mathbb{E}\left[S_{t-h} \dfrac{\partial S_t}{\partial \theta_i}+S_t \dfrac{\partial S_{t-h}}{\partial \theta_i}\right] .
\end{align*}
Now we write 
\begin{align*}\sqrt{n} {\hat D}^{-1/2}\mathbf{\hat{r}}_m    & =  {\hat D}^{-1/2} \sqrt{n}  (  \mathbf{\hat{r}}_m - \mathbf{r}^0_m ) + {\hat D}^{-1/2} \sqrt{n}  \mathbf{r}^0_m  \\ 
& =   D^{-1/2} \sqrt{n}  (  \mathbf{\hat{r}}_m - \mathbf{r}^0_m ) + { D}^{-1/2} \sqrt{n}  \mathbf{r}^0_m  + \mathrm{o}_{\mathbb P} (1) \ . 
\end{align*}
Then it holds that 
\begin{align}\label{opp}
 n\mathbf{\hat{r}}_m'\hat{D}^{-1}\mathbf{\hat{r}}_m & = \big  ( \sqrt{n} {\hat D}^{-1/2}\mathbf{\hat{r}}_m \big) '  \times  \big  ( \sqrt{n} {\hat D}^{-1/2}\mathbf{\hat{r}}_m \big) \nonumber \\ 
 & = n (  \mathbf{\hat{r}}_m - \mathbf{r}^0_m )' D^{-1} (  \mathbf{\hat{r}}_m - \mathbf{r}^0_m )  + 2 n (  \mathbf{\hat{r}}_m - \mathbf{r}^0_m )' D^{-1}  \mathbf{r}^0_m  + n { \mathbf{r}^0_m}'D^{-1}  \mathbf{r}^0_m  + \mathrm{o}_{\mathbb P} (1) .
\end{align}
By the ergodic theorem, $ (  \mathbf{\hat{r}}_m - \mathbf{r}^0_m )' D^{-1}  \mathbf{r}^0_m = \mathrm{o}_{\mathbb P} (1)$. By Lemma 17.1 in \cite{vdv}, the convergence \eqref{tt} implies that 
$$ (  \mathbf{\hat{r}}_m - \mathbf{r}^0_m )' D^{-1} (  \mathbf{\hat{r}}_m - \mathbf{r}^0_m )   \xrightarrow[n\to\infty]{\mathrm{d}} \sum_{i=1}^m \lambda_i Z_i^2$$ 
where $(Z_i)_{1\le i\le m}$ are i.i.d. with $\mathcal N(0,1)$ laws and the $\lambda_i$'s are the eigenvalues of the matrix $D^{-1/2}D^\ast D^{-1/2}$. 
Reporting these convergences in \eqref{opp}, we deduce that 
\begin{align*}
\mathbf{\hat{r}}_m'\hat{D}^{-1}\mathbf{\hat{r}}_m
 & =  (\mathbf{\hat{r}}_m - \mathbf{r}^0_m )' D^{-1} (  \mathbf{\hat{r}}_m - \mathbf{r}^0_m )  + 2  (  \mathbf{\hat{r}}_m - \mathbf{r}^0_m )' D^{-1}  \mathbf{r}^0_m  +  { \mathbf{r}^0_m}'D^{-1}  \mathbf{r}^0_m  + \mathrm{o}_{\mathbb P} (1)  \\ 
 & =  { \mathbf{r}^0_m}'D^{-1}  \mathbf{r}^0_m  + \mathrm{o}_{\mathbb P} (1) 
\end{align*}
and the remark is proved. \zak
%%%%%%%
\subsection{Proof of Corollary \ref{thm-Port-test-connubis}}
Note that if the model is correct we have %almost surely
$${\hat{r}}_0 = \dfrac{1}{n} \sum\limits_{t = h +1}^n [\underline{\varepsilon}_t'\tilde{H}_t^{-1}({\hat{\theta}_n})\underline{\varepsilon}_t-d]^2 \xrightarrow[n\to\infty]{\mathrm{a.s.}} \mathbb{E}[\underline{\varepsilon}_t'{H}_t^{-1}\underline{\varepsilon}_t-d]^2=
\mathbb{E}\left[\eta_{t}'\eta_t-d\right]^2.$$
From \eqref{DLTrm} we have $\sqrt{n}(\hat{r}_0-r_0)=\mathrm{o}_{\mathbb{P}}(1)$. Applying the central limit theorem to the process $([\underline{\varepsilon}_t'{H}_t^{-1}\underline{\varepsilon}_t-d]^2)_{t\in\mathbb{Z}}$, we obtain
\begin{align*}
\sqrt{n}\left(\hat{r}_0-r_0)\right)&=\frac{1}{\sqrt{n}}\sum_{t=1}^n\left( [\underline{\varepsilon}_t'{H}_t^{-1}\underline{\varepsilon}_t-d]^2-\mathbb{E}[\underline{\varepsilon}_t'{H}_t^{-1}\underline{\varepsilon}_t-d]^2 \right)+\mathrm{o}_{\mathbb{P}}(1)
\xrightarrow[n\rightarrow \infty]{\text{in law}}\mathcal{N}\left(0,\Phi\right).
\end{align*}
So we have $\sqrt{n}(\hat{r}_0-r_0)=\mathrm{O}_{\mathbb{P}}(1)$ and
 $\sqrt{n}(r_0-\mathbb{E}[\underline{\varepsilon}_t'{H}_t^{-1}\underline{\varepsilon}_t-d]^2 )=\mathrm{O}_{\mathbb{P}}(1)$. Now,
 using \eqref{DLTrm} and the ergodic theorem, we have
$$n\left(\frac{\hat{r}_h}{\hat{r}_0}- \frac{\hat{r}_h}{\mathbb{E}[\underline{\varepsilon}_t'{H}_t^{-1}\underline{\varepsilon}_t-d]^2 }\right)=\sqrt{n}\hat{r}_h\frac{\sqrt{n}\left(\mathbb{E}[\underline{\varepsilon}_t'{H}_t^{-1}\underline{\varepsilon}_t-d]^2 -\hat{r}_0\right)}{\mathbb{E}[\underline{\varepsilon}_t'{H}_t^{-1}\underline{\varepsilon}_t-d]^2 \hat{r}_0}=\mathrm{O}_{\mathbb{P}}(1),$$
which means $\sqrt{n}\hat{\rho}(h)={\sqrt{n}\hat{r}_h}/{\mathbb{E}[\underline{\varepsilon}_t'{H}_t^{-1}\underline{\varepsilon}_t-d]^2 }+\mathrm{O}_{\mathbb{P}}(n^{-1/2}).$
For $h=1,\dots,m$, it follows that
\begin{align}\label{hat_rho}
\sqrt{n}\hat{\rho}_m=\frac{\sqrt{n}\mathbf{\hat{r}}_m}{\mathbb{E}[\underline{\varepsilon}_t'{H}_t^{-1}\underline{\varepsilon}_t-d]^2 }+\mathrm{o}_{\mathbb{P}}(1).
\end{align}
Thus from \eqref{hat_rho} the asymptotic distribution of the sum of squared residuals autocorrelations $\sqrt{n}\hat{\rho}_m$ depends on the  distribution of $\sqrt{n}\mathbf{\hat{r}}_m$.
Consequently we have
$$\lim_{n\rightarrow\infty}\mathrm{Var}\left(\sqrt{n}\hat{\rho}_m\right)=
\lim_{n\rightarrow\infty}\mathrm{Var}\left(\frac{\sqrt{n}\mathbf{\hat{r}}_m}{\mathbb{E}[\underline{\varepsilon}_t'{H}_t^{-1}\underline{\varepsilon}_t-d]^2 }\right)=: D_{\hat{\rho}} = \frac{D}{\left(\mathbb{E}[\underline{\varepsilon}_t'{H}_t^{-1}\underline{\varepsilon}_t-d]^2 \right)^2}.$$
Thus the first result \eqref{rho1-connu} of Corollary \ref{thm-Port-test-connubis} is proved.

The proof the second result \eqref{rho2-connu} of Corollary \ref{thm-Port-test-connubis} is the same that the one given for Theorem \ref{thm-Port-test-connu} and the proof is completed.
\zak
\subsection{Proof of Theorem \ref{thmD-inconnu}}
We follow the arguments and the different steps that we used in the proof of Theorem \ref{thmD-connu}. As in the case where $\underline{\delta}_0$ was known, the proof is decomposed in the following points which will be treated in separate subsections.
%
%The proof of theorem presents similarities with the case where the power is known. For these reasons we only present the differences. As the the proof of Theorem \ref{thmD-connu} we decomposed in 3 steps.
\begin{enumerate}[$\qquad (i)$]
	\item Asymptotic impact of unknown initials values on the statistic ${\mathbf{\hat{r}}}_m$.
	\item Asymptotic distribution of $\sqrt{n}{\mathbf{\hat{r}}}_m$.
	\item Invertibility of the matrix ${\cal D}$.
\end{enumerate}  
There are many similarities with the proof of Theorem \ref{thmD-connu}. We only indicates where the fact that the power is estimated has an importance is our reasoning.

\textbf{$(i)$ Asymptotic impact of unknown initials values on the statistic  ${\mathbf{\hat{r}}}_m$}

The proof of the asymptotic impact of the initial values on the statistic ${\mathbf{\hat{r}}}_m$ is the same than the  one where $\underline{\delta}_0$ was known. It suffices to adapt this step by replacing $\theta$ by $\vartheta$ and $H_t$ by $\mathcal{H}_t$.

\textbf{$(ii)$ Asymptotic distribution of  $\sqrt{n}{\mathbf{\hat{r}}}_m$}

The asymptotic distribution of $\sqrt{n}\underline{\mathbf{\hat{r}}}_m$ is similar to that the one when the power $\underline{\delta}_0$ is assumed to be known. We adapt this step by replacing again $\theta$ by $\vartheta$ and $H_t$ by $\mathcal{H}_t$. The only difference resides in the estimations of the derivatives  when we differentiate with respect to $\delta_i$, $i=1,\dots,d$.

For instance the $(h,i)$-th element of the matrix ${\cal C}_m$ denoted by ${\cal C}(h,i)$ is given by
\begin{align*}
{\cal C}(h,i) = \mathbb{E}\left[S_{t-h} \dfrac{\partial S_t}{\partial \vartheta_i}\right] = -\mathbb{E}\left[S_{t-h}
\text{Tr}\left({\cal H}_{0t}^{-1}\dfrac{\partial {\cal H}_t(\vartheta_0)}{\partial\vartheta_i}\right)\right]
%\\&
= - \mathbb{E}\left[S_{t-h}\mathbf{h}'_t(i)\text{vec}(I_d)\right].
\end{align*}
Consequently we have
\begin{equation}\label{CmInc}
%\dfrac{\partial {\mathbf{r}}_m(\theta_0)}{\partial \theta'} \underset{n\to\infty}{\longrightarrow} 
{\cal C}_m := [{\cal C}(h,i)]_{1\leq h\leq m,1\leq i\leq s_0}
= - \mathbb{E}\left[\left(\mathbb{S}_{t-1:t-m}\right)\left(\mathbf{h}'_t\text{vec}(I_d)\right)'\right].
\end{equation}
%
%where $\mathbb{S}_{t-1:t-m} = (S_{t-1},\ldots, S_{t-m})'$. 

\textbf{$(iii)$ Invertibility of the matrix ${\cal D}$}

The proof of the invertibility of matrix ${\cal D}$ needs to have some modifications compared to the case where the power $\underline{\delta}_0$ is assumed to be known. The start of the proof stay identical, it suffices only  to replace $H_t$ by $\mathcal{H}_t$ and $\theta$ by $\vartheta$. 
We rewrite $\mathbf{h}_t(i)$ as follow
\begin{align*}
\mathbf{h}_t(i) = \left[\text{vec}\left({\cal H}_{0t}^{-1/2}\dfrac{\partial {\cal H}_t(\vartheta_0)}{\partial \vartheta_i}{\cal H}_{0t}^{-1/2}\right)\right] = \left[{\cal H}_{0t}^{-1/2}\otimes {\cal H}_{0t}^{-1/2}\right] \text{vec}\left(\dfrac{\partial {\cal H}_t(\vartheta_0)}{\partial \vartheta_i}\right) = \mathbf{H}_t\mathbf{d}_t(i),
\end{align*}
where $\mathbf{H}_t =\left[{\cal H}_{0t}^{-1/2}\otimes {\cal H}_{0t}^{-1/2}\right] $ and $\mathbf{d}_t(i) = \text{vec}\left({\partial {\cal H}_t(\vartheta_0)}/{\partial \vartheta_i}\right)$. 
Thus we define the matrix of size $d^2\times s_0$,  $\mathbf{d}_t = (\mathbf{d}_t(1) \vert \ldots \vert \mathbf{d}_t({s_0}))$ such that $\mathbf{h}_t = \mathbf{H}_t\mathbf{d}_t$.
To study the invertibility of the matrix ${\cal D}$ we let $V = \mathbb{S}_{t-1:t-m}  S_t - {\cal C}_m\mathcal{J}^{-1}{\partial l_t(\vartheta_0)}/{\partial \vartheta}$ such that $\mathbb{E}[VV']={\cal D}$.
We can also rewrite the vector $V$ as
\begin{equation*}
V = \mathbb{S}_{t-1:t-m}  S_t + {\cal C}_m\mathcal{J}^{-1}\mathbf{d}_t'\mathbf{H}_t'\text{vec}(s_t).
\end{equation*}
If the matrix $\mathbb{E}[VV']$ is singular, then there exists a vector $\lambda = (\lambda_1,\ldots, \lambda_{m})'$ not equal to zero such that 
\begin{equation}\label{lambdaVinc}
\lambda'V = \lambda'\mathbb{S}_{t-1:t-m}  S_t + \mu  \mathbf{d}_t'\mathbf{H}_t'\text{vec}(s_t) = 0, \quad \mbox{a.s.},
\end{equation}
with $\mu = \lambda'{\cal C}_m\mathcal{J}^{-1}$. We have $\mu \neq 0$, else $\lambda'\mathbb{S}_{t-1:t-m} S_t = 0$ almost surely, that implies there exists $j\in \{1,\ldots, m\}$ such that $S_{t-j}$ be mesurable respect to the $\sigma$-field $\{S_r, t-1\leq r\leq t-m\}$ with $r\neq t-j$. That is impossible because the $S_t$ are independent and not degenerated. 
Consequently \eqref{lambdaVinc} becomes 
\begin{align}
\mu'\mathbf{d}_t' &= \sum\limits_{i=1}^{s_0} \mu_i\mathbf{d}_t(i) = \sum\limits_{i=1}^{s_0}\mu_i\text{vec}\left(\dfrac{\partial {\cal H}_t(\vartheta_0)}{\partial \vartheta_i}\right) = \sum\limits_{i=1}^{s_0}\mu_i\dfrac{\partial}{\partial \vartheta_i}\left[(D_{0t}\otimes D_{0t})\text{vec}(R_0)\right],\quad \mbox{a.s.}
\label{mud1}
\end{align}
We can then rewrite \eqref{mud1} in order to separate the derivatives of the matrix $\mathcal{H}_t$ when we differentiate with respect to the vectors $\theta$ and $\underline{\delta}$. It follows that 
\begin{align}
\mu'\mathbf{d}_t'&=\sum\limits_{i=1}^{s_2}\mu_i\dfrac{\partial \left[(D_{0t}\otimes D_{0t})\text{vec}(R_0)\right]}{\partial \theta_i} + \sum\limits_{i=s_2+1}^{s_0}\mu_i\dfrac{\partial (D_t\otimes D_t)}{\partial \delta_i}\text{vec}(R_0),\quad \mbox{a.s.}\nonumber
\\
&=\sum\limits_{i=1}^{s_1}\mu_i\dfrac{\partial (D_{0t}\otimes D_{0t})}{\partial \theta_i}\text{vec}(R_0) + \sum\limits_{i = s_1 + 1}^{s_2}\mu_i (D_{0t} \otimes D_{0t})\dfrac{\partial \text{vec}(R_0)}{\partial \theta_i} + \sum\limits_{i=s_2+1}^{s_0}\mu_i\dfrac{\partial (D_t\otimes D_t)}{\partial \delta_i}\text{vec}(R_0)=0\quad \mbox{a.s.}
\label{mud3}
\end{align}
Since the vectors $\partial \text{vec}(R_0)/\partial \theta_i$, $i = s_1+1,\ldots, s_2$ are linearly independent, the vector  $(\mu_{s_1+1}, \ldots, \mu_{s_2})'$ is null and thus Equation \eqref{mud3} yields
\begin{equation}\label{mud}
\sum\limits_{i=1}^{s_1} \mu_i \dfrac{\partial (D_{0t}\otimes D_{0t})}{\partial \theta_i}\text{vec}(R_0)+ \sum\limits_{i=s_2+1}^{s_0}\mu_i\dfrac{\partial (D_t\otimes D_t)}{\partial \delta_i}\text{vec}(R_0)=0,\quad \mbox{a.s.}
\end{equation}
The rows $1,d+1,\ldots,d^2$ of the Equation  \eqref{mud}
%$
%\sum\limits_{i=s_2+1}^{s_0} \mu_i {\partial (D_{0t}\otimes D_{0t})}/{\partial \delta_i}\text{vec}(R_0)=0,\quad \mbox{a.s.}$ 
yield
\begin{align}\label{derht-delta}
\sum\limits_{i=1}^{s_1}\mu_i\dfrac{\partial\underline{h}_t(\vartheta_0)}{\partial \theta_i}+\sum\limits_{i=s_2+1}^{s_0}\mu_i\dfrac{\partial\underline{h}_t(\vartheta_0)}{\partial \delta_i} = 0,\quad \text{a.s.}
\end{align}
We have for $i_1 = 1,\ldots, d$ and $i=1,\ldots,s_1$
\begin{align}\label{derDinc}
\dfrac{\partial h_{i_1,t}(\vartheta_0)}{\partial  \theta_{i}} =\dfrac{\partial \left(h_{i_1,t}^{\delta_{0,i_1}/2}\right)^{2/\delta_{0,i_1}}}{\partial \theta_{i}}(\vartheta_0) = \dfrac{2}{\delta_{0,i_1}} h_{i_1,0t} \times \dfrac{1}{h_{i_1,0t}^{\delta_{0,i_1}/2}} \dfrac{\partial h_{i_1,t}^{\delta_{0,i_1}/2}}{\partial\theta_i}(\vartheta_0),
\end{align}
where the derivatives involved in \eqref{derDinc} are defined for all $\vartheta\in\Delta$ recursively by
\begin{equation*}
\dfrac{\partial h_{i_1,t}^{\delta_{i_1}/2}(\vartheta)}{\partial \theta} = c_t(\vartheta) + \sum\limits_{i_2=1}^d\sum\limits_{i=1}^p B_i(i_1,i_2)\dfrac{\partial h_{i_2,t-i}^{\delta_{i_2}/2}(\vartheta)}{\partial \theta},
\end{equation*}
with
\begin{equation}\label{truc1}
\begin{aligned}
c_t(\vartheta) &= \left(0,\dots,1,0,\dots, \left(\varepsilon_{i_1,t-1}^+\right)^{\delta_{i_1}},0, \dots,\left(\varepsilon_{i_d,t-1}^+\right)^{\delta_{i_d}},0,\dots, \left(\varepsilon_{i_1,t-q}^+\right)^{\delta_{i_1}},0, \dots,\left(\varepsilon_{i_d,t-q}^+\right)^{\delta_{i_d}},\right.\\ &\qquad\left. 0,
,\dots, \left(\varepsilon_{i_1,t-1}^-\right)^{\delta_{i_1}},0, \dots,\left(\varepsilon_{i_d,t-1}^-\right)^{\delta_{i_d}},0,\dots, \left(\varepsilon_{i_1,t-q}^-\right)^{\delta_{i_1}},0, \dots,\left(\varepsilon_{i_d,t-q}^-\right)^{\delta_{i_d}},\right.\\ &\qquad\left. 0,
,\dots, h_ {i_1,t-1}^{\delta_{i_1}/2},0, \dots,h_ {i_d,t-1}^{\delta_{i_d}/2},0,\dots, h_ {i_1,t-p}^{\delta_{i_1}/2},0, \dots,h_ {i_d,t-p}^{\delta_{i_d}/2},\dots, 0\right)'.
\end{aligned}
\end{equation}
So we can focus  on the derivatives with respect to $\underline{\delta}$:
%\begin{align}\label{Dtau}
%\dfrac{\partial D_t(i_1, i_1)}{\partial\delta_{j}}  & = h_{i_1,t}^{1/2}\left[ -\boldsymbol{\delta}_{j,i_1} \dfrac{1}{\delta_{i_1}^2} \log\left(h_{i_1,t}^{\delta_{i_1}/2}\right) +\dfrac{1}{\delta_{i_1}} \dfrac{1}{h_{i_1,t}^{\delta_{i_1}/2}} \dfrac{\partial h_ {i_1,t}^{\delta_{i_1}/2}}{\partial \delta_{j}}\right],
%\end{align}
%%
%where $\boldsymbol{\delta}_{j,i_1}$ denotes the Kronecker symbol.
%
%
\begin{align}\label{derDinc1}
\dfrac{\partial h_{i_1,t}}{\partial \delta_j} &= \dfrac{2}{\delta_{i_1}}h_{i_1,t}\left[-\boldsymbol{\delta}_{j,i_1}\dfrac{1}{\delta_{i_1}}\log\left(h_{i_1,t}^{\delta_{i_1}/2}\right) +  \dfrac{1}{h_{i_1,t}^{\delta_{i_1}/2}}\dfrac{\partial h_{i_1,t}^{\delta_{i_1}/2}}{\partial \delta_j}\right], \qquad j=1,\ldots, d
\end{align}
with
\begin{align}\nonumber
\dfrac{\partial h_{i_1,t}^{\delta_{i_1}/2}}{\partial \delta_j} &= \sum\limits_{i=1}^q\left[A_i^+(i_1,j)\log(\varepsilon_{j,t-i}^+)(\varepsilon_{j,t-i}^+)^{\delta_j} + A_i^-(i_1,j)\log(\varepsilon_{j,t-i}^-)(\varepsilon_{j,t-i}^-)^{\delta_j}\right]
%\\ &\qquad 
+ \sum\limits_{i_2=1}^d \sum\limits_{i=1}^p B_i(i_1,i_2)\dfrac{\partial h_{i_2,t-i}^{\delta_{i_2}/2}}{\partial \delta_j}\\
&=A_1^+(i_1,j)\log(\varepsilon_{j,t-1}^+)(\varepsilon_{j,t-1}^+)^{\delta_j} + A_1^-(i_1,j)\log(\varepsilon_{j,t-1}^-)(\varepsilon_{j,t-1}^-)^{\delta_j} + R_{t-2},\label{truc2}
\end{align}
where $\boldsymbol{\delta}_{j,i_1}$ denotes the Kronecker symbol. We also remind that
\begin{equation}\label{truc3}
h_{i_1,t}^{\delta_{0,i_1}/2} (\vartheta_0)= \sum\limits_{i_2=1}^d\left[A_{01}^+(i_1,i_2)(\varepsilon_{i_2,t-1}^+)^{\delta_{0,i_2}} + A_{01}^-(i_1,i_2)(\varepsilon_{i_2,t-1}^-)^{\delta_{0,i_2}}\right] + R_{t-2}.
\end{equation}
The distribution of $\eta_t$ is non-degenerated, so Equation \eqref{lambdaVinc} becomes 
\begin{equation*}
\lambda'V = \lambda'\mathbb{S}_{t-1:t-m} + \mu'\mathbf{d}'\mathbf{H}'\mathds{1} = 0,\quad \mbox{a.s.} 
\end{equation*}
In view of \eqref{derDinc} and \eqref{derDinc1}, we can finally write
%\begin{equation}\label{lambdaVIncbis}
\begin{align}\nonumber
\lambda'V = \lambda'\mathbb{S}_{t-1:t-m}h_{i_1,t}^{\delta_{i_1}/2}(\vartheta_0)& + \sum\limits_{i=1}^{s_1}\mu_i^\ast \dfrac{\partial h_{i_1,t}^{\delta_{i_1}/2}(\vartheta_0)}{\partial \theta_i} + \sum\limits_{i=1}^d \mu_{i+s_2}^\ast \dfrac{\partial h_{i_1,t}^{\delta_{i_1}/2}(\vartheta_0)}{\partial \delta_i}\\&\hspace*{2cm} - \mu_{i_1+ s_2}^\ast h_{i_1,t}^{\delta_{i_1}/2}(\vartheta_0)\log\left(h_{i_1,t}^{\delta_{i_1}/2}(\vartheta_0)\right) = 0, \quad \mbox{a.s.}\label{lambdaVIncbis}
\end{align}
%\end{equation}
where  $\mu_i^\ast = 2\mu_i / \delta_{0,i_1}$, $\mu_{i+s_2}^\ast=2\mu_{i+s_2}/\delta_{0,i_1}$ and when $i=i_1$ we have $\mu_{i_1 + s_2}^\ast = 2\mu_{i_1 + s_2} / \delta_{0,i_1}^2$.
Recall that $\underline{\varepsilon}_t^+ = \mathcal{H}_t^{1/2}\eta_t^+$ and $\underline{\varepsilon}_t^- = \mathcal{H}_t^{1/2}\eta_t^-$ and we  decomposed  \eqref{lambdaVIncbis} in four terms.
The first one leads to
\begin{align}\label{part1}
\lambda'\mathbb{S}_{t-1:t-m}h_{i_1,t}^{\delta_{i_1}/2}(\vartheta_0) &= \left\{\sum\limits_{i_2=1}^d \left[A_{01}^+(i_1,i_2)\left(\sum\limits_{j_1=1}^d\mathcal{H}_{0,t-1}^{1/2}(i_2,j_1)\eta_{j_1,t-1}^+\right)^{\delta_{0,i_2}} \right.\right.\nonumber\\
&\qquad\qquad \left.\left. + A_{01}^-(i_1,i_2)\left(\sum\limits_{j_1=1}^d\mathcal{H}_{0,t-1}^{1/2}(i_2,j_1)\eta_{j_1,t-1}^-\right)^{\delta_{0,i_2}}\right]\right\}R_{t-2}\nonumber\\
&\qquad + \left\{\sum\limits_{i_2=1}^d \left[A_{01}^+(i_1,i_2)\left(\sum\limits_{j_1=1}^d\mathcal{H}_{0,t-1}^{1/2}(i_2,j_1)\eta_{j_1,t-1}^+\right)^{\delta_{0,i_2}} \right.\right.\nonumber\\
&\qquad\qquad \left.\left. + A_{01}^-(i_1,i_2)\left(\sum\limits_{j_1=1}^d\mathcal{H}_{0,t-1}^{1/2}(i_2,j_1)\eta_{j_1,t-1}^-\right)^{\delta_{0,i_2}}\right]\right\}\nonumber\\
&\qquad \times \left(\lambda_{1}\sum\limits_{i=1}^d \eta_{i,t-1}^2\right) + \left(\lambda_{1}\sum\limits_{i=1}^d \eta_{i,t-1}^2\right)R_{t-2} + R_{t-2},
\end{align}
by using \eqref{truc3} and the fact that  
\begin{align*}
\lambda'\mathbb{S}_{t-1:t-m}  &=  \lambda_1S_{t-1} + R_{t-2}
= \lambda_{1}\sum\limits_{i=1}^d \eta_{i,t-1}^2 + R_{t-2}.
\end{align*}
Using \eqref{truc1}, the second term of  \eqref{lambdaVIncbis} can be rewritten
%The second term of \eqref{lambdaV2} can also be rewritten as
%\small{
\begin{align}\label{part2}
\mu^{\ast'}\dfrac{\partial h_{i_1,t}^{\delta_{0,i_1}/2}(\vartheta_0)}{\partial \theta} &= \sum\limits_{i_2=1}^d \left[\mu_{i_1+i_2d}^\ast(\varepsilon_{i_2,t-1}^+)^{\delta_{0,i_2}} + \mu_{i_1+(i_2+q)d^2}^\ast(\varepsilon_{i_2,t-1}^-)^{\delta_{0,i_2}}\right] + R_{t-2}\nonumber\\
&= \sum\limits_{i_2=1}^d \left[\mu_{i_1+i_2d}^\ast \left(\sum\limits_{j_1=1}^d{\cal H}_{0,t-1}^{1/2}(i_2,j_1)\eta_{j_1,t-1}^+\right)^{\delta_{0,i_2}} \right.\nonumber\\
&\qquad\qquad \left.  + \mu_{i_1 + (i_2+q)d^2}^\ast\left(\sum\limits_{j_1=1}^d {\cal H}_{0,t-1}^{1/2}(i_2,j_1)\eta_{j_1,t-1}^-\right)^{\delta_{0,i_2}}\right]+R_{t-2},
\end{align}
%}
%\normalsize
where the vector $\mu^\ast=(\mu^\ast_{1},\dots,\mu^\ast_{s_1})'$.

Now using \eqref{truc2} the third term of the equation \eqref{lambdaVIncbis} can be rewritten as 
\begin{align}\label{part3}
\sum\limits_{i=1}^d&\mu_{i+s_2}^\ast  \dfrac{\partial h_{i_1,t}^{\delta_{i_1}/2}(\vartheta_0)}{\partial \delta_i} = \sum\limits_{i_2 = 1}^d \mu_{i_2 + s_2}^\ast \left[A_{01}^+(i_1,i_2)\log\left(\sum\limits_{j_1=1}^d \mathcal{H}_{0,t-1}^{1/2}(i_2,j_1)\eta_{j_1,t-1}^+\right)\left(\sum\limits_{j_1=1}^d\mathcal{H}_{0,t-1}^{1/2}(i_2,j_1)\eta_{j_1,t-1}^+\right)^{\delta_{0,i_2}}\right.\nonumber\\
&\qquad \left. A_1^-(i_1,i_2)\log\left(\sum\limits_{j_1=1}^d \mathcal{H}_{0,t-1}^{1/2}(i_2,j_1)\eta_{j_1,t-1}^-\right)\left(\sum\limits_{j_1=1}^d\mathcal{H}_{0,t-1}^{1/2}(i_2,j_1)\eta_{j_1,t-1}^-\right)^{\delta_{0,i_2}}\right] + R_{t-2}.
\end{align}
Finally by using \eqref{truc3}, the last term of \eqref{lambdaVIncbis} can be rewritten as
\begin{align}\label{part4}
\mu_{i_1+s_2}^\ast h_{i_1,t}^{\delta_{i_1}/2}(\vartheta_0)\log\left(h_{i_1,t}^{\delta_{i_1}/2}(\vartheta_0)\right) & = \mu_{i_1 + s_2}^\ast\left[\sum\limits_{i_2=1}^d A_{01}^+(i_1,i_2)\left(\sum\limits_{j_1=1}^d \mathcal{H}_{0,t-1}^{1/2}(i_2,j_1)\eta_{j_1,t-1}^+\right)^{\delta_{0,i_2}} \right.\nonumber\\
&\qquad \left.+ \sum\limits_{i_2=1}^d A_{01}^-(i_1,i_2)\left(\sum\limits_{j_1=1}^d \mathcal{H}_{0,t-1}^{1/2}(i_2,j_1)\eta_{j_1,t-1}^-\right)^{\delta_{0,i_2}}  + R_{t-2} \right]\nonumber \\
&\qquad \times \log\left[R_{t-2} + \sum\limits_{i_2=1}^d A_{01}^+(i_1,i_2)\left(\sum\limits_{j_1=1}^d \mathcal{H}_{0,t-1}^{1/2}(i_2,j_1)\eta_{j_1,t-1}^+\right)^{\delta_{0,i_2}} \right.\nonumber\\
&\qquad \left.+ \sum\limits_{i_2=1}^d A_{01}^-(i_1,i_2)\left(\sum\limits_{j_1=1}^d \mathcal{H}_{0,t-1}^{1/2}(i_2,j_1)\eta_{j_1,t-1}^-\right)^{\delta_{0,i_2}} \right].
\end{align}
Combining Equations \eqref{part1}, \eqref{part2}, \eqref{part3} and \eqref{part4} and by separating the non negative terms and the non positive terms, Equation \eqref{lambdaVIncbis} is equivalent to the two equations  
\begin{align}\label{EqPosInc}
&\left[\sum\limits_{i_2=1}^d A_{01}^+(i_1,i_2)\left(\sum\limits_{j_1=1}^d \mathcal{H}_{0,t-1}^{1/2}(i_2,j_1)\eta_{j_1,t-1}^+\right)^{\delta_{0,i_2}}\right]\left[\lambda_{1}\sum\limits_{i_2=1}^d\left(\eta_{i_2,t-1}^+\right)^2 + R_{t-2}\right] \nonumber\\
&\quad + R_{t-2}\sum\limits_{i_2=1}^d\left(\eta_{i_2,t-1}^+\right)^2 + R_{t-2} + \sum\limits_{i_2=1}^d \mu_{i_1 + (i_2+q)d}^\ast \left(\sum\limits_{j_1=1}^d \mathcal{H}_{0,t-1}^{1/2}(i_2,j_1)\eta_{j_1,t-1}^+\right)^{\delta_{0,i_2}}\nonumber\\
&\quad + \sum\limits_{i_2=1}^d \mu_{i_2 + s_2}^\ast A_{01}^+(i_1,i_2)\log\left(\sum\limits_{j_1=1}^d\mathcal{H}_{0,t-1}^{1/2}(i_2,j_1)\eta_{j_1,t-1}^+\right)\left(\sum\limits_{j_1=1}^d \mathcal{H}_{0,t-1}^{1/2}(i_2,j_1)\eta_{j_1,t-1}^+\right)^{\delta_{0,i_2}}\nonumber\\
&\quad - \mu_{i_1+s_2}^\ast \left[\sum\limits_{i_2=1}^d A_{01}^+(i_1,i_2)\left(\sum\limits_{j_1=1}^d\mathcal{H}_{0,t-1}^{1/2}(i_2,j_1)\eta_{j_1,t-1}^+\right)^{\delta_{0,i_2}} + R_{t-2}\right]\nonumber\\
&\quad \times \log\left[R_{t-2} + \sum\limits_{i_2=1}^dA_{01}^+(i_1,i_2) \left(\sum\limits_{j_1=1}^d\mathcal{H}_{0,t-1}^{1/2}(i_2,j_1)\eta_{j_1,t-1}^+\right)^{\delta_{0,i_2}}\right] = 0,\quad \mbox{a.s.}
\end{align} 
\begin{align}\label{EqNegInc}
&\left[\sum\limits_{i_2=1}^d A_{01}^-(i_1,i_2)\left(\sum\limits_{j_1=1}^d \mathcal{H}_{0,t-1}^{1/2}(i_2,j_1)\eta_{j_1,t-1}^-\right)^{\delta_{0,i_2}}\right]\left[\lambda_{1}\sum\limits_{i_2=1}^d\left(\eta_{i_2,t-1}^-\right)^2 + R_{t-2}\right] \nonumber\\
&\quad + R_{t-2}\sum\limits_{i_2=1}^d\left(\eta_{i_2,t-1}^-\right)^2 + R_{t-2} + \sum\limits_{i_2=1}^d \mu_{i_1 + (i_2+q)d^2}^\ast \left(\sum\limits_{j_1=1}^d \mathcal{H}_{0,t-1}^{1/2}(i_2,j_1)\eta_{j_1,t-1}^-\right)^{\delta_{0,i_2}}\nonumber\\
&\quad + \sum\limits_{i_2=1}^d \mu_{i_2 + s_2}^\ast A_{01}^-(i_1,i_2)\log\left(\sum\limits_{j_1=1}^d\mathcal{H}_{0,t-1}^{1/2}(i_2,j_1)\eta_{j_1,t-1}^-\right)\left(\sum\limits_{j_1=1}^d \mathcal{H}_{0,t-1}^{1/2}(i_2,j_1)\eta_{j_1,t-1}^-\right)^{\delta_{0,i_2}}\nonumber\\
&\quad - \mu_{i_1+s_2}^\ast \left[\sum\limits_{i_2=1}^d A_{01}^-(i_1,i_2)\left(\sum\limits_{j_1=1}^d\mathcal{H}_{0,t-1}^{1/2}(i_2,j_1)\eta_{j_1,t-1}^-\right)^{\delta_{0,i_2}} + R_{t-2}\right]\nonumber\\
&\quad \times \log\left[R_{t-2} + \sum\limits_{i_2=1}^dA_{01}^-(i_1,i_2) \left(\sum\limits_{j_1=1}^d\mathcal{H}_{0,t-1}^{1/2}(i_2,j_1)\eta_{j_1,t-1}^-\right)^{\delta_{i_2}}\right] = 0,\quad \mbox{a.s.}
\end{align} 
When $d=1$, from \eqref{EqPosInc} and \eqref{EqNegInc} we retrieve an equation of the following form obtained by \cite{Nous-PortTestUni}
\begin{equation*}
a\vert y \vert^{\underline{\delta} + 2} + [b + c(\vert y \vert^{\underline{\delta}})]\log[b + c(\vert y \vert^{\underline{\delta}})]+[d + e\log(\vert y \vert)] \vert y \vert^{\underline{\delta}} + f y^2 + g = 0
\end{equation*}
which cannot have more than 11 positive roots or more than 11 negative roots, except if $a=b=c=d=e=f=g=0$. 

When $d\geq2$ and also from \eqref{EqPosInc} and \eqref{EqNegInc}, for a fixed component, we obtain an equation of the form 
\begin{align*}
&\sum\limits_{i=1}^d a_i\vert y\vert^{\delta_i + 2} + \sum\limits_{i=1}^d b_i \vert y\vert^{\delta_i+1} + \sum\limits_{i=1}^d c_i \vert y\vert^{\delta_i}+ \sum\limits_{i=1}^dd_i\log(\vert y\vert)\vert y\vert^{\delta_i}\\ 
&\qquad  + \left(e + \sum\limits_{i=1}^d e_i \vert y\vert^{\delta_i}\right)\log\left(f + \sum\limits_{i=1}^df_i\vert y\vert^{\delta_i}\right) + g y^2 + h\vert y\vert + k = 0.
\end{align*}
Note that an equation of this form can not have more than $11d+1$ non negative roots or more than $11d+1$ non positive roots for $d\geq 2$, unless $a_i = b_i = c_i = d_i = e_i = f_i = e = f = g = h = k = 0$.\\
By the assumption $\mathbf{A9}'$, Equations \eqref{EqPosInc} and \eqref{EqNegInc} imply that \\ 
$ \lambda_{1}\left[\sum_{i_2=1}^d A_{01}^+(i_1,i_2) + A_{01}^-(i_1,i_2)\right] = 0$ and $\mu_{i+s_2}^\ast\left[\sum_{i_2=1}^d A_{01}^+(i_1,i_2) + A_{01}^-(i_1,i_2)\right] = 0$ for all $i=1,\ldots, d$. But under the assumption $\mathbf{A4}$, if $p>0$, $\mathcal{A}_0(1)^+ + \mathcal{A}_0^- \neq 0$. It is impossible to have $A_{01}^+(i_1,i) = A_{01}^+(i_1,i) = 0$, for all $i=1,\ldots,d$. Then, there exists an $i_0$ such that $A_{01}(i_1,i_0)^+ + A_{01}(i_1,i_0)^- \neq 0$ and we then have $\lambda_{1} = 0$ and $\mu_{i_0 + s_2}^\ast = 0$.

In the general case, Equation \eqref{mud} necessarily leads
\begin{equation*}
A_{01}^+(i_1,i_0) + A_{01}^-(i_1,i_0) = \cdots = A_{0q}^+(i_1,i_0) + A_{0q}^-(i_1,i_0) = 0, \qquad \forall i_0,i_1 = 1,\ldots, d,
\end{equation*}
that is impossible under Assumption $\mathbf{A4}$ and then $\lambda = 0$. This is in contradiction with $\lambda'V = 0$, almost surely, that leads that the assumption of non invertibility of matrix ${\cal D}$ is absurd.
\zak
\subsection{Proof of Theorem \ref{thm-Port-test-inconnu}}
The proof is the same to that  of Theorem \ref{thm-Port-test-connu} in the case where the power 
$\underline{\delta}_0$ is assumed to be known.  
\zak
\subsection{Proof of Corollary \ref{thm-Port-test-inconnubis}}
The proof is the same to that  of Corollary \ref{thm-Port-test-connubis}  in the case where the power $\underline{\delta}_0$ is assumed to be known.  
\zak

\begin{center}
{\bf Acknowledgements}
\end{center}
We sincerely thank the
anonymous reviewers and editor for helpful remarks.
\begin{center}
{\bf Data availability statement}
\end{center}
The data used in Section \ref{donnees reelles} are available in the supporting information of this article. It is openly available from %the website Yahoo Finance: http://fr.finance.yahoo.com/.
the website of the European Central Bank: http://www.ecb.int/stats/exchange/eurofxref/html/index.en.html.

\begin{center}
{\bf Supporting information}
\end{center}
Additional Supporting Information may be found online in the supporting information tab for this article. 
%%%% Debut Biblio %%%%%%

\end{document}